\documentclass[a4paper,12pt]{article}
\usepackage{amsmath,amssymb,amsthm,graphicx}
\usepackage{enumitem}
\usepackage{color}
\usepackage{epsfig}
\usepackage{graphics}
\usepackage[font=small]{caption}
\usepackage{subfig}
\usepackage[hang,flushmargin]{footmisc} 
\usepackage{float}
\usepackage{booktabs}
\usepackage{natbib}
\usepackage{setspace}
\usepackage{mathrsfs}
\usepackage{ulem}
\usepackage[left=3cm,right=3cm,bottom=3cm,top=3cm]{geometry}

\newcommand{\reals}{\mathbb{R}}

\newcommand{\pr}{\mathbb{P}}        
\newcommand{\ex}{\mathbb{E}}        
\newcommand{\var}{\textnormal{Var}} 

\newcommand{\normal}{N}        

\newcommand{\ind}{\boldsymbol{1}} 

\theoremstyle{plain}

\newtheorem{theorem}{Theorem}[section]

\newtheorem{lemma}{Lemma}

\begin{document}

\begin{center}
{\LARGE \textbf{Clustering with}} \\[0.35cm]
{\LARGE \textbf{Statistical Error Control}}
\end{center}
\vspace{-15pt}

\begin{center}
\begin{minipage}[c][2cm][c]{5cm}
\centering {\large Michael Vogt\renewcommand{\thefootnote}{1}\footnotemark[1]} \\[0.15cm] {University of Bonn}      
\end{minipage}
\begin{minipage}[c][2cm][c]{5cm}
\centering {\large Matthias Schmid\renewcommand{\thefootnote}{2}\footnotemark[2]} \\[0.15cm] {University of Bonn}      
\end{minipage}
\end{center}
\vspace{-25pt}

\footnotetext[1]{Corresponding author. Address: Department of Economics and Hausdorff Center for Mathema\-tics, University of Bonn, 53113 Bonn, Germany. Email: \texttt{michael.vogt@uni-bonn.de}.}
\renewcommand{\thefootnote}{2}
\footnotetext[2]{Address: Department of Medical Biometry, Informatics and Epidemiology, University of Bonn, 53105 Bonn, Germany. \texttt{matthias.schmid@imbie.uni-bonn.de}.}

\renewcommand{\thefootnote}{\arabic{footnote}}
\setcounter{footnote}{0}

\renewcommand{\abstractname}{}
\begin{abstract}
\noindent This paper presents a clustering approach that allows for rigorous statistical error control similar to a statistical test. We develop estimators for both the unknown number of clusters and the clusters themselves. The estimators depend on a tuning parameter $\alpha$ which is similar to the significance level of a statistical hypothesis test. By choosing $\alpha$, one can control the probability of overestimating the true number of clusters, while the probability of underestimation is asymptotically negligible. In addition, the probability that the estimated clusters differ from the true ones is controlled. In the theoretical part of the paper, formal versions of these statements on statistical error control are derived in a standard model setting with convex clusters. A simulation study and two applications to temperature and gene expression microarray data complement the theoretical analysis. 
\end{abstract}
\vspace{5pt}

\renewcommand{\baselinestretch}{1.2}\normalsize

\noindent \textbf{Key words:} Cluster analysis; number of clusters; multiple statistical testing; statistical error control; $k$-means clustering. 

\noindent \textbf{AMS 2010 subject classifications:} 62H30; 62H15; 62E20. 

\numberwithin{equation}{section}
\allowdisplaybreaks[2]

\vspace{-5pt}

\section{Introduction}\label{sec-intro}

In a wide range of applications, the aim is to cluster a large number of subjects into a small number of groups. Prominent examples are the clustering of genes in microarray analysis \citep{jiang}, the clustering of temperature curves using data recorded on a spatial grid \citep{Fovell1993, DeGaetano2001}, and the clustering of consumer profiles on the basis of survey data \citep{Wedel2000}.

A major challenge in cluster analysis is to estimate the unknown number of groups $K_0$ from a sample of data. A common approach is to compute a criterion function which measures the quality of the clustering for different cluster numbers $K$. An estimator of $K_0$ is then obtained by optimizing the criterion function over $K$. Prominent examples of this approach are the Hartigan index  \citep{hartigan1975}, the silhouette statistic \citep{rousseeuw1987} and the gap statistic \citep{gapstatistic}.

Another common way to estimate $K_0$ is based on statistical test theory. Roughly speaking, one can distinguish between two types of test-based procedures: The {\it first type} relies on a statistical test which either checks whether some clusters can be merged or whether a cluster can be subdivided. Given a set of clusters, the test is repeatedly applied until no clusters can be merged or split any more. The number of remaining clusters serves as an estimator of $K_0$. Classical examples of methods that proceed in this way are discussed in \citet[Chapter 3.5]{Gordon1999} who terms them ``local methods''. Obviously, these methods involve a multiple testing problem. However, the employed critical values do not properly control for the fact that multiple tests are performed. The significance level $\alpha$ used to carry out the tests thus cannot be interpreted strictly. Put differently, the procedures do not allow for rigorous statistical error control.

Test-based approaches of the {\it second type} proceed by sequentially testing a model with $K$ clusters against one with $K+1$ clusters. The smallest number $K$ for which the test does not reject serves as an estimator of $K_0$. Most work in this direction has been done in the framework of Gaussian mixture models; see \cite{McLachlan2014} for an overview. However, deriving a general theory for testing a mixture with $K$ components against one with $K' > K$ components has turned out to be a very challenging problem; see \cite{GhoshSen1985} and \cite{Hartigan1985} for a description of the main technical issues involved. Many results are therefore restricted to the special case of testing a homogeneous model against a mixture with $K=2$ clusters; see \cite{LiuShao2004} and \cite{LiChenMarriott2009} among many others. More general test procedures often lack a complete theoretical foundation or are based on very restrictive conditions.

Only recently, there have been some advances in developing a general theory for testing $K$ against $K' > K$ clusters under reasonably weak conditions. In a mixture model setup, \cite{LiChen2010} and \cite{ChenLiFu2012} have constructed a new expectation-maximization (EM) procedure to approach this testing problem. Outside the mixture model context, \cite{Maitra2012} have developed a bootstrap procedure to test a model with $K$ groups against one with $K^\prime > K$ groups. These papers derive the theoretical properties of the proposed tests under the null hypothesis of $K$ clusters, where $K$ is a pre-specified fixed number. However, they do not formally investigate the properties of a procedure which estimates $K_0$ by sequentially applying the tests. In particular, they do not analyze whether such a sequential procedure may allow for a rigorous interpretation of the significance level $\alpha$ that is used to carry out the tests.

The main contribution of this paper is to construct an estimator $\widehat{K}_0$ of $K_0$ which allows for rigorous statistical error control in the following sense: For any pre-specified significance level $\alpha \in (0,1)$, the proposed estimator $\widehat{K}_0 = \widehat{K}_0(\alpha)$ has the property that 
\begin{align}
\pr \big( \widehat{K}_0 > K_0 \big) & = \alpha + o(1) \label{err-control-1a}, \\
\pr \big( \widehat{K}_0 < K_0 \big) & = o(1). \label{err-control-1b} 
\end{align}
According to this, the probability of overestimating $K_0$ is controlled by the level $\alpha$, while the probability of underestimating $K_0$ is asymptotically negligible. By picking $\alpha$, we can thus control the probability of choosing too many clusters, while, on the other hand, we can ignore the probability of choosing too few clusters (at least asymptotically).

We show how to construct an estimator $\widehat{K}_0$ with the properties \eqref{err-control-1a} and \eqref{err-control-1b} in a standard model setting with convex clusters which is introduced in Section \ref{sec-model}. Our estimation approach is developed in Section \ref{sec-est}. As we will see, the proposed procedure does not only provide us with an estimator of $K_0$. It also yields estimators of the groups themselves which allow for statistical error control similarly to $\widehat{K}_0$. Our approach is based on the following general strategy: 
\begin{enumerate}[label = (\roman*), leftmargin=0.85cm]
\item Construct a statistical test which, for any given number $K$, checks the null hypothesis that there are $K$ clusters in the data. 
\item Starting with $K = 1$, sequentially apply this test until it does not reject the null hypothesis of $K$ clusters any more. 
\item Define the estimator $\widehat{K}_0$ of $K_0$ as the smallest number $K$ for which the test does not reject the null.
\end{enumerate}
This strategy is discussed in detail in Section \ref{subsec-est-approach}. It is generic in the sense that it can be employed with different test statistics. For our theoretical analysis, we apply it with a specific statistic which is introduced in Section \ref{subsec-est-stat}. For this specific choice, we derive the statements \eqref{err-control-1a} and \eqref{err-control-1b} on statistical error control under suitable regularity conditions. Some alternative choices of the test statistic are discussed in Section \ref{sec-ext}. In the following, we refer to our estimation procedure as {\it CluStErr} (``{\it Clu}stering with {\it St}atistical {\it Err}or Control'').

The theoretical properties of our estimators, in particular the statements \eqref{err-control-1a} and \eqref{err-control-1b}, are derived in Section \ref{sec-asym}. As we will see there, our theory is valid under quite general conditions. First of all, as opposed to many other studies from the clustering literature including those from a Gaussian mixture context, we do not restrict the random variables in our model to be Gaussian. For our theory to work, we merely require them to satisfy a set of moment conditions. Secondly, our approach is essentially free of tuning parameters, the only choice parameter being the significance level $\alpha$. Thirdly, to apply our method, we of course need to compute critical values for the underlying test. However, as opposed to other test-based methods, we do not have to estimate or bootstrap the critical values by a complicated procedure. They can rather be easily computed analytically. This makes our method particularly simple to implement in practice.

We complement the theoretical analysis of the paper by a simulation study and two applications on temperature and microarray data in Section \ref{sec-app}. The R code to reproduce the numerical examples is contained in the add-on package {\bf CluStErr} \citep{CluStErrPackage}, which implements the CluStErr method and which is part of the supplemental materials of the paper.

\section{Model}\label{sec-model}

Suppose we measure $p$ features on $n$ different subjects. In particular, for each subject $i \in \{1,\ldots,n\}$, we observe the vector $\boldsymbol{Y}_i = (Y_{i1}\,\ldots,Y_{ip})^\top$, where $Y_{ij}$ denotes the measurement of the $j$-th feature for the $i$-th subject. Our data sample thus has the form $\{ \boldsymbol{Y}_i: 1 \le i \le n \}$. 
Both the number of subjects $n$ and the number of features $p$ are assumed to tend to infinity, with $n$ diverging much faster than $p$. This reflects the fact that $n$ is much larger than $p$ in the applications we have in mind. When clustering the genes in a typical microarray data set, for instance, the number of genes $n$ is usually a few thousands, whereas the number of tissue samples $p$ is not more than a few tenths. The exact technical conditions on the sizes of $n$ and $p$ are laid out in Section \ref{subsec-asym-ass}.

The data vectors $\boldsymbol{Y}_i$ of the various subjects $i=1,\ldots,n$ are supposed to satisfy the model
\begin{equation}\label{model-eq1}
\boldsymbol{Y}_i = \boldsymbol{\mu}_i + \boldsymbol{e}_i,
\end{equation}
where $\boldsymbol{\mu}_i = (\mu_{i1},\ldots,\mu_{ip})^\top$ is a deterministic signal vector and $\boldsymbol{e}_i = (e_{i1},\ldots,e_{ip})^\top$ is the noise vector. The subjects in our sample are assumed to belong to $K_0$ different classes. More specifically, the set of subjects $\{1,\ldots,n\}$ can be partitioned into $K_0$ groups $G_1,\ldots,G_{K_0}$ such that for each $k = 1,\ldots,K_0$, 
\begin{equation}\label{model-eq2}
\boldsymbol{\mu}_i = \boldsymbol{m}_k \quad \text{ for all } i \in G_k,   
\end{equation}
where $\boldsymbol{m}_k \in \reals^p$ are vectors with $\boldsymbol{m}_k \ne \boldsymbol{m}_{k^\prime}$ for $k \ne k^\prime$. Hence, the members of each group $G_k$ all have the same signal vector $\boldsymbol{m}_k$.

Equations \eqref{model-eq1} and \eqref{model-eq2} specify a model with convex spherical clusters which underlies the $k$-means and many other Euclidean distance-based clustering algorithms. This framework has been employed extensively in the literature and is useful in a wide range of applications, which is also illustrated by the examples in Section \ref{sec-app}. It is thus a suitable baseline model for developing our ideas on clustering with statistical error control. We now discuss the two model equations \eqref{model-eq1} and \eqref{model-eq2} in detail.

\vspace{8pt}

\noindent \textbf{Details on equation (\ref{model-eq1}).} 
The noise vector $\boldsymbol{e}_i = (e_{i1},\ldots,e_{ip})^\top$ is assumed to consist of entries $e_{ij}$ with the additive component structure $e_{ij} = \alpha_i + \varepsilon_{ij}$. Equation \eqref{model-eq1} for the $i$-th subject thus writes as
\begin{equation}\label{model-eq1a} 
\boldsymbol{Y}_i = \boldsymbol{\mu}_i + \boldsymbol{\alpha}_i + \boldsymbol{\varepsilon}_i,
\end{equation}
where $\boldsymbol{\alpha}_i = (\alpha_i,\ldots,\alpha_i)^\top$ and $\boldsymbol{\varepsilon}_i = (\varepsilon_{i1},\ldots,\varepsilon_{ip})^\top$. 
Here, $\alpha_i$ is a subject-specific random intercept term.  Moreover, the terms $\varepsilon_{ij}$ are standard idiosyncratic noise variables with $\ex[\varepsilon_{ij}] = 0$. We assume the error terms $\varepsilon_{ij}$ to be i.i.d.\ both across $i$ and $j$. The random intercepts $\alpha_i$, in contrast, are allowed to be dependent across subjects $i$ in an arbitrary way. 

In general, the components of \eqref{model-eq1a} may depend on the sample size $p$. The exact formulation of the model equation \eqref{model-eq1a} for the $i$-th subject thus reads
$\boldsymbol{Y}_{i,p} = \boldsymbol{\mu}_{i,p} + \boldsymbol{\alpha}_{i,p} + \boldsymbol{\varepsilon}_{i,p}$,
where $\boldsymbol{Y}_{i,p} = (Y_{i1,p},\ldots,Y_{ip,p})^\top$, $\boldsymbol{\mu}_{i,p} = (\mu_{i1,p},\ldots,\mu_{ip,p})^\top$, $\boldsymbol{\alpha}_{i,p} = (\alpha_{i,p},\ldots,\alpha_{i,p})^\top$ and $\boldsymbol{\varepsilon}_{i,p} = (\varepsilon_{i1,p},\ldots,\varepsilon_{ip,p})^\top$. However, to keep the notation simple, we suppress this dependence on $p$ and write the model for the $i$-th subject as \eqref{model-eq1a}.

If we drop the random intercept $\boldsymbol{\alpha}_i$ from \eqref{model-eq1a}, the signal vector $\boldsymbol{\mu}_i$ is equal to the mean $\ex[\boldsymbol{Y}_i]$. In the general equation \eqref{model-eq1a} in contrast, $\boldsymbol{\mu}_i$ is only identified up to an additive constant. 
To identify $\boldsymbol{\mu}_i$ in \eqref{model-eq1a}, we impose the normalization constraint $p^{-1} \sum\nolimits_{j=1}^p \mu_{ij} = 0$ for each $i$. We thus normalize the entries of $\boldsymbol{\mu}_i$ to be zero on average for each $i$. Under the technical conditions specified in Section \ref{subsec-asym-ass}, the constraint $p^{-1} \sum\nolimits_{j=1}^p \mu_{ij} = 0$ implies that $\alpha_i = \lim_{p \rightarrow \infty} p^{-1} \sum\nolimits_{j=1}^p Y_{ij}$ almost surely, which in turn identifies the signal vector $\boldsymbol{\mu}_i$. 

\vspace{8pt}

\noindent \textbf{Details on equation (\ref{model-eq2}).} This equation specifies the group structure in our model. 
We assume the number of groups $K_0$ to be fixed, implying that the groups $G_k = G_{k,n}$ depend on the sample size $n$. Keeping the number of classes $K_0$ fixed while letting the size of the classes $G_{k,n}$ grow is a reasonable assumption: It reflects the fact that in most applications, we expect the number of groups $K_0$ to be very small as compared to the total number of subjects $n$. To keep the notation simple, we suppress the dependence of the classes $G_{k,n}$ on the sample size $n$ and denote them by $G_k$ throughout the paper. 

\vspace{8pt}

\noindent In the remainder of this section, we discuss two special cases of model \eqref{model-eq1}--\eqref{model-eq2} which are relevant for our applications in Section \ref{sec-app}. 

\newpage

\noindent \textbf{A model for the clustering of time series data.} Suppose we observe time series $\boldsymbol{Y}_i = (Y_{i1},\ldots,Y_{ip})^\top$ of length $p$ for $n$ different subjects $i$. The time series $\boldsymbol{Y}_i$ of the $i$-th subject is assumed to follow the time trend model
\begin{equation}\label{trend-model}
Y_{ij} = \mu_i(t_j) + \alpha_i + \varepsilon_{ij} \quad (1 \le j \le p), 
\end{equation}
where $\mu_i(\cdot)$ is an unknown nonparametric trend function and $t_1 < \ldots < t_p$ are the observed time points. The deterministic design points $t_j$ are supposed to be the same across subjects $i$ and are normalized to lie in the unit interval. An important example is the equidistant design $t_j = j/p$. However, it is also possible to allow for non-equidistant designs. To identify the trend function $\mu_i(\cdot)$ in \eqref{trend-model}, we suppose that $\int_0^1 \mu_i(w) dw = 0$ for each $i$, which is a slight modification of the identification constraint stipulated in \eqref{model-eq1a}. Analogous to our general model, we impose a group structure on the observed time series: There are $K_0$ groups of time series $G_1,\ldots,G_{K_0}$ such that $\mu_i(\cdot) = m_k(\cdot)$ for all $i \in G_k$. Hence, the members of each class $G_k$ all have the same time trend function $m_k(\cdot)$.

\vspace{8pt}

\noindent \textbf{A model for the clustering of genes in microarray experiments.} In a microarray experiment, the expression levels of $n$ different genes are often measured in $p$ different tissue samples (obtained, e.g., from $p$ different patients). For each gene $i$, we observe the vector $\boldsymbol{Y}_i = (Y_{i1},\ldots,Y_{ip})^\top$, where $Y_{ij}$ is the measured expression level of gene $i$ for tissue sample $j$. The vector $\boldsymbol{Y}_i$ of gene $i$ is supposed to satisfy the model equation \eqref{model-eq1a}, which componentwise reads as 
\begin{equation}\label{gene-model-1a}
Y_{ij} = \mu_{ij} + \alpha_i + \varepsilon_{ij} \quad (1 \le j \le p). 
\end{equation}
Here, $\mu_{ij}$ can be regarded as the true expression level of gene $i$ for tissue $j$, whereas $Y_{ij}$ is the measured expression level corrupted by the noise term $\alpha_i + \varepsilon_{ij}$. 

Most microarray experiments involve different types of tissues, for example tumor ``cases'' versus healthy ``controls'', or different tumor (sub)types. We therefore suppose that there are $T$ different types of tissues in our sample and order them according to their type (which is known by experimental design). More specifically, the tissues $j$ of type $t$ are labelled by $j_{t-1} \le j < j_t$, where $1 = j_0 < j_1 < \ldots < j_{T-1} < j_T = p+1$. If the patients from which tissues are obtained constitute samples of sufficiently homogeneous populations, it is natural to assume that the true expression level $\mu_{ij}$ of gene $i$ is the same for tissues $j$ of the same type, i.e., $\mu_{ij} = \mu_{ij^\prime}$ for $j_{t-1} \le j,j^\prime < j_t$. The signal vector $\boldsymbol{\mu}_i$ thus has a piecewise constant structure for each $i$; see Figures \ref{fig:IlluGene} and \ref{fig:ClustersGene} in Section \ref{sec-app} for an illustration. 

As in our general model, we suppose that there are $K_0$ groups of genes $G_1,\ldots,G_{K_0}$ such that $\boldsymbol{\mu}_i = \boldsymbol{m}_k$ for all $i \in G_k$ and some vector $\boldsymbol{m}_k$. The genes of each class $G_k$ thus have the same (co-)expression profile $\boldsymbol{m}_k$.

\section{Estimation Method}\label{sec-est}

We now present our approach to estimate the unknown groups $G_1,\ldots,G_{K_0}$ and their unknown number $K_0$ in model \eqref{model-eq1}--\eqref{model-eq2}. Section \ref{subsec-est-approach} gives an overview of the general method, while Sections \ref{subsec-est-stat}--\ref{subsec-est-sigma} fill in the details.

\subsection{The general method}\label{subsec-est-approach}

To construct our method, we proceed in two steps: In the first step, we specify an algorithm that clusters the set of subjects $\{1,\ldots,n\}$ into $K$ groups for any given number $K$ (which may or may not coincide with the true number of classes $K_0$). Let $\{ \widehat{G}_k^{[K]}: 1 \le k \le K \}$ be the $K$ clusters produced by the algorithm when the number of clusters is $K$. For $K = 1$, we trivially set $\widehat{G}_1^{[1]} = \{1,\ldots,n\}$. For our theory to work, we require the clustering algorithm to consistently estimate the class structure $\{G_k:1 \le k \le K_0\}$ when $K = K_0$. More specifically, we require the estimators $\{ \widehat{G}_k^{[K_0]}: 1 \le k \le K_0\}$ to have the property that 
\begin{equation}\label{prop-clusters}
\pr \Big( \big\{ \widehat{G}_k^{[K_0]}: 1 \le k \le K_0 \big\} = \big\{ G_k: 1 \le k \le K_0 \big\} \Big) \rightarrow 1. 
\end{equation}
This is a quite weak restriction which is satisfied by a wide range of clustering algorithms under our regularity conditions. As shown in Section \ref{subsec-est-kmeans}, it is for example satisfied by a $k$-means type algorithm. Moreover, it can be shown to hold for a number of hierarchical clustering algorithms, in particular for agglomerative algorithms with single, average and complete linkage. Our estimation method can be based on any clustering algorithm that has the consistency property \eqref{prop-clusters}.

In the second step, we construct a test for each $K$ which checks whether the data can be well described by the $K$ clusters $\widehat{G}_1^{[K]},\ldots,\widehat{G}_K^{[K]}$. We thereby test whether the number of clusters is equal to $K$. More formally, we use the $K$-cluster partition $\{ \widehat{G}_k^{[K]}: 1 \le k \le K \}$ to construct a statistic $\widehat{\mathcal{H}}^{[K]}$ that allows us to test the hypothesis $H_0: K = K_0$ versus $H_1: K < K_0$. For any given number of clusters $K$, our test is defined as $T_{\alpha}^{[K]} = \ind(\widehat{\mathcal{H}}^{[K]} > q(\alpha))$, where $q(\alpha)$ is the $(1-\alpha)$-quantile of a known distribution which will be specified later on. We reject $H_0$ at the level $\alpha$ if $T_{\alpha}^{[K]} = 1$, i.e., if $\widehat{\mathcal{H}}^{[K]} > q(\alpha)$. A detailed construction of the statistic $\widehat{\mathcal{H}}^{[K]}$ along with a precise definition of the quantile $q(\alpha)$ is given in Section \ref{subsec-est-stat}.


To estimate the classes $G_1,\ldots,G_{K_0}$ and their number $K_0$, we proceed as follows: For each $K=1,2,\ldots$, we check whether $\widehat{\mathcal{H}}^{[K]} \le q(\alpha)$ and stop as soon as this criterion is satisfied. Put differently, we carry out our test for each $K =1,2,\ldots$ until it does not reject $H_0$ any more. Our estimator of $K_0$ is defined as the smallest number $K$ for which $\widehat{\mathcal{H}}^{[K]} \le q(\alpha)$, that is, for which the test does not reject $H_0$. Formally speaking, we define 
\begin{equation}\label{def-K0-hat}
\widehat{K}_0 = \min \big\{ K =1,2,\ldots \, \big| \, \widehat{\mathcal{H}}^{[K]} \le q(\alpha) \big\}. 
\end{equation}
Moreover, we estimate the class structure $\{ G_k : 1 \le k \le K_0 \}$ by the partition $\{ \widehat{G}_k: 1 \le k \le \widehat{K}_0 \}$, where we set $\widehat{G}_k = \widehat{G}_k^{[\widehat{K}_0]}$. The definition \eqref{def-K0-hat} can equivalently be written as
\begin{equation}\label{def'-K0-hat}
\widehat{K}_0 = \min \big\{ K =1,2,\ldots \, \big| \, \widehat{p}^{[K]} > \alpha \big\}, 
\end{equation}  
where $\widehat{p}^{[K]}$ is the $p$-value corresponding to the statistic $\widehat{\mathcal{H}}^{[K]}$. The heuristic idea behind \eqref{def'-K0-hat} is as follows: Starting with $K=1$, we successively test whether the data can be well described by a model with $K$ clusters, in particular by the partition $\{ \widehat{G}_k^{[K]}: 1 \le k \le K\}$. For each $K$, we compute the $p$-value $\widehat{p}^{[K]}$ which expresses our confidence in a model with $K$ clusters. We stop as soon as $\widehat{p}^{[K]} > \alpha$, that is, as soon as we have enough statistical confidence in a model with $K$ groups.


As shown in Section \ref{sec-asym}, under appropriate regularity conditions, our statistic $\widehat{\mathcal{H}}^{[K]}$ has the property that 
\begin{equation}\label{prop-stat}
\pr\Big( \widehat{\mathcal{H}}^{[K]} \le q(\alpha) \Big) = 
\begin{cases}
o(1) & \text{for } K < K_0 \\
(1 - \alpha) + o(1) & \text{for } K = K_0.
\end{cases} 
\end{equation}
Put differently, $\pr(T_{\alpha}^{[K]} = 0) \rightarrow 1-\alpha$ for $K = K_0$ and $\pr(T_{\alpha}^{[K]} = 1) \rightarrow 1$ for $K < K_0$. Hence, our test is asymptotically of level $\alpha$. Moreover, it detects the alternative $H_1: K < K_0$ with probability tending to $1$, that is, its power against $H_1$ is asymptotically equal to $1$. From \eqref{prop-stat}, it follows that 
\begin{align}
\pi_>(\alpha) & := \pr \big( \widehat{K}_0 > K_0 \big) = \alpha + o(1) \label{error-control-1a} \\
\pi_<(\alpha) & := \pr \big( \widehat{K}_0 < K_0 \big) = o(1). \label{error-control-1b} 
\end{align}
Hence, the  probability of overestimating $K_0$ is asymptotically bounded by $\alpha$, while the probability of underestimating $K_0$ is asymptotically negligible. By picking $\alpha$, we can thus control the probability of choosing too many clusters similarly to the type-I-error probability of a test. Moreover, we can asymptotically ignore the probability of choosing too few clusters similarly to the type-II-error probability of a test. In finite samples, there is of course a trade-off between the probabilities of under- and overestimating $K_0$: By decreasing the significance level $\alpha$, we can reduce the probability of overestimating $K_0$, since $\alpha^\prime + o(1) = \pi_>(\alpha^\prime) \le \pi_>(\alpha) = \alpha + o(1)$ for $\alpha^\prime < \alpha$. However, we pay for this by increasing the probability of underestimating $K_0$, since $\pi_<(\alpha^\prime) \ge \pi_<(\alpha)$ for $\alpha^\prime < \alpha$. This can also be regarded as a trade-off between the size and the power of the test on which $\widehat{K}_0$ is based. Taken together, the two statements \eqref{error-control-1a} and \eqref{error-control-1b} yield that
\begin{equation}\label{error-control-1c}
\pr \big( \widehat{K}_0 \ne K_0 \big) = \alpha + o(1), 
\end{equation}
i.e., the probability that the estimated number of classes $\widehat{K}_0$ differs from the true number of classes $K_0$ is asymptotically equal to $\alpha$. With the help of \eqref{error-control-1c} and the consistency property \eqref{prop-clusters} of the estimated clusters, we can further show that 
\begin{equation}\label{error-control-2}
\pr \Big( \big\{ \widehat{G}_k: 1 \le k \le \widehat{K}_0 \big\} \ne \big\{ G_k: 1 \le k \le K_0 \big\} \Big) = \alpha + o(1), 
\end{equation}
i.e., the probability of making a classification error is asymptotically equal to $\alpha$ as well. The statements \eqref{error-control-1a}--\eqref{error-control-2} give a mathematically precise description of the statistical error control that can be performed by our method.

\subsection{Construction of the statistic $\boldsymbol{\widehat{\mathcal{H}}^{[K]}}$}\label{subsec-est-stat}

To construct the statistic $\widehat{\mathcal{H}}^{[K]}$, we use the following notation: 
\begin{enumerate}[label=(\roman*),leftmargin=0.85cm]

\item Let $Y_{ij}^* = Y_{ij} - \alpha_i$ be the observations adjusted for the random intercepts $\alpha_i$ and set $\widehat{Y}_{ij} = Y_{ij} - \overline{Y}_i$ with $\overline{Y}_i = p^{-1} \sum\nolimits_{j=1}^p Y_{ij}$. The variables $\widehat{Y}_{ij}$ serve as approximations of $Y_{ij}^*$, since under standard regularity conditions
\begin{align*} 
\widehat{Y}_{ij} 
 & = \mu_{ij} + \varepsilon_{ij} - \frac{1}{p} \sum\limits_{j=1}^p \mu_{ij} - \frac{1}{p} \sum\limits_{j=1}^p \varepsilon_{ij} \\
 & = \mu_{ij} + \varepsilon_{ij} + O_p(p^{-1/2}) =  Y_{ij}^* + O_p(p^{-1/2}).
\end{align*}

\item For any set $S \subseteq \{1,\ldots,n\}$, let $m_{j,S} = (\#S)^{-1} \sum\nolimits_{i \in S} \mu_{ij}$ be the average of the signals $\mu_{ij}$ with $i \in S$ and estimate it by $\widehat{m}_{j,S} = (\#S)^{-1} \sum\nolimits_{i \in S} \widehat{Y}_{ij}$. We use the notation $m_{j,k}^{[K]} = m_{j,\widehat{G}_k^{[K]}}$ and $\widehat{m}_{j,k}^{[K]} = \widehat{m}_{j,\widehat{G}_k^{[K]}}$ to denote the average of the signals in the cluster $\widehat{G}_k^{[K]}$ and its estimator, respectively. 

\item For any cluster $\widehat{G}_k^{[K]}$, we define cluster-specific residuals by setting $\widehat{\varepsilon}_{ij}^{[K]} = \widehat{Y}_{ij} - \widehat{m}_{j,k}^{[K]}$ for $i \in \widehat{G}_k^{[K]}$ and $1 \le j \le p$.

\item Let $\widehat{\sigma}^2$ be an estimator of the error variance $\sigma^2 = \ex[\varepsilon_{ij}^2]$. Moreover, let $\widehat{\kappa}$ be an estimator of the parameter $\kappa = (\ex[ \{ (\varepsilon_{ij}/\sigma)^2 - 1 \}^2])^{1/2}$, which serves as a normalization constant later on. See Section \ref{subsec-est-sigma} for a detailed construction of the estimators $\widehat{\sigma}^2$ and $\widehat{\kappa}$. 

\end{enumerate}
With this notation at hand, we define the statistic
\begin{equation}\label{delta-hat}
\widehat{\Delta}_i^{[K]} = \frac{1}{\sqrt{p}} \sum\limits_{j=1}^p \Big\{ \Big(\frac{\widehat{\varepsilon}_{ij}^{[K]}}{\widehat{\sigma}}\Big)^2 - 1 \Big\} \Big/ \widehat{\kappa} 
\end{equation}
for each subject $i$. This is essentially a scaled version of the residual sum of squares for the $i$-th subject when the number of clusters is $K$. Intuitively, $\widehat{\Delta}_i^{[K]}$ measures how well the data of the $i$-th subject are described when the sample of subjects is partitioned into the $K$ clusters $\widehat{G}_1^{[K]}, \ldots,\widehat{G}_K^{[K]}$. 
The individual statistics $\widehat{\Delta}_i^{[K]}$ are the building blocks of the overall statistic $\widehat{\mathcal{H}}^{[K]}$.

Before we move on with the construction of $\widehat{\mathcal{H}}^{[K]}$, we have a closer look at the stochastic behaviour of the statistics $\widehat{\Delta}_i^{[K]}$. To do so, we consider the following stylized situation: We assume that the variables $\varepsilon_{ij}$ are i.i.d.\ normally distributed with mean $0$ and variance $\sigma^2$. Moreover, we neglect the estimation error in the expressions $\widehat{Y}_{ij}$, $\widehat{m}_{j,k}^{[K]}$, $\widehat{\sigma}^2$ and $\widehat{\kappa}$. In this situation,  
\[ \widehat{\Delta}_i^{[K]} = \frac{1}{\sqrt{p}} \sum\limits_{j = 1}^p \Big\{ \frac{(\varepsilon_{ij} + d_{ij})^2}{\sigma^2} - 1 \Big\} \Big/ \kappa \]
for any $i \in S = \widehat{G}_k^{[K]}$, where $d_{ij} = \mu_{ij} - (\#S)^{-1} \sum\nolimits_{i^\prime \in S} \mu_{i^\prime j}$ is the difference between the signal $\mu_{ij}$ of the $i$-th subject and the average signal in the cluster $S$. We now give a heuristic discussion of the behaviour of $\widehat{\Delta}_i^{[K]}$ in the following two cases: 
\begin{itemize}[leftmargin=1.75cm]

\item[$K=K_0$:] By condition \eqref{prop-clusters}, $\widehat{G}_k^{[K_0]}$ consistently estimates $G_k$. Neglecting the estimation error in $\widehat{G}_k^{[K_0]}$, we obtain that 
\[ \widehat{\Delta}_i^{[K_0]} = \frac{1}{\sqrt{p}} \sum\limits_{j = 1}^p \Big\{ \frac{\varepsilon_{ij}^2}{\sigma^2} - 1 \Big\} \Big/ \kappa. \]
Since $\varepsilon_{ij}/\sigma$ is standard normal, $\kappa = \sqrt{2}$ and thus
\begin{equation}\label{delta-K0}
\widehat{\Delta}_i^{[K_0]} \sim \frac{\chi_p^2 - p}{\sqrt{2p}} 
\end{equation}
for each $i$. Hence, the individual statistics $\widehat{\Delta}_i^{[K_0]}$ all have a rescaled $\chi^2$-distribution.  

\item[$K < K_0$:] If we pick $K$ smaller than the true number of classes $K_0$, the clusters $\{ \widehat{G}_k^{[K]}: 1 \le k \le K \}$ cannot provide an appropriate approximation of the true class structure $\{ G_k: 1 \le k \le K_0 \}$. In particular, there is always a cluster $S = \widehat{G}_k^{[K]}$ which contains subjects from at least two different classes. For simplicity, let $S = G_{k_1} \cup G_{k_2}$. For any $i \in S$, it holds that 
\begin{align*}
\widehat{\Delta}_i^{[K]} 
 & = \frac{1}{\sqrt{p}} \sum\limits_{j = 1}^p \Big\{ \frac{(\varepsilon_{ij} + d_{ij})^2}{\sigma^2} - 1 \Big\} \Big/ \kappa \\
 & = \frac{1}{\sqrt{p}} \sum\limits_{j = 1}^p \frac{d_{ij}^2}{\sigma^2 \kappa} + O_p(1)   
\end{align*}
under our regularity conditions from Section \ref{subsec-asym-ass}. Moreover, it is not difficult to see that for at least one $i \in S$, $p^{-1/2} \sum\nolimits_{j = 1}^p d_{ij}^2 \ge c \sqrt{p}$ for some small constant $c > 0$. This implies that for some $i \in S$,
\begin{equation}\label{delta-K}
\widehat{\Delta}_i^{[K]}  \ge c \sqrt{p} \quad \text{for some } c > 0 \text{ with prob.\ tending to } 1,
\end{equation} 
i.e., the statistic $\widehat{\Delta}_i^{[K]}$ has an explosive behaviour.
\end{itemize}
According to these heuristic considerations, the statistics $\widehat{\Delta}_i^{[K]}$ exhibit a quite different behaviour depending on whether $K < K_0$ or $K = K_0$. When $K < K_0$, the statistic $\widehat{\Delta}_i^{[K]}$ has an explosive behaviour at least for some subjects $i$. This mirrors the fact that a partition with $K < K_0$ clusters cannot give a reasonable approximation to the true class structure. In particular, it cannot describe the data of all subjects $i$ in an accurate way, resulting in an explosive behaviour of the (rescaled) residual sum of squares $\widehat{\Delta}_i^{[K]}$ for some subjects $i$. When $K = K_0$ in contrast, $\{ \widehat{\Delta}_i^{[K_0]}: 1 \le i \le n \}$ is a collection of (approximately) independent random variables that (approximately) have a rescaled $\chi^2$-distribution. Hence, all statistics $\widehat{\Delta}_i^{[K_0]}$ have a stable, non-explosive behaviour. This reflects the fact that the partition $\{ \widehat{G}_k^{[K_0]}: 1 \le k \le K_0 \}$ is an accurate estimate of the true class structure and thus yields a moderate residual sum of squares $\widehat{\Delta}_i^{[K_0]}$ for all subjects $i$.


Since the statistics $\widehat{\Delta}_i^{[K]}$ behave quite differently depending on whether $K = K_0$ or $K < K_0$, they can be used to test $H_0: K = K_0$ versus $H_1: K < K_0$. In particular, testing $H_0$ versus $H_1$ can be achieved by testing the hypothesis that $\widehat{\Delta}_i^{[K]}$ are i.i.d.\ variables with a rescaled $\chi^2$-distribution against the alternative that at least one $\widehat{\Delta}_i^{[K]}$ has an explosive behaviour. We now construct a statistic $\widehat{\mathcal{H}}^{[K]}$ for this testing problem. A natural approach is to take the maximum of the individual statistics $\widehat{\Delta}_i^{[K]}$: Define 
\begin{equation}\label{H-hat-max}
\widehat{\mathcal{H}}^{[K]} = \max_{1 \le i \le n} \widehat{\Delta}_i^{[K]} 
\end{equation}
and let $q(\alpha)$ be the $(1-\alpha)$-quantile of $\mathcal{H} = \max_{1 \le i \le n} Z_i$, where $Z_i$ are independent random variables with the distribution $(\chi_p^2 - p)/\sqrt{2p}$.

Our heuristic discussion from above, in particular formula \eqref{delta-K0}, suggests that for $K = K_0$, 
\[ \pr\Big( \widehat{\mathcal{H}}^{[K_0]} \le q(\alpha) \Big) \approx (1 - \alpha). \]
Moreover, for $K < K_0$, we can show with the help of \eqref{delta-K} and some additional considerations that $\widehat{\mathcal{H}}^{[K]} \ge c \sqrt{p}$ for some $c > 0$ with probability tending to $1$. The quantile $q(\alpha)$, in contrast, can be shown to grow at the rate $\sqrt{\log n}$. Since $\sqrt{\log n } = o(\sqrt{p})$ under our conditions from Section \ref{subsec-asym-ass}, $\widehat{\mathcal{H}}^{[K]}$ diverges faster than the quantile $q(\alpha)$, implying that  
\[ \pr\Big( \widehat{\mathcal{H}}^{[K]} \le q(\alpha) \Big) = o(1) \]
for $K < K_0$. This suggests that $\widehat{\mathcal{H}}^{[K]}$ has the property \eqref{prop-stat} and thus is a reasonable statistic to test the hypothesis $H_0: K = K_0$ versus $H_1: K < K_0$.

In this paper, we restrict attention to the maximum statistic $\widehat{\mathcal{H}}^{[K]}$ defined in \eqref{H-hat-max}. In principle though, we may work with any statistic that satisfies the higher-order property \eqref{prop-stat}. In Section \ref{sec-ext}, we discuss some alternative choices of $\widehat{\mathcal{H}}^{[K]}$.


\subsection{A $\boldsymbol{k}$-means clustering algorithm}\label{subsec-est-kmeans}

We now construct a $k$-means type clustering algorithm which has the consistency property \eqref{prop-clusters}. Since its introduction by \cite{Cox1957} and \cite{Fisher1958}, the $k$-means algorithm has become one of the most popular tools in cluster analysis.
Our version of the algorithm mainly differs from the standard one in the choice of the initial values. To ensure the consistency property \eqref{prop-clusters}, we pick initial clusters $\mathscr{C}_1^{[K]},\ldots,\mathscr{C}_K^{[K]}$ for each given $K$ as follows:

\vspace{8pt}

\noindent \textbf{Choice of the starting values.} Let $i_1,\ldots,i_K$ be indices which (with probability tending to $1$) belong to $K$ different classes $G_{k_1},\ldots,G_{k_K}$ in the case that $K \le K_0$ and to $K_0$ different classes in the case that $K > K_0$. We explain how to obtain such indices below. With these indices at hand, we compute the distance measures $\widehat{\rho}_k(i) = \widehat{\rho}(i_k,i)$ for all $1 \le i \le n$ and $1 \le k \le K$, where 
\[ \widehat{\rho}(i,i^\prime) = \frac{1}{p}\sum\limits_{j=1}^p \big( \widehat{Y}_{ij} - \widehat{Y}_{i^\prime j} \big)^2. \]
The starting values $\mathscr{C}_1^{[K]},\ldots,\mathscr{C}_K^{[K]}$ are now defined by assigning the index $i$ to cluster $\mathscr{C}_k^{[K]}$ if $\widehat{\rho}_k(i) = \min_{1 \le k^\prime \le K} \widehat{\rho}_{k^\prime}(i)$. 

The indices $i_1,\ldots,i_K$ in this construction are computed as follows: For $K=2$, pick any index $i_1 \in \{1,\ldots,n\}$ and calculate $i_2 = \text{arg} \max_{1 \le i \le n} \widehat{\rho}(i_1,i)$. Next suppose we have already constructed the indices $i_1,\ldots,i_{K-1}$ for the case of $K-1$ clusters and compute the corresponding starting values $\mathscr{C}_1^{[K-1]},\ldots,\mathscr{C}_{K-1}^{[K-1]}$ as described above. Calculate the maximal within-cluster distance $\widehat{\rho}_{\max}(k) = \max_{i \in \mathscr{C}_k^{[K-1]}} \widehat{\rho}_k(i)$ for each $1 \le k \le K-1$ and let $\mathscr{C}_{k^*}^{[K-1]}$ be a cluster with $\widehat{\rho}_{\max}(k^*) \ge \widehat{\rho}_{\max}(k)$ for all $k$. Define $i_K = \text{arg} \max_{i \in \mathscr{C}_{k^*}^{[K-1]}} \widehat{\rho}_{k^*}(i)$. 
\vspace{8pt}

\noindent \textbf{The $\boldsymbol{k}$-means algorithm.} Let the number of clusters $K$ be given and denote the starting values by $C_k^{(0)} := \mathscr{C}_k^{[K]}$ for $1 \le k \le K$. The $r$-th iteration of our $k$-means algorithm proceeds as follows: 
\begin{itemize}[leftmargin=1.5cm]
\item[Step $r$:] Let $C_1^{(r-1)},\ldots,C_K^{(r-1)}$ be the clusters from the $(r-1)$-th iteration step. Compute cluster means $m_{j,k}^{(r)} = (\#C_k^{(r-1)})^{-1} \sum\nolimits_{i \in C_k^{(r-1)}} \widehat{Y}_{ij}$ and calculate the distance measures $\widehat{\rho}_k^{(r)}(i) = p^{-1} \sum\nolimits_{j=1}^p (\widehat{Y}_{ij} - m_{j,k}^{(r)})^2$ for all $1 \le i \le n$ and $1 \le k \le K$. Define updated groups $C_1^{(r)},\ldots,C_K^{(r)}$ by assigning the index $i$ to the cluster $C_k^{(r)}$ if $\widehat{\rho}_k^{(r)}(i) = \min_{1 \le k^\prime \le K} \widehat{\rho}_{k^\prime}^{(r)}(i)$. 
\end{itemize}
Repeat this algorithm until the estimated groups do not change any more. For a given sample of data, this is guaranteed to happen after finitely many steps. The resulting $k$-means estimators are denoted by $\{ \widehat{G}_k^{[K]}: 1 \le k \le K \}$. In Section \ref{sec-asym}, we formally show that these estimators have the consistency property \eqref{prop-clusters} under our regularity conditions.

\subsection{Estimation of $\boldsymbol{\sigma^2}$ and $\boldsymbol{\kappa}$}\label{subsec-est-sigma}

In practice, the error variance $\sigma^2$ and the normalization constant $\kappa$ are unknown and need to be estimated from the data at hand. We distinguish between two different estimation approaches, namely a difference- and a residual-based approach. 

\vspace{8pt}

\noindent \textbf{Difference-based estimators.} To start with, consider the time trend model from Section \ref{sec-model}, where $\mu_{ij} = \mu_i(j/p)$ with some trend function $\mu_i(\cdot)$. Supposing that the functions $\mu_i(\cdot)$ are Lipschitz continuous, we get that $Y_{ij} - Y_{i,j-1} = \{ \varepsilon_{ij} - \varepsilon_{i,j-1} \} + \{ \mu_i(j/p) - \mu_i((j-1)/p) \} = \{ \varepsilon_{ij} - \varepsilon_{i,j-1} \} + O(p^{-1})$. This motivates to estimate the error variance $\sigma^2 = \ex[\varepsilon_{ij}^2]$ by 
\[ \widehat{\sigma}^2_{\text{Lip}} = \frac{1}{n(p-1)} \sum\limits_{i=1}^n \sum\limits_{j=2}^p \frac{(Y_{ij} - Y_{i,j-1})^2}{2}. \]
Similarly, the fourth moment $\vartheta = \ex[\varepsilon_{ij}^4]$ can be estimated by 
\[ \widehat{\vartheta}_{\text{Lip}} = \frac{1}{n(p-1)} \sum\limits_{i=1}^n \sum\limits_{j=2}^p \frac{(Y_{ij} - Y_{i,j-1})^4}{2} - 3 (\widehat{\sigma}_{\text{Lip}}^2)^2, \]
which in turn allows us to estimate the parameter $\kappa$ by
\[ \widehat{\kappa}_{\text{Lip}} = \Big( \frac{\widehat{\vartheta}_{\text{Lip}}}{(\widehat{\sigma}_{\text{Lip}}^2)^2} - 1 \Big)^{1/2}. \]
Difference-based estimators of this type have been considered in the context of nonparametric regression by \cite{MuellerStadtmueller1988} and \cite{HallKay1990} among others. Under the technical conditions \ref{A1}--\ref{A3} from Section \ref{subsec-asym-ass}, it is straightforward to show that $\widehat{\sigma}^2_{\text{Lip}} = \sigma^2 + O_p((np)^{-1/2} + p^{-2})$ and $\widehat{\kappa}_{\text{Lip}} = \kappa + O_p((np)^{-1/2} + p^{-2})$. The estimators $\widehat{\sigma}^2_{\text{Lip}}$ and $\widehat{\kappa}_{\text{Lip}}$ are particularly suited for applications where the trend functions $\mu_i(\cdot)$ can be expected to be fairly smooth. This ensures that the unknown first differences $(\varepsilon_{ij} - \varepsilon_{i,j-1})$ can be sufficiently well approximated by the terms $(Y_{ij} - Y_{i,j-1})$. 

A similar difference-based estimation strategy can be used in the model for gene expression microarray data from Section \ref{sec-model}. In this setting, the signal vectors $\boldsymbol{\mu}_i$ have a piecewise constant structure. In particular, $\mu_{ij} = \mu_{ij^\prime}$ for $j_{t-1} \le j,j^\prime < j_t$, where $j_t$ are known indices with $1 = j_0 < j_1 < \ldots < j_{T-1} < j_T = p+1$. This implies that $Y_{ij} - Y_{i,j-1} = \varepsilon_{ij} - \varepsilon_{i,j-1}$ for $j_{t-1} < j < j_t$. Similarly as before, we may thus estimate $\sigma^2$, $\vartheta$ and $\kappa$ by
\begin{align*} 
\widehat{\sigma}^2_{\text{pc}} & = \frac{1}{n (p-T)} \sum\limits_{i=1}^n \sum\limits_{j=1}^p \ind_j \frac{(Y_{ij} - Y_{i,j-1})^2}{2} \\
\widehat{\vartheta}_{\text{pc}} & = \frac{1}{n (p-T)} \sum\limits_{i=1}^n \sum\limits_{j=1}^p \ind_j \frac{(Y_{ij} - Y_{i,j-1})^4}{2} - 3 (\widehat{\sigma}_{\text{pc}}^2)^2  
\end{align*}
and $\widehat{\kappa}_{\text{pc}} = (\widehat{\vartheta}_{\text{pc}}/(\widehat{\sigma}_{\text{pc}}^2)^2 - 1)^{1/2}$, where $\ind_j = \ind( j \notin \{ j_0,j_1,\ldots,j_T \} )$. It is not difficult to see that under the conditions \ref{A1}--\ref{A3}, $\widehat{\sigma}^2_{\text{pc}} = \sigma^2 + O_p((np)^{-1/2})$ and $\widehat{\kappa}_{\text{pc}} = \kappa + O_p((np)^{-1/2})$.

\vspace{8pt}

\noindent \textbf{Residual-based estimators.} Let $\{ \widehat{G}_k^{[K]}: 1 \le k \le K \}$ be the $k$-means estimators from Section \ref{subsec-est-kmeans} for a given $K$. Moreover, let $\widehat{\varepsilon}_{ij}^{[K]}$ be the cluster-specific residuals introduced at the beginning of Section \ref{subsec-est-stat} and denote the vector of residuals for the $i$-th subject by $\widehat{\boldsymbol{\varepsilon}}_i^{[K]} = (\widehat{\varepsilon}_{i1}^{[K]},\ldots,\widehat{\varepsilon}_{ip}^{[K]})^\top$. With this notation at hand, we define the residual sum of squares for $K$ clusters by 
\begin{equation}\label{RSS-K}
\text{RSS}(K) = \frac{1}{np} \sum\limits_{k=1}^K \sum\limits_{i \in \widehat{G}_k^{[K]}} \| \widehat{\boldsymbol{\varepsilon}}_i^{[K]} \|^2, 
\end{equation}
where $\| \cdot \|$ denotes the usual Euclidean norm for vectors. 
$\text{RSS}(K)$ can be shown to be a consistent estimator of $\sigma^2$ for any fixed $K \ge K_0$. The reason is the following: For any $K \ge K_0$, the $k$-means estimators $\widehat{G}_k^{[K]}$ have the property that
\begin{equation}\label{add-prop-clusters}
\pr \Big( \widehat{G}_k^{[K]} \subseteq G_{k^\prime} \text{ for some } 1 \le k^\prime \le K_0 \Big) \rightarrow 1 
\end{equation}
for $1 \le k \le K$ under the technical conditions \ref{A1}--\ref{A3} from Section \ref{subsec-asym-ass}. Hence, with probability tending to $1$, the estimated clusters $\widehat{G}_k^{[K]}$ contain elements from only one class $G_{k^\prime}$. The residuals $\widehat{\varepsilon}_{ij}^{[K]}$ should thus give a reasonable approximation to the unknown error terms $\varepsilon_{ij}$. This in turn suggests that the residual sum of squares $\text{RSS}(K)$ should be a consistent estimator of $\sigma^2$ for $K \ge K_0$. 

Now suppose we know that $K_0$ is not larger than some upper bound $K_{\max}$. In this situation, we may try to estimate $\sigma^2$ by $\widetilde{\sigma}^2_{\text{RSS}} = \text{RSS}(K_{\max})$. Even though consistent, this is a very poor estimator of $\sigma^2$. The issue is the following: The larger $K_{\max}$, the smaller the residual sum of squares $\text{RSS}(K_{\max})$ tends to be. This is a natural consequence of the way in which the $k$-means algorithm works. Hence, if $K_{\max}$ is much larger than $K_0$, then $\widetilde{\sigma}^2_{\text{RSS}} = \text{RSS}(K_{\max})$ tends to strongly underestimate $\sigma^2$. To avoid this issue, we replace the naive estimator $\widetilde{\sigma}^2_{\text{RSS}}$ by a refined version: 
\begin{enumerate}[label=(\roman*),leftmargin=0.85cm]

\item Split the data vector $\boldsymbol{Y}_i = (Y_{i1},\ldots,Y_{ip})^\top$ into the two parts $\boldsymbol{Y}_i^{A} = (Y_{i1},Y_{i3},\ldots)^\top$ and $\boldsymbol{Y}_i^{B} = (Y_{i2},Y_{i4},\ldots)^\top$. Moreover, let $\overline{Y}_i^A$ be the empirical mean of the entries in the vector $\boldsymbol{Y}_i^A$ and define $\widehat{\boldsymbol{Y}}_i^A = (Y_{i1} - \overline{Y}_i^A,Y_{i3} - \overline{Y}_i^A,\ldots)^\top$. Finally, set $\mathcal{Y}^A = \{ \widehat{\boldsymbol{Y}}_i^A: 1 \le i \le n \}$ and analogously define $\mathcal{Y}^B = \{ \widehat{\boldsymbol{Y}}_i^B: 1 \le i \le n \}$. Importantly, under our technical conditions, the random vectors in $\mathcal{Y}^A$ are independent from those in $\mathcal{Y}^B$.  

\item Apply the $k$-means algorithm with $K = K_{\max}$ to the sample $\mathcal{Y}^A$ and denote the resulting estimators by $\{ \widehat{G}_k^A: 1 \le k \le K_{\max} \}$. These estimators can be shown to have the property \eqref{add-prop-clusters}, provided that we impose the following condition: Let $\boldsymbol{m}_k$ be the class-specific signal vector of the class $G_k$ and define the vectors $\boldsymbol{m}_k^A$ and $\boldsymbol{m}_k^B$ in the same way as above. Assume that 
\begin{equation}\label{ass-sigma-RSS}
\boldsymbol{m}_k^A \ne \boldsymbol{m}_{k^\prime}^A \text{ for } k \ne k^\prime. 
\end{equation}
According to this assumption, the signal vectors $\boldsymbol{m}_k$ and $\boldsymbol{m}_{k^\prime}$ of two different classes can be distinguished from each other only by looking at their odd entries $\boldsymbol{m}_k^A$ and $\boldsymbol{m}_{k^\prime}^A$. It goes without saying that this is not a very severe restriction. 

\item Compute cluster-specific residuals from the data sample $\mathcal{Y}^B$,
\[ \widehat{\boldsymbol{\varepsilon}}_i^B = \widehat{\boldsymbol{Y}}_i^B - \frac{1}{\#\widehat{G}_k^A} \sum\limits_{i^\prime \in \widehat{G}_k^A} \widehat{\boldsymbol{Y}}_{i^\prime}^B \quad \text{for } i \in \widehat{G}_k^A, \]
and define  
\[ \widehat{\sigma}^2_{\text{RSS}} = \frac{1}{n \lfloor p/2 \rfloor} \sum\limits_{k=1}^{K_{\max}} \sum\limits_{i \in \widehat{G}_k^A} \| \widehat{\boldsymbol{\varepsilon}}_i^B \|^2. \]
In contrast to the naive estimator $\widetilde{\sigma}^2_{\text{RSS}}$, the refined version $\widehat{\sigma}^2_{\text{RSS}}$ does not tend to strongly underestimate $\sigma^2$. The main reason is that the residuals $\widehat{\boldsymbol{\varepsilon}}_i^B$ are computed from the random vectors $\widehat{\boldsymbol{Y}}_i^B$ which are independent of the estimated clusters $\widehat{G}_k^A$.

\end{enumerate}

\noindent Writing $\widehat{\boldsymbol{\varepsilon}}_i^B = (\widehat{\varepsilon}_{i1}^B,\ldots,\widehat{\varepsilon}_{i\lfloor p/2 \rfloor}^B)^\top$, we can analogously estimate the fourth error moment $\vartheta = \ex[\varepsilon_{ij}^4]$ by
\[ \widehat{\vartheta}_{\text{RSS}} = \frac{1}{n \lfloor p/2 \rfloor} \sum\limits_{k=1}^{K_{\max}} \sum\limits_{i \in \widehat{G}_k^A} \sum\limits_{j=1}^{\lfloor p/2 \rfloor} \big( \widehat{\varepsilon}_{ij}^B \big)^4 \]
and set $\widehat{\kappa}_{\text{RSS}} = (\widehat{\vartheta}_{\text{RSS}}/(\widehat{\sigma}^2_{\text{RSS}})^2 - 1)^{1/2}$. Under the conditions \ref{A1}--\ref{A3}, it can be shown that 
\begin{equation}\label{rate-sigma-RSS}
\widehat{\sigma}^2_{\text{RSS}} = \sigma^2 + O_p(p^{-1}) \quad \text{and} \quad \widehat{\kappa}_{\text{RSS}} = \kappa + O_p(p^{-1}). 
\end{equation}
A sketch of the proof is given in the Supplementary Material.

\section{Asymptotics}\label{sec-asym}

In this section, we investigate the asymptotic properties of our estimators. We first list the assumptions needed for the analysis and then summarize the main results.

\subsection{Assumptions}\label{subsec-asym-ass}

To formulate the technical conditions that we impose on model \eqref{model-eq1}--\eqref{model-eq2}, we denote the size, i.e., the cardinality of the class $G_k$ by $n_k = \# G_k$. Moreover, we use the shorthand $a_\nu \ll b_\nu$ to express that $a_\nu/b_\nu \le c \nu^{-\delta}$ for sufficiently large $\nu$ with some $c > 0$ and a small $\delta > 0$. Our assumptions read as follows: 
\begin{enumerate}[label=(C\arabic*),leftmargin=1.05cm]
\item \label{A1} The errors $\varepsilon_{ij}$ are identically distributed and independent across both $i$ and $j$ with $\ex[\varepsilon_{ij}] = 0$ and $\ex[|\varepsilon_{ij}|^{\theta}] \le C < \infty$ for some $\theta > 8$. 
\item \label{A2} The class-specific signal vectors $\boldsymbol{m}_k = (m_{1,k},\ldots,m_{p,k})^\top$ differ across groups in the following sense: There exists a constant $\delta_0 > 0$ such that
\[ \frac{1}{p} \sum\limits_{j=1}^p \big( m_{j,k} - m_{j,k^\prime} \big)^2 \ge \delta_0 \]
for any pair of groups $G_k$ and $G_{k^\prime}$ with $k \ne k^\prime$. Moreover, $|m_{j,k}| \le C$ for all $k$ and $j$, where $C > 0$ is a sufficiently large constant. 
\item \label{A3} Both $n$ and $p$ tend to infinity. The group sizes $n_k = \# G_k$ are such that $p \ll n_k \ll p^{(\theta/4) - 1}$ for all $1 \le k \le K_0$, implying that $p \ll n \ll p^{(\theta/4) - 1}$.
\end{enumerate}
We briefly comment on the above conditions. By imposing \ref{A1}, we restrict the noise terms $\varepsilon_{ij}$ to be i.i.d. Yet the error terms $e_{ij} = \alpha_i + \varepsilon_{ij}$ of our model may be dependent across subjects $i$, as we do not impose any restrictions on the random intercepts $\alpha_i$. This is important, for instance, when clustering the genes in a microarray data set, where we may expect different genes $i$ to be correlated. \ref{A2} is a fairly harmless condition, which requires the signal vectors to differ in an $L_2$-sense across groups. By \ref{A3}, the group sizes $n_k$ and thus the total number of subjects $n$ are supposed to grow faster than the number of features $p$. We thus focus attention on applications where $n$ is (much) larger than $p$. However, $n$ should not grow too quickly as compared to $p$. Specifically, it should not grow faster than $p^{(\theta/4) - 1}$, where $\theta$ is the number of existing error moments $\ex[|\varepsilon_{ij}|^{\theta}] < \infty$. As can be seen, the bound $p^{(\theta/4) - 1}$ on the growth rate of $n$ gets larger with increasing $\theta$. In particular, if all moments of $\varepsilon_{ij}$ exist, $n$ may grow as quickly as any polynomial of $p$. Importantly, \ref{A3} allows the group sizes $n_k$ to grow at different rates (between $p$ and $p^{(\theta/4) - 1}$). Put differently, it allows for strongly heterogeneous group sizes. Our estimation methods are thus able to deal with situations where some groups are much smaller than others.

\subsection{Main results}

Our first result shows that the maximum statistic $\widehat{\mathcal{H}}^{[K]} = \max_{1 \le i \le n} \widehat{\Delta}_i^{[K]}$ has the property \eqref{prop-stat} and thus is a reasonable statistic to test the hypothesis $H_0: K = K_0$ versus $H_1: K < K_0$.
\begin{theorem}\label{theo1}
Assume that the estimated clusters have the consistency property \eqref{prop-clusters}. Moreover, let $\widehat{\sigma}^2$ and $\widehat{\kappa}$ be any estimators with $\widehat{\sigma}^2 = \sigma^2 + O_p(p^{-(1/2 + \delta)})$ and $\widehat{\kappa} = \kappa + O_p(p^{-\delta})$ for some $\delta > 0$. Under \ref{A1}--\ref{A3}, the statistic $\widehat{\mathcal{H}}^{[K]}$ has the property \eqref{prop-stat}, that is,  
\begin{equation*}
\pr\Big( \widehat{\mathcal{H}}^{[K]} \le q(\alpha) \Big) = 
\begin{cases}
o(1) & \text{for } K < K_0 \\
(1 - \alpha) + o(1) & \text{for } K = K_0.
\end{cases} 
\end{equation*}
\end{theorem}
\noindent This theorem is the main stepping stone to derive the central result of the paper, which describes the asymptotic properties of the estimators $\widehat{K}_0$ and $\{ \widehat{G}_k: 1 \le k \le \widehat{K}_0 \}$.     
\begin{theorem}\label{theo2}
Under the conditions of Theorem \ref{theo1}, it holds that 
\begin{equation*}
\pr \big( \widehat{K}_0 > K_0 \big) = \alpha + o(1) \quad \text{and} \quad \pr \big( \widehat{K}_0 < K_0 \big) = o(1), 
\end{equation*}
implying that $\pr(\widehat{K}_0 \ne K_0) = \alpha + o(1)$. Moreover,  
\begin{equation*}
\pr \Big( \big\{ \widehat{G}_k: 1 \le k \le \widehat{K}_0 \big\} \ne \big\{ G_k: 1 \le k \le K_0 \big\} \Big) = \alpha + o(1). 
\end{equation*}
\end{theorem}
\noindent Theorem \ref{theo2} holds true for any clustering algorithm with the consistency property \eqref{prop-clusters}. The next result shows that this property is fulfilled, for example, by the $k$-means algorithm from Section \ref{subsec-est-kmeans}. 
\begin{theorem}\label{theo3}
Under \ref{A1}--\ref{A3}, the $k$-means estimators $\{ \widehat{G}_k^{[K_0]}: 1 \le k \le K_0 \}$ from Section \ref{subsec-est-kmeans} satisfy \eqref{prop-clusters}, that is, 
\begin{equation*}
\pr \Big( \big\{ \widehat{G}_k^{[K_0]}: 1 \le k \le K_0 \big\} = \big\{ G_k: 1 \le k \le K_0 \big\} \Big) \rightarrow 1. 
\end{equation*}
\end{theorem} 
\noindent The proofs of Theorems \ref{theo1}--\ref{theo3} can be found in the Supplementary Material.

\section{Applications and Simulation Study}\label{sec-app}

\subsection{Clustering of temperature curves}\label{subsec-app-3}

\begin{figure}[!t]
\vspace{-1.5cm}
\centering
\includegraphics[scale = 0.6]{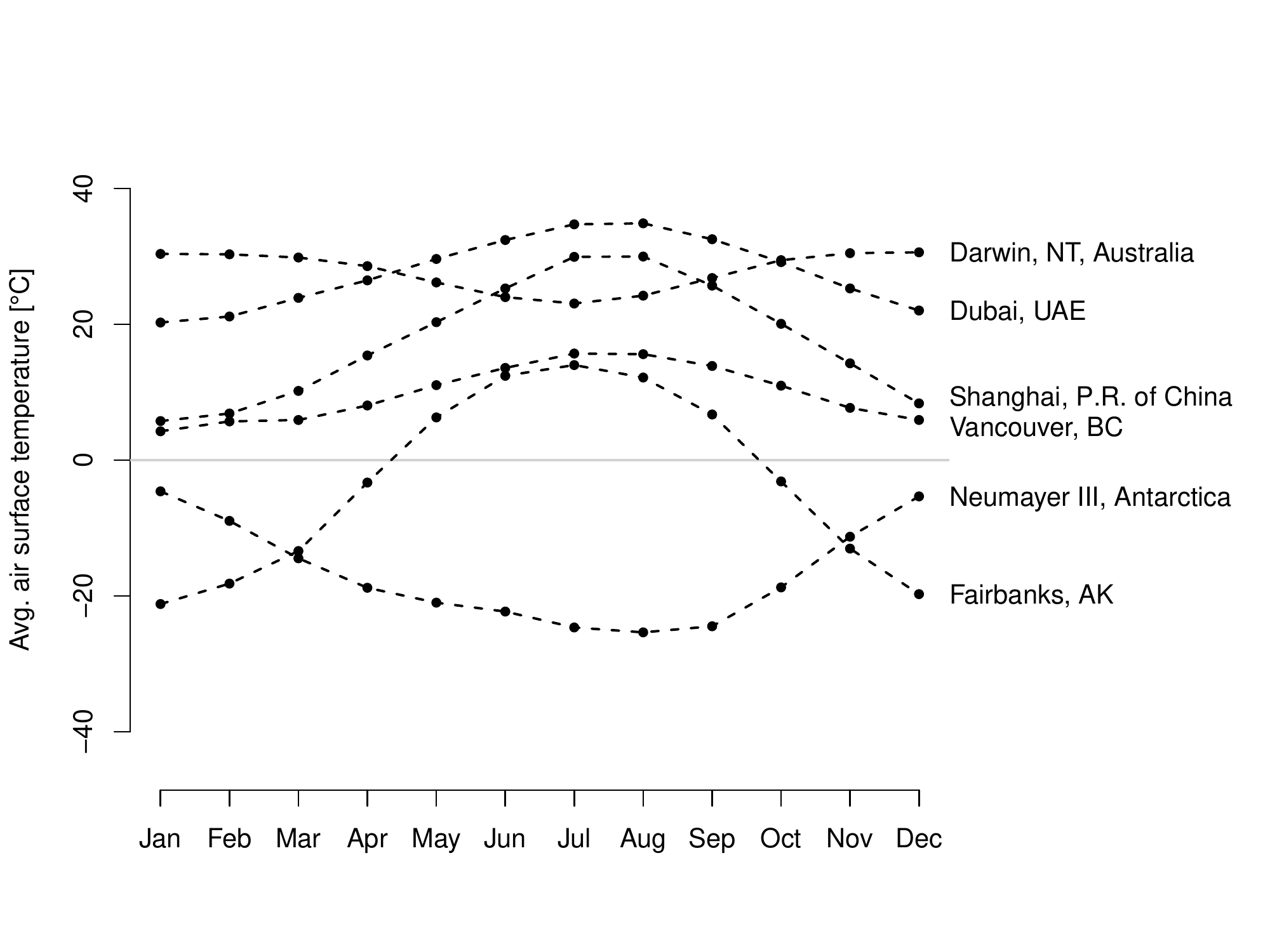}
\vspace{-0.5cm}

\caption{Analysis of the Berkeley Earth temperature data. The plot depicts the average land surface temperature curves at various example locations on earth. 
} \label{fig:ExampleTemp}
\end{figure}

Our first application is concerned with the analysis of a data set on land surface temperatures that was collected by the investigators of the Berkeley Earth project \citep{berkeley}. The data, which are publicly available at \texttt{http://berkeley\linebreak earth.org/data}, contain measurements on a grid of worldwide locations that is defined on a one degree (longitude) by one degree (latitude) basis. For each grid point, the data set contains a monthly land surface temperature profile. This profile is a vector with twelve entries, the first entry specifying the average temperature of all Januaries from 1951 to 1980, the second entry specifying the average temperature of all Februaries from 1951 to 1980, and so on. The temperature profiles at various example locations on earth are shown in Figure \ref{fig:ExampleTemp}. As grid points containing $100\%$ sea surface are not taken into account, the overall number of grid points is equal to $n = \text{24,311}$. A detailed description of the derivation of the data can be found in \cite{berkeley}. Our analysis is based on the Berkeley Earth source file from April 19, 2016.

The aim of our analysis is to cluster the $\text{24,311}$ grid points in order to obtain a set of climate regions characterized by distinct temperature profiles. For this purpose, we impose the time trend model \eqref{trend-model} on the data and apply the CluStErr method to them, setting $n = \text{24,311}$, $p = 12$ and $\alpha = 0.05$. To estimate the error variance $\sigma^2$ and the normalization constant $\kappa$, we apply the difference-based estimators $\widehat{\sigma}^2_{\text{Lip}}$ and $\widehat{\kappa}_{\text{Lip}}$ from Section \ref{subsec-est-sigma}, thus making use of the smoothness of the temperature curves illustrated in Figure \ref{fig:ExampleTemp}.

\begin{figure}[!t]
\vspace{-0.4cm}
\centering
\includegraphics[scale = 0.6]{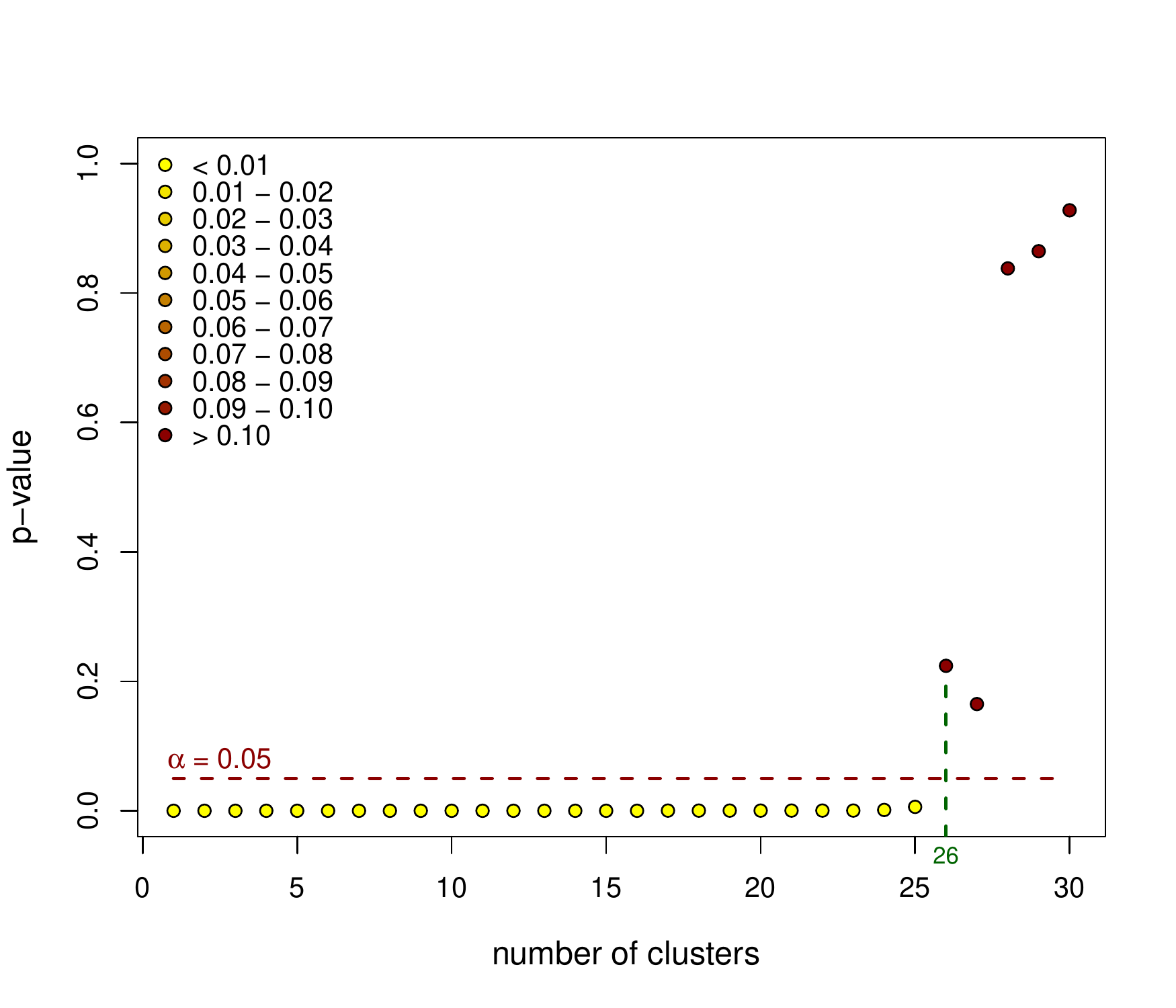}
\caption{Analysis of the Berkeley Earth temperature data. The plot depicts the $p$-values $\widehat{p}^{[K]}$ corresponding to the test statistics $\widehat{\mathcal{H}}^{[K]}$ as a function of $K$. The horizontal dashed line specifies the significance level $\alpha = 0.05$, and the vertical dashed line indicates that the estimated number of clusters is $\widehat{K}_0 = 26$.} \label{fig:QuantileTemp}
\end{figure}

\begin{figure}[!ht]
\vspace{-0.8cm}

\hspace{-0.9cm}
\includegraphics[scale = 0.4]{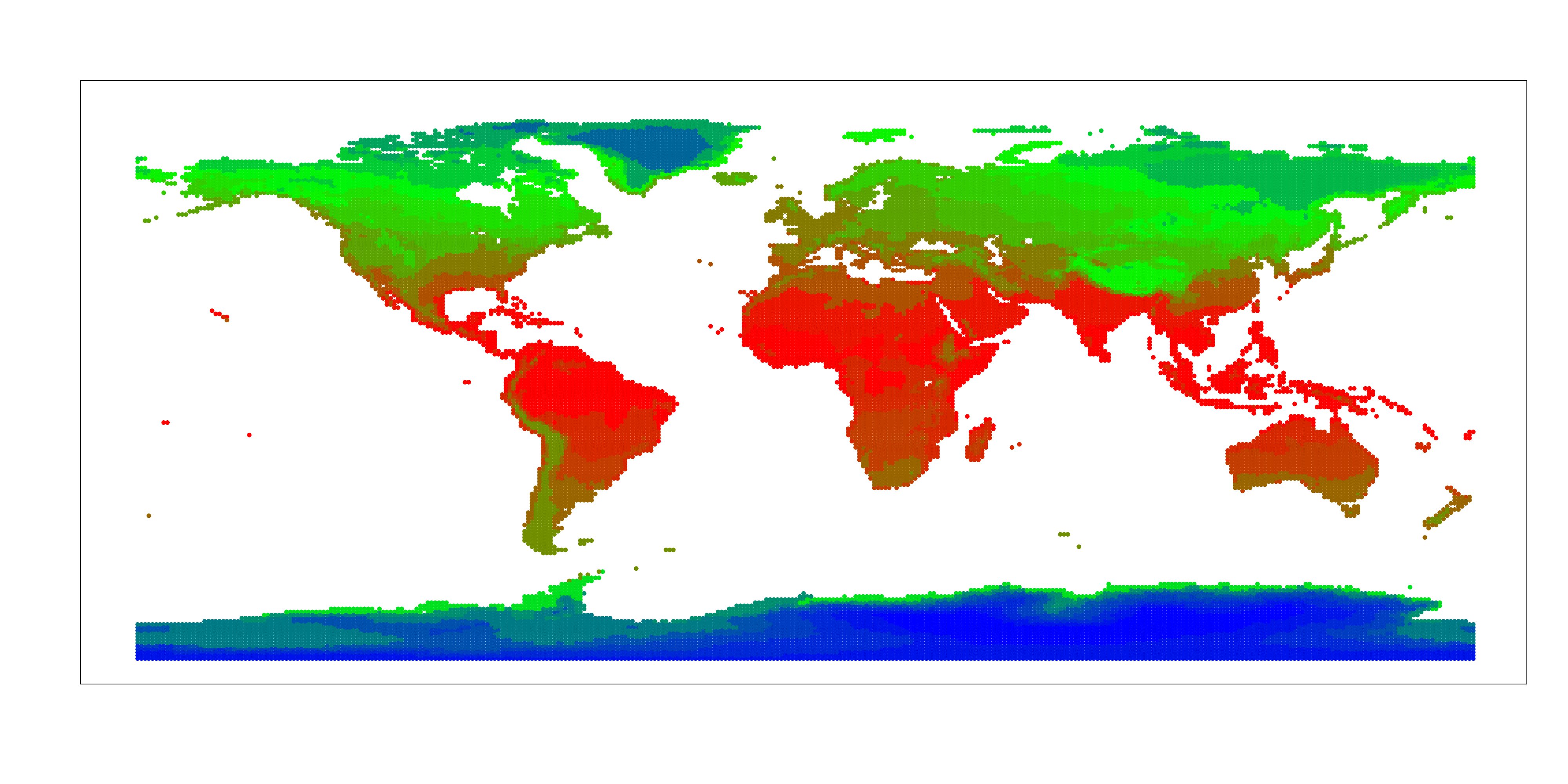}
\vspace{-1.25cm}

\caption{Visualization of the $\widehat{K}_0 = 26$ clusters obtained from the analysis of the Berkeley Earth temperature data. Each shade of color refers to one cluster.
} \label{fig:ClustersTemp}
\end{figure}

\newpage

The estimation results are presented in Figures~\ref{fig:QuantileTemp} and \ref{fig:ClustersTemp}. Figure \ref{fig:QuantileTemp} depicts the $p$-values $\widehat{p}^{[K]}$ corresponding to the test statistics $\widehat{\mathcal{H}}^{[K]}$ for different numbers of clusters $K$. It shows that the $p$-value $\widehat{p}^{[K]}$ remains below the $\alpha = 0.05$ threshold for any $K < 26$ but jumps across this threshold for $K = 26$. The CluStErr algorithm thus estimates the number of clusters to be equal to $\widehat{K}_0 = 26$, suggesting that there are $26$ distinct climate regions. The sizes of the estimated clusters range between $244$ and $\text{2,110}$; the error variance is estimated to be $\widehat{\sigma}^2 = 16.25$. Figure~\ref{fig:ClustersTemp} uses a spatial grid to visualize the 26 regions and demonstrates the plausibility of the obtained results. For example, mountain ranges such as the Himalayas and the South American Andes, but also tropical climates in Africa, South America and Indonesia are easily identified from the plot. Of note, the results presented in Figure \ref{fig:ClustersTemp} show a remarkable similarity to the most recent modification of the K\"oppen-Geiger classification, which is one of the most widely used classification systems in environmental research \citep{kgs}. In particular, the overall number of climate regions defined in \cite{kgs} is equal to 29, which is similar to the cluster number $\widehat{K}_0 = 26$ identified by the CluStErr algorithm. Thus, although our example is purely illustrative, and although expert classification systems account for additional characteristics such as precipitation and vegetation, Figure \ref{fig:ClustersTemp} confirms the usefulness of the CluStErr method.

\subsection{Clustering of gene expression data}\label{subsec-app-1}

Our second application is concerned with the analysis of gene expression data, which has become a powerful tool for the understanding of disease processes in biomedical research \citep{jiang}. A popular approach to measure gene expression is to carry out microarray experiments. These experiments simultaneously quantify the expression levels of $n$ genes across $p$ samples of patients with different clinical conditions, such as tumor stages or disease subtypes. In the analysis of microarray data, clustering of the $n$ genes is frequently used to detect genes with similar cellular function and to discover groups of ``co-expressed'' genes showing similar expression patterns across clinical conditions \citep{chipman, jiang, dhaese}.

\begin{figure}[!t]
\vspace{-0.5cm}
\centering
\includegraphics[scale = 0.5]{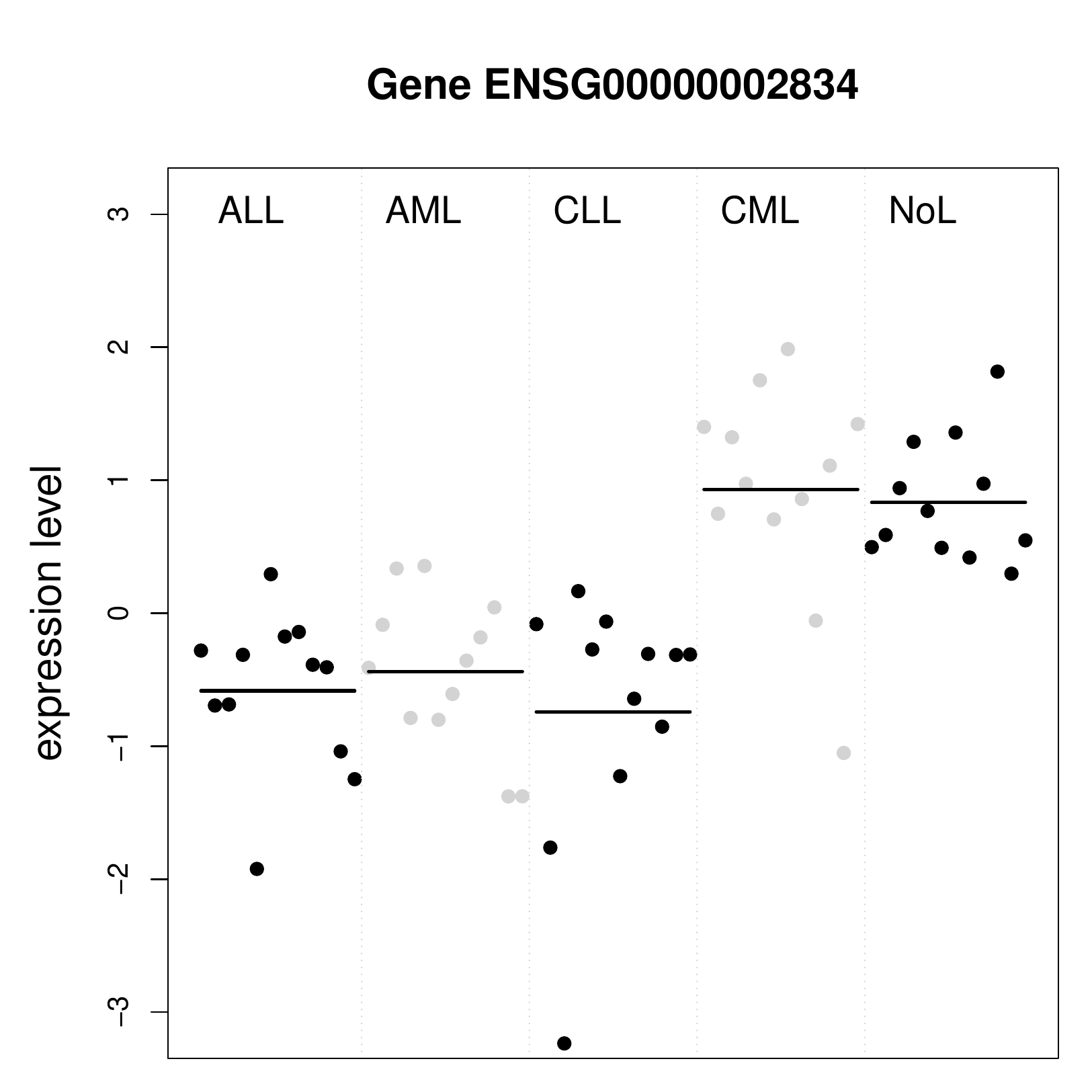}
\caption{Analysis of the MILE study gene expression data. The plot depicts the expression levels of a randomly selected gene after normalization and standardization. The label of the gene (``ENSG00000002834'') refers to its {\it Ensembl} gene ID to which the original Affymetrix probesets were mapped \citep{leukemiasEset}. Horizontal lines represent the average gene expression levels across the five tissue types.
} \label{fig:IlluGene}
\end{figure}

\begin{figure}[!t]
\phantom{upperboundary}
\vspace{-1.5cm}
\centering
\includegraphics[scale = 0.6]{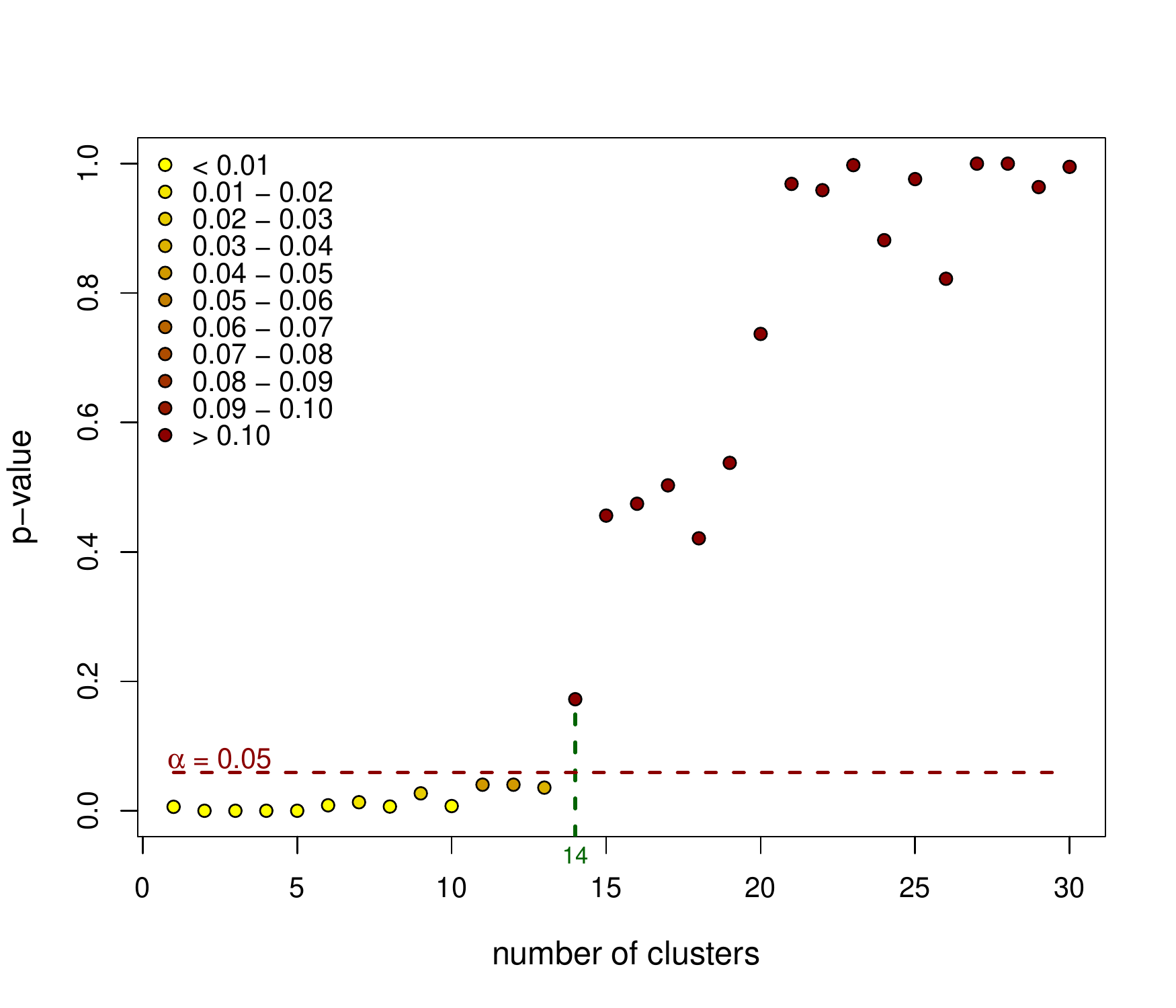}
\caption{Analysis of the MILE study gene expression data. The plot depicts the $p$-values $\widehat{p}^{[K]}$ corresponding to the test statistics $\widehat{\mathcal{H}}^{[K]}$ as a function of $K$. The dashed vertical line indicates that the number of clusters is estimated to be $\widehat{K}_0 = 14$. 
}
\label{fig:QuantileGene}
\end{figure}

\begin{figure}[!hp]
\centering
\includegraphics[scale = 0.55]{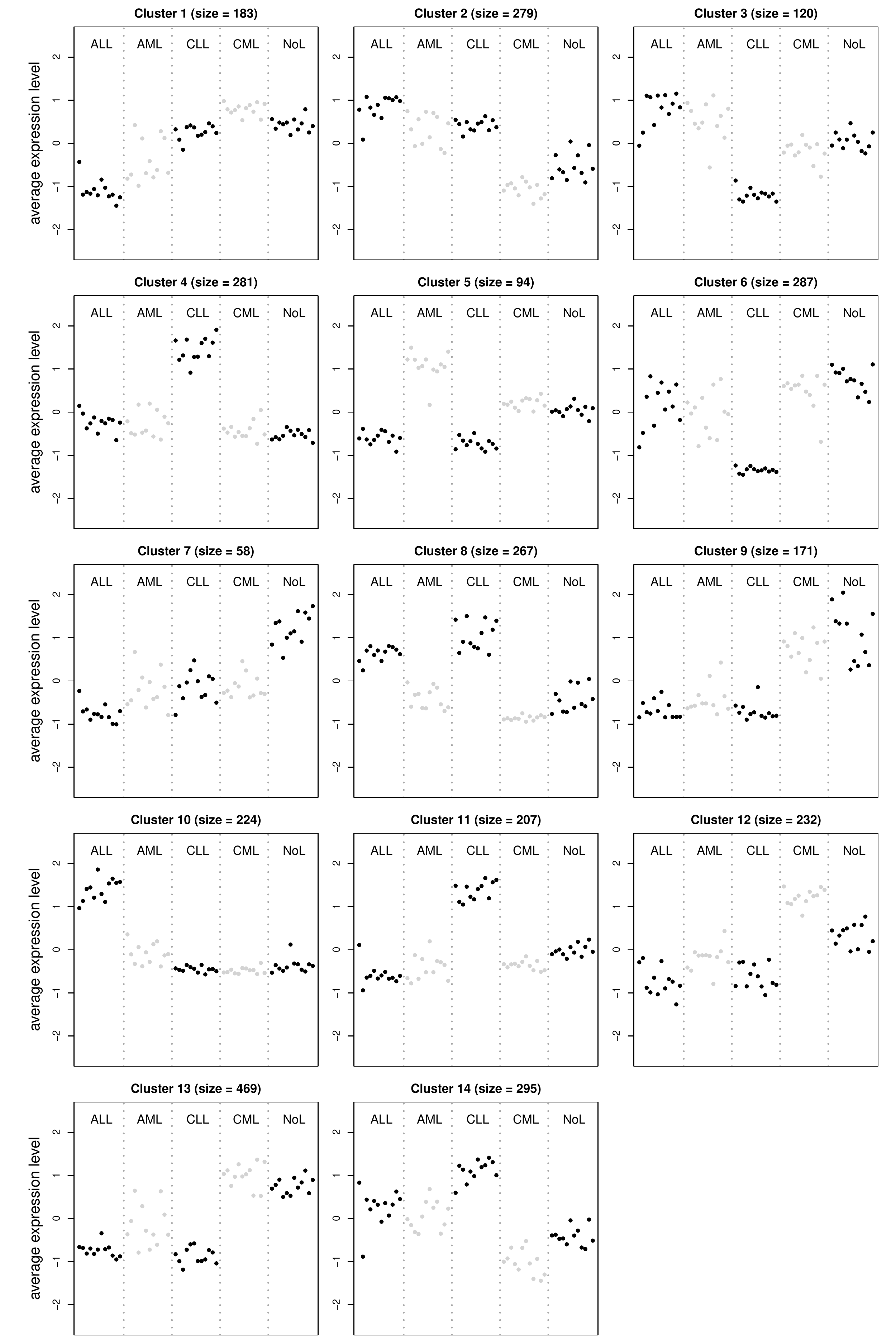}
\vspace{0.25cm}

\caption{Visualization of the cluster centres obtained from the analysis of the MILE study gene expression data. The dots represent the cluster centres $\boldsymbol{\widehat{m}}_k = (\# \widehat{G}_k)^{-1} \sum\nolimits_{i \in \widehat{G}_k} \boldsymbol{\widehat{Y}}_i$, which estimate the cluster-specific signal vectors $\boldsymbol{m}_k = (\# G_k)^{-1} \sum\nolimits_{i \in G_k} \boldsymbol{\mu}_i$.
Gene ENSG00000002834, whose expression profile is visualized in Figure \ref{fig:IlluGene}, is an element of cluster $\#13$. Note the similarity of the patterns in cluster $\#13$ and Figure \ref{fig:IlluGene}.} \label{fig:ClustersGene}
\end{figure}

In what follows, we analyze a set of gene expression data that was collected during the first stage of the Microarray Innovations in Leukemia (MILE) study \citep{haferlach}. The data set contains expression level measurements for $\text{20,172}$ genes and is publicly available as part of the Bioconductor package {\bf leukemiasEset} \citep{leukemiasEset}. The gene expression levels were measured using Affymetrix HG-U133 Plus 2.0 microarrays. For statistical analysis, the raw expression data were normalized using the Robust Multichip Average (RMA) method, followed by an additional gene-wise standardization of the expression levels. For details on data collection and pre-processing, we refer to \cite{haferlach} and \cite{leukemiasEset}.

The data of the MILE study were obtained from $p = 60$ bone marrow samples of patients that were untreated at the time of diagnosis. Of these patients, $48$ were either diagnosed with acute lymphoblastic leukemia (ALL, 12 patients), acute myeloid leukemia (AML, 12 patients), chronic lymphocytic leukemia (CLL, 12 patients), or chronic myeloid leukemia (CML, 12 patients). The other $12$ samples were obtained from non-leukemia (NoL) patients. From a biomedical point of view, the main interest focuses on the set of ``differentially expressed'' genes, that is, on those genes that show a sufficient amount of variation in their expression levels across the five tissue types (ALL, AML, CLL, CML, NoL). To identify the set of these genes, we run a univariate ANOVA for each gene and discard those with Bonferroni-corrected $p$-values $\ge 0.01$ in the respective overall $F$-tests. Application of this procedure results in a sample of $n = \text{3,167}$ univariately significant genes.

The aim of our analysis is to cluster the $n = \text{3,167}$ genes into groups whose members have similar expression patterns across the five tissue types (ALL, AML, CLL, CML, NoL). To do so, we impose model \eqref{gene-model-1a} from Section \ref{sec-model} on the data. The measured expression profiles $\boldsymbol{Y}_i = (Y_{i1},\ldots,Y_{ip})^\top$ of the various genes $i =1,\ldots,n$ are thus assumed to follow the model equation $\boldsymbol{Y}_i = \boldsymbol{\mu}_i + \boldsymbol{\alpha}_i + \boldsymbol{\varepsilon}_i$. The signal vectors $\boldsymbol{\mu}_i$ are supposed to have a piecewise constant structure  after the patients have been ordered according to their tissue type (ALL, AML, CLL, CML, NoL). For illustration, the expression profile $\boldsymbol{Y}_i$ of a randomly selected gene is plotted in Figure \ref{fig:IlluGene}.

To cluster the genes, we apply the CluStErr algorithm with the significance level $\alpha = 0.05$ and the difference-based estimators $\widehat{\sigma}^2_{\text{pc}}$ and $\widehat{\kappa}_{\text{pc}}$ from Section \ref{subsec-est-sigma}, thus exploiting the piecewise constant structure of the signal vectors. The estimation results are presented in Figures \ref{fig:QuantileGene} and \ref{fig:ClustersGene}. The plot in Figure \ref{fig:QuantileGene} depicts the $p$-values $\widehat{p}^{[K]}$ corresponding to the test statistics $\widehat{\mathcal{H}}^{[K]}$ as a function of the cluster number $K$. It shows that the estimated number of clusters is $\widehat{K}_0 = 14$. The estimated sizes of the $14$ clusters range between $58$ and $469$. Moreover, the estimated error variance is $\widehat{\sigma}^2 = 0.442$.  In Figure \ref{fig:ClustersGene}, the cluster centres $\boldsymbol{\widehat{m}}_k = (\# \widehat{G}_k)^{-1} \sum\nolimits_{i \in \widehat{G}_k} \boldsymbol{\widehat{Y}}_i$ are presented, which estimate the cluster-specific signal vectors $\boldsymbol{m}_k = (\# G_k)^{-1} \sum\nolimits_{i \in G_k} \boldsymbol{\mu}_i$. All clusters show a distinct separation of at least one tissue type, supporting the assumption of piecewise constant signals $\boldsymbol{m}_k$ and indicating that the genes contained in the clusters are co-expressed differently across the five groups. For example, cluster $\#2$ separates CML and NoL samples from ALL, AML and CLL samples, whereas cluster $\#4$ separates CLL samples from the other tissue types. Thus, each of the $14$ clusters represents a specific pattern of co-expressed gene profiles.

\subsection{Simulation study}\label{subsec-app-2}

\begin{figure}[!th]
\hspace{-0.6cm}
\includegraphics[scale = 0.675]{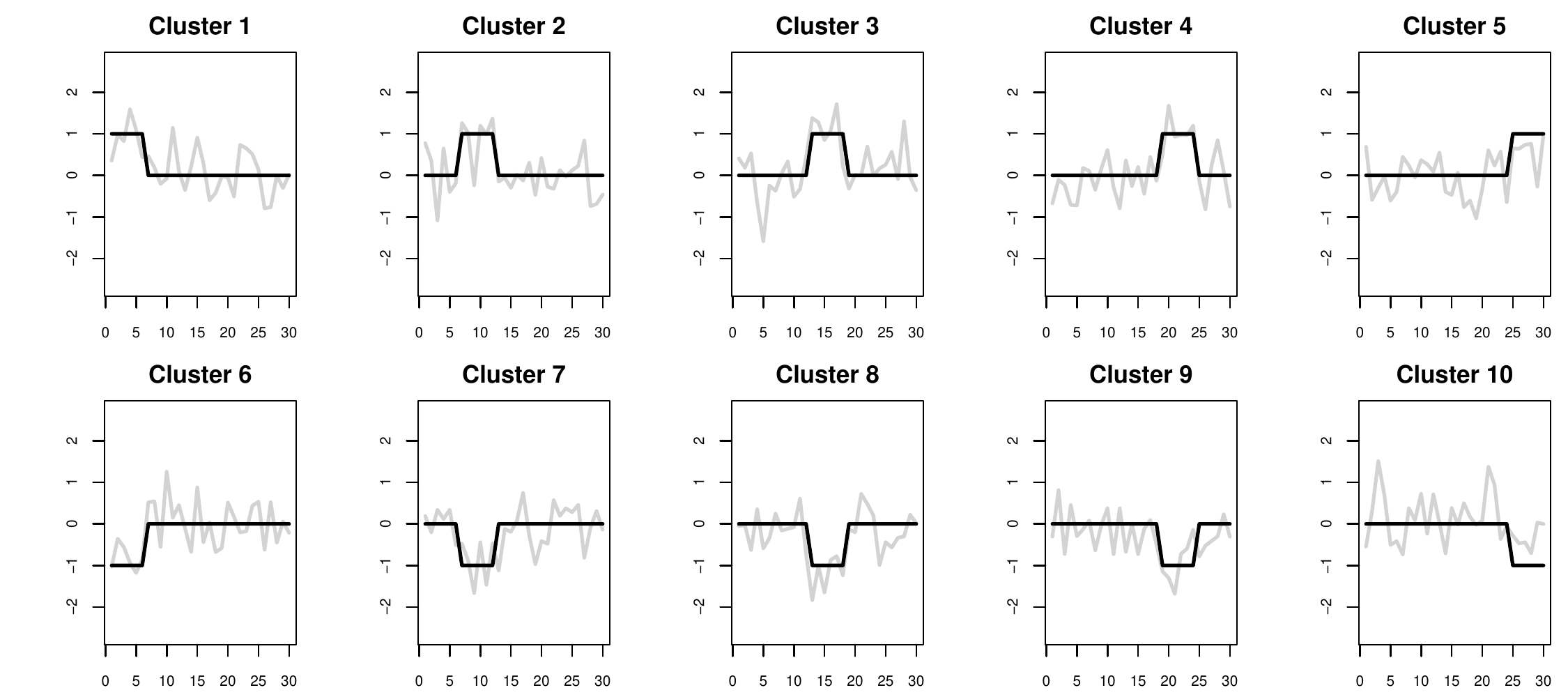}
\caption{Visualization of the cluster-specific signal vectors $\boldsymbol{m}_k$ for the simulation study. The black lines represent the signal vectors $\boldsymbol{m}_k$ for $k=1,\ldots,K_0=10$. The gray lines depict the data vectors $\boldsymbol{Y}_i = \boldsymbol{m}_k + \boldsymbol{\varepsilon}_i$ of a randomly selected member $i$ of the $k$-th cluster for each $k$. All plots are based on a setting with $p = 30$ and noise-to-signal ratio $\textnormal{NSR} = 1.5$.
} \label{fig:IlluSimu}
\end{figure}

To explore the properties of the CluStErr method more systematically, we carry out a simulation study which splits into two main parts. The first part investigates the finite sample behaviour of CluStErr, whereas the second part compares CluStErr with several competing methods. The simulation design is inspired by the analysis of the gene expression data in Section \ref{subsec-app-1}. It is based on model \eqref{gene-model-1a} from Section \ref{sec-model}. The data vectors $\boldsymbol{Y}_i$ have the form $\boldsymbol{Y}_i = \boldsymbol{\mu}_i + \boldsymbol{\varepsilon}_i$ with piecewise constant signal profiles $\boldsymbol{\mu}_i$. We set the number of clusters to $K_0 = 10$ and define the cluster-specific signal vectors $\boldsymbol{m}_k$ by 
\begin{align*}
\boldsymbol{m}_1 & = (\boldsymbol{1},\boldsymbol{0},\boldsymbol{0},\boldsymbol{0},\boldsymbol{0})^\top , & \boldsymbol{m}_6 & = (-\boldsymbol{1},\boldsymbol{0},\boldsymbol{0},\boldsymbol{0},\boldsymbol{0})^\top , \\
\boldsymbol{m}_2 & = (\boldsymbol{0},\boldsymbol{1},\boldsymbol{0},\boldsymbol{0},\boldsymbol{0})^\top ,& \boldsymbol{m}_7 & = (\boldsymbol{0},-\boldsymbol{1},\boldsymbol{0},\boldsymbol{0},\boldsymbol{0})^\top ,\\
\boldsymbol{m}_3 & = (\boldsymbol{0},\boldsymbol{0},\boldsymbol{1},\boldsymbol{0},\boldsymbol{0})^\top ,& \boldsymbol{m}_8 & = (\boldsymbol{0},\boldsymbol{0},-\boldsymbol{1},\boldsymbol{0},\boldsymbol{0})^\top ,\\
\boldsymbol{m}_4 & = (\boldsymbol{0},\boldsymbol{0},\boldsymbol{0},\boldsymbol{1},\boldsymbol{0})^\top ,& \boldsymbol{m}_9 & = (\boldsymbol{0},\boldsymbol{0},\boldsymbol{0},-\boldsymbol{1},\boldsymbol{0})^\top ,\\
\boldsymbol{m}_5 & = (\boldsymbol{0},\boldsymbol{0},\boldsymbol{0},\boldsymbol{0},\boldsymbol{1})^\top  ,& \boldsymbol{m}_{10} & = (\boldsymbol{0},\boldsymbol{0},\boldsymbol{0},\boldsymbol{0},-\boldsymbol{1})^\top,  
\end{align*}
where $\boldsymbol{1} = (1,\ldots,1)$ and $\boldsymbol{0} = (0,\ldots,0)$ are vectors of length $p/5$. A graphical illustration of the signal vectors $\boldsymbol{m}_k$ is provided in Figure \ref{fig:IlluSimu}. The error terms $\varepsilon_{ij}$ are assumed to be i.i.d.\ normally distributed with mean $0$ and variance $\sigma^2$. In the course of the simulation study, we consider different values of $n$, $p$ and $\sigma^2$ as well as different cluster sizes. To assess the noise level in the simulated data, we consider the ratios between the error variance $\sigma^2$ and the ``variances'' of the signals $\boldsymbol{m}_k$. In particular, we define the noise-to-signal ratios $\text{NSR}_k(\sigma^2) = \sigma^2/\var(\boldsymbol{m}_k)$, where $\var(\boldsymbol{m}_k)$ denotes the empirical variance of the vector $\boldsymbol{m}_k$. Since $\var(\boldsymbol{m}_k) \approx 0.16$ is the same for all $k$ in our design, we obtain that $\text{NSR}_k(\sigma^2) = \text{NSR}(\sigma^2) \approx \sigma^2/0.16$ for all $k$.  
 
\vspace{10pt}

\noindent \textbf{Finite sample properties of CluStErr.} In this part of the simulation study, we analyze a design with equally sized clusters and set the sample size to $(n,p) = (1000,30)$. Three different noise-to-signal ratios $\textnormal{NSR}$ are considered, in particular $\textnormal{NSR} = 1$, $1.5$ and $2$. Since $\sigma^2 \approx 0.16 \, \text{NSR}$, the corresponding error variances amount to $\sigma^2 \approx 0.16$, $0.25$ and $0.32$, respectively. The noise-to-signal ratio $\textnormal{NSR} = 1$ mimics the noise level in the application on gene expression data from Section \ref{subsec-app-1}, where the estimated noise-to-signal ratios all lie between $0.6$ and $1$. The ratios $\textnormal{NSR} = 1.5$ and $\textnormal{NSR} = 2$ are used to investigate how the CluStErr method behaves when the noise level increases. We implement the CluStErr algorithm with $\alpha = 0.05$ and the difference-based estimators $\widehat{\sigma}^2_{\text{pc}}$ and $\widehat{\kappa}_{\text{pc}}$ from Section \ref{subsec-est-sigma}. For each of the three noise-to-signal ratios under consideration, we simulate $B = 1000$ samples and compute the estimate $\widehat{K}_0$ for each sample.

\begin{figure}[!t]
\centering
\includegraphics[scale = 0.8]{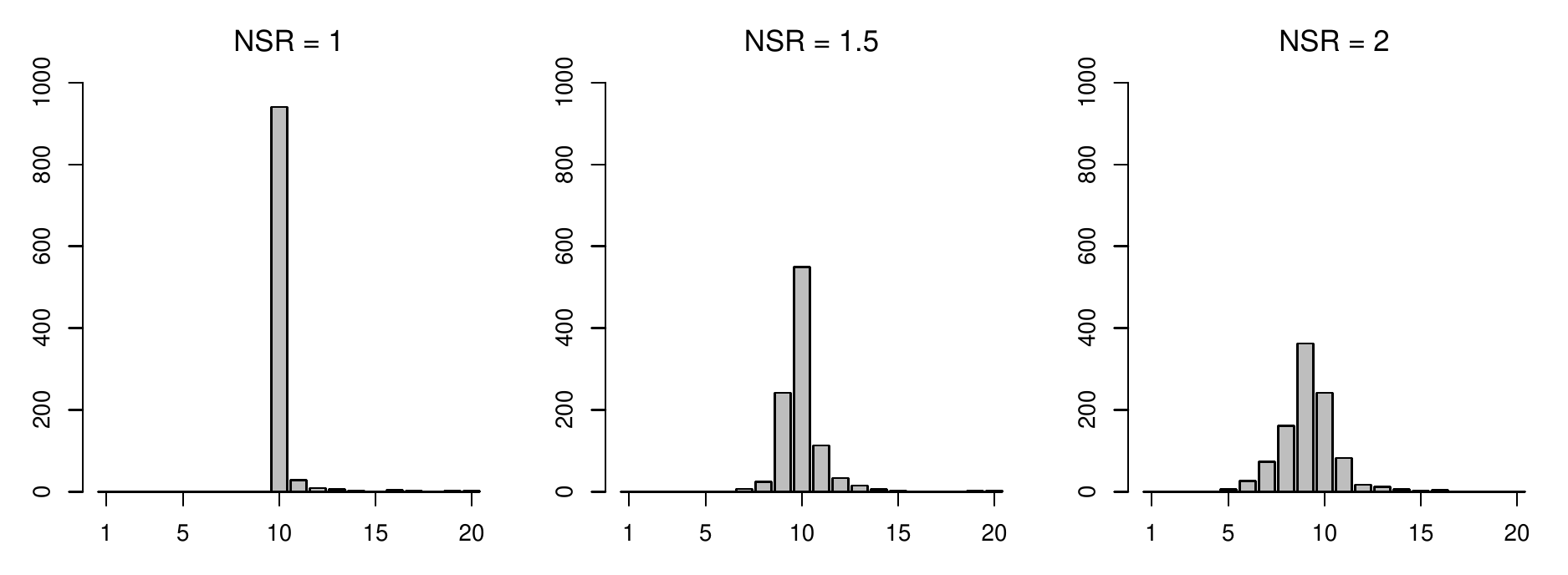} 
\caption{Histograms of the estimates $\widehat{K}_0$ in the three simulation scenarios with the noise-to-signal ratios $\textnormal{NSR} = 1$, $1.5$ and $2$.}\label{sim-fig1}

\vspace{-0.25cm}

\end{figure}

The simulation results are presented in Figure \ref{sim-fig1}. Each panel shows a histogram of the estimates $\widehat{K}_0$ for a specific noise-to-signal ratio. For the ratio level $\textnormal{NSR} = 1$, the CluStErr method produces very accurate results: About $95\%$ of the estimates are equal to the true value $K_0 = 10$ and most of the remaining estimates take the value $11$. For the ratio level $\textnormal{NSR} = 1.5$, the estimation results are also quite precise: Most of the estimates take a value between $9$ and $11$ with around $55\%$ of them being equal to the true value $K_0 = 10$. Only for the highest noise-to-signal ratio $\textnormal{NSR} = 2$, the estimation results are less accurate. In this case, the noise level in the data is too high for the method to produce precise results. As one can see, the estimates have a strong downward bias, which can be explained as follows: When there is too much noise in the data, the test procedure on which the estimator $\widehat{K}_0$ is based does not have enough power to detect the alternative $H_1: K < K_0$. As a result, our repeated test procedure stops too soon, thus underestimating the true number of clusters.

\begin{table}[!b]
\caption{Estimates of the probabilities $\pr(\widehat{K}_0 < K_0)$, $\pr(\widehat{K}_0 = K_0)$ and $\pr(\widehat{K}_0 > K_0)$ in the three simulation scenarios with the noise-to-signal ratios $\textnormal{NSR} = 1$, $1.5$ and $2$.}\label{sim-table1}
\scriptsize
\centering
\begin{tabular}{@{}cccc@{}} \toprule
 & $\textnormal{NSR} = 1$ & $\textnormal{NSR} = 1.5$  & $\textnormal{NSR} = 2$ \\ 
\cmidrule{1-4}
$\pr(\widehat{K}_0 < K_0)$ & 0.000 & 0.274 & 0.631 \\
$\pr(\widehat{K}_0 = K_0)$ & 0.941 & 0.549 & 0.242 \\
$\pr(\widehat{K}_0 > K_0)$ & 0.059 & 0.177 & 0.127 \\
\bottomrule
\end{tabular}
\end{table}

According to our theoretical results, the estimator $\widehat{K}_0$ allows for statistical error control in the following sense: It has the property that $\pr(\widehat{K}_0 > K_0) = \alpha + o(1)$ and $\pr(\widehat{K}_0 < K_0) = o(1)$, implying that $\pr(\widehat{K}_0 = K_0) = (1-\alpha) + o(1)$. Setting $\alpha$ to $0.05$, we should thus observe that $\widehat{K}_0$ equals $K_0 = 10$ in approximately $95\%$ of the simulations and overestimates $K_0$ in about $5\%$ of them. Table \ref{sim-table1} shows that this is indeed the case for the lowest noise-to-signal ratio $\textnormal{NSR} = 1$. In this situation, the probability $\pr(\widehat{K}_0 > K_0)$ of overestimating $K_0$ is around $5\%$, while the probability $\pr(\widehat{K}_0 < K_0)$ of underestimating $K_0$ is $0\%$, implying that $\pr(\widehat{K}_0 = K_0)$ is about $95\%$. For the two higher ratio levels $\textnormal{NSR} = 1.5$ and $2$, in contrast, the estimated values of the probabilities $\pr(\widehat{K}_0 < K_0)$, $\pr(\widehat{K}_0 = K_0)$ and $\pr(\widehat{K}_0 > K_0)$ do not accurately match the values predicted by the theory. This is due to the fact that the statistical error control of the CluStErr method is asymptotic in nature. Table \ref{sim-table2} illustrates this fact by reporting the estimated values of the probabilities $\pr(\widehat{K}_0 < K_0)$, $\pr(\widehat{K}_0 = K_0)$ and $\pr(\widehat{K}_0 > K_0)$ for the noise-to-signal ratio $\textnormal{NSR} = 1.5$ and various sample sizes $(n,p) = (1000,30)$, $(1500,40)$, $(2000,50)$, $(2500,60)$, $(3000,70)$. As one can clearly see, the estimated probabilities gradually approach the values predicted by the theory as the sample size increases.

\begin{table}[!t]
\caption{Estimates of the probabilities $\pr(\widehat{K}_0 < K_0)$, $\pr(\widehat{K}_0 = K_0)$ and $\pr(\widehat{K}_0 > K_0)$ in the simulation scenario with $\textnormal{NSR} = 1.5$ and five different sample sizes $(n,p)$.}\label{sim-table2}
\scriptsize
\centering

\begin{tabular}{@{}cccccc@{}} \toprule
$(n,p)$ & $(1000,30)$  & $(1500,40)$  & $(2000,50)$  & $(2500,50)$  & $(3000,60)$ \\ 
\cmidrule{1-6}
$\pr(\widehat{K}_0 < K_0)$ & 0.274 & 0.037 & 0.002 & 0.000 & 0.000 \\
$\pr(\widehat{K}_0 = K_0)$ & 0.549 & 0.795 & 0.893 & 0.918 & 0.956 \\
$\pr(\widehat{K}_0 > K_0)$ & 0.177 & 0.168 & 0.105 & 0.082 & 0.044 \\
\bottomrule
\end{tabular}

\end{table}


To summarize, our simulations on the finite sample behaviour of the CluStErr method indicate the following: 
(i) The method produces accurate estimates of $K_0$ as long as the noise level in the data is not too high. 
(ii) For sufficiently large sample sizes, it controls the probability of under- and overestimating the number of clusters $K_0$ as predicted by the theory. (iii) For smaller sample sizes, however, the error control is not fully accurate.

It is important to note that (iii) is not a big issue: Even in situations where the error control is not very precise, the CluStErr method may still produce accurate estimates of $K_0$. This is illustrated by our simulations. Inspecting the histogram of Figure \ref{sim-fig1} with $\textnormal{NSR} = 1.5$, for example, the estimated probability $\pr(\widehat{K}_0 = K_0)$ is seen to be only around $55\%$ rather than $95\%$. Nevertheless, most of the estimates take a value between $9$ and $11$. Hence, in most of the simulations, the CluStErr method yields a reasonable approximation to the true number of clusters. From a heuristic perspective, the CluStErr method can indeed be expected to produce satisfying estimation results even in smaller samples when the error control is not very precise. This becomes clear when regarding CluStErr as a thresholding procedure. For $K = 1,2,\ldots$, it checks whether the statistic $\widehat{\mathcal{H}}^{[K]}$ is below a certain threshold level $q$ and stops as soon as this is the case. For this approach to work, it is crucial to pick the threshold level $q$ appropriately. Our theoretical results suggest that the choice $q = q(\alpha)$ for a common significance level such as $\alpha = 0.05$ should be appropriate. Of course, this choice guarantees precise error control only for sufficiently large sample sizes. Nevertheless, in smaller samples, the threshold level $q = q(\alpha)$ can still be expected to be of the right order of magnitude, thus resulting in reasonable estimates of $K_0$. 

\vspace{10pt}

\noindent \textbf{Comparison of CluStErr with competing methods.} We now compare the CluStErr method to other criteria for selecting the number of clusters $K_0$, in particular to (i) the gap statistic \citep{gapstatistic}, (ii) the silhouette statistic \citep{rousseeuw1987}, and (iii) the Hartigan index \citep{hartigan1975}. As before, we set the sample size to $(n,p) = (1000,30)$ and consider the three noise-to-signal ratios $\textnormal{NSR} = 1$, $1.5$ and $2$. In addition to a ``balanced'' scenario with clusters of the same size $n/K_0$ each, we also consider an ``unbalanced'' scenario with clusters of sizes $1 + 18k$ for $k = 1,\ldots,K_0$. For each design, we simulate $B =100$ samples and compare the estimated cluster numbers obtained from the CluStErr method with those produced by the three competing algorithms.

The CluStErr estimates are computed as described in the first part of the simulation study. The three competing methods are implemented with a $k$-means algorithm as the underlying clustering method. To compute the values of the gap statistic, we employ the {\tt clusGap} function contained in the R package {\bf cluster} \citep{clusGap}. The number of clusters is estimated by the function {\tt maxSE} with the option {\tt Tibs2001SEmax}. We thus determine the number of clusters as suggested in \cite{gapstatistic}. To compute the silhouette and Hartigan statistics, we apply the R package {\bf NbClust} \citep{NbClust}.

\begin{table}[!t]
\caption{Results of the comparison study. The entries of the table give the numbers of simulations (out of a total of $100$) for which a certain estimate of $K_0$ is obtained. The first line in part (a) of the table, for example, has to be read as follows: The CluStErr estimate $\widehat{K}_0$ is equal to the true $K_0 = 10$ in $95$ out of $100$ simulations, and it is equal to $K =11, 12, 13$ in $2,1,2$ simulations, respectively.}\label{sim-table3} 
\scriptsize
\centering
\subfloat[balanced scenario]{\label{sim-table3a}
\begin{tabular}{@{}llccccccccccccccc@{}} 
\toprule
    &  & \multicolumn{15}{c}{Estimated number of clusters} \\
\cmidrule{3-17}
NSR & Method & 1 & 2 & 3 & 4 & 5 & 6 & 7 & 8 & 9 & 10 & 11 & 12 & 13 & 14 & $\ge 15$ \\
\cmidrule{1-17}
  1 & CluStErr   &  0 &   0 &   0 &   0 &   0 &   0 &   0 &   0 &   0 &  95 &   2 &   1 &   2 &   0 &   0 \\ 
    & Gap        &  3 &   5 &   0 &   0 &   1 &   3 &  11 &  25 &  36 &  16 &   0 &   0 &   0 &   0 &   0 \\ 
    & Silhouette &  0 &   0 &   0 &   0 &   0 &   0 &   0 &   5 &  19 &  43 &  20 &  11 &   1 &   0 &   1 \\ 
    & Hartigan   &  0 &   0 &   0 &   0 &   0 &   0 &   0 &   1 &   4 &  40 &  31 &  10 &   5 &   6 &   3 \\ 
1.5 & CluStErr   &  0 &   0 &   0 &   0 &   0 &   0 &   0 &   1 &  22 &  59 &  10 &   3 &   3 &   1 &   1 \\ 
    & Gap        & 13 &  26 &   0 &   0 &   0 &   0 &   2 &  13 &  17 &  29 &   0 &   0 &   0 &   0 &   0 \\ 
    & Silhouette &  0 &   0 &   0 &   0 &   0 &   0 &   0 &   2 &  12 &  54 &  23 &   8 &   1 &   0 &   0 \\ 
    & Hartigan   &  0 &   0 &   0 &   0 &   0 &   0 &   0 &   1 &   3 &  61 &  18 &   9 &   2 &   4 &   2 \\ 
  2 & CluStErr   &  0 &   0 &   0 &   0 &   1 &   2 &  10 &  13 &  32 &  31 &   7 &   2 &   1 &   0 &   1 \\ 
    & Gap        & 22 &  26 &   2 &   0 &   0 &   0 &   1 &   6 &   9 &  34 &   0 &   0 &   0 &   0 &   0 \\ 
    & Silhouette &  0 &   0 &   0 &   0 &   0 &   0 &   0 &   0 &   6 &  66 &  26 &   2 &   0 &   0 &   0 \\ 
    & Hartigan   &  0 &   0 &   0 &   0 &   0 &   0 &   0 &   0 &  10 &  68 &  18 &   2 &   2 &   0 &   0 \\ 
\bottomrule
\end{tabular}
}

\subfloat[unbalanced scenario]{\label{sim-table3b}
\begin{tabular}{@{}llccccccccccccccc@{}} 
\toprule
    &  & \multicolumn{15}{c}{Estimated number of clusters} \\
\cmidrule{3-17}
NSR & Method & 1 & 2 & 3 & 4 & 5 & 6 & 7 & 8 & 9 & 10 & 11 & 12 & 13 & 14 & $\ge 15$ \\
\cmidrule{1-17}
  1 & CluStErr   &  0 &   0 &   0 &   0 &   0 &   0 &   0 &   0 &   0 &  94 &   4 &   2 &   0 &   0 &   0 \\ 
    & Gap        &  0 &   1 &   4 &   9 &  17 &  11 &  21 &  20 &  16 &   1 &   0 &   0 &   0 &   0 &   0 \\ 
    & Silhouette &  0 &   0 &   0 &   0 &   0 &   0 &  10 &  26 &  37 &  19 &   6 &   2 &   0 &   0 &   0 \\ 
    & Hartigan   &  0 &   0 &   1 &   2 &   7 &   8 &  21 &  21 &  12 &  11 &   7 &   2 &   2 &   1 &   5 \\ 
1.5 & CluStErr   &  0 &   0 &   0 &   0 &   0 &   0 &   0 &   3 &  15 &  55 &  20 &   5 &   2 &   0 &   0 \\ 
    & Gap        &  2 &   4 &   4 &   3 &   8 &  15 &  21 &  21 &  22 &   0 &   0 &   0 &   0 &   0 &   0 \\ 
    & Silhouette &  0 &   0 &   0 &   0 &   0 &   1 &   7 &  20 &  47 &  24 &   1 &   0 &   0 &   0 &   0 \\ 
    & Hartigan   &  0 &   0 &   0 &   1 &  14 &  20 &  12 &  14 &  20 &   8 &   7 &   4 &   0 &   0 &   0 \\ 
  2 & CluStErr   &  0 &   0 &   0 &   0 &   0 &   0 &   3 &  17 &  25 &  31 &  15 &   3 &   6 &   0 &   0 \\ 
    & Gap        & 14 &   4 &   5 &   1 &   6 &   6 &  19 &  26 &  19 &   0 &   0 &   0 &   0 &   0 &   0 \\ 
    & Silhouette &  0 &   0 &   0 &   0 &   0 &   0 &   5 &  28 &  48 &  19 &   0 &   0 &   0 &   0 &   0 \\ 
    & Hartigan   &  0 &   0 &   1 &   4 &  15 &  15 &  15 &  18 &  17 &   7 &   5 &   0 &   0 &   2 &   1 \\ 
\bottomrule
\end{tabular}
}


\end{table}

\newpage

The results of the comparison study are presented in Table \ref{sim-table3}. Part (a) of the table provides the results for the balanced scenario with equal cluster sizes. As can be seen, the CluStErr method clearly outperforms its competitors in the setting with $\textnormal{NSR} = 1$. In the scenario with $\textnormal{NSR} = 1.5$, it also performs well in comparison to the other methods. Only for the highest noise-to-signal ratio $\textnormal{NSR} = 2$, it produces estimates of $K_0$ with a strong downward bias and is outperformed by the silhouette and Hartigan statistics. Part (b) of Table \ref{sim-table3} presents the results for the unbalanced scenario where the clusters strongly differ in size. In this scenario, all of the three competing methods substantially underestimate the number of clusters. The CluStErr method, in contrast, provides accurate estimates of $K_0$ in the two designs with $\textnormal{NSR} = 1$ and $\textnormal{NSR} = 1.5$. Only in the high-noise design with $\textnormal{NSR} = 2$, it produces estimates with a substantial downward bias, which nevertheless is much less pronounced than that of its competitors.

To summarize, the main findings of our comparison study are as follows: 
(i) The CluStErr method performs well in comparison to its competitors as long as the noise-to-signal ratio is not too high. It is however outperformed by its competitors in a balanced setting when the noise level is high. In Section \ref{sec-ext}, we discuss some modifications of the CluStErr method to improve its behaviour in the case of high noise. 
(ii) The CluStErr method is able to deal with both balanced and unbalanced cluster sizes, whereas its competitors perform less adequately in unbalanced settings.

The findings (i) and (ii) can heuristically be explained as follows: The CluStErr method is based on the test statistic $\widehat{\mathcal{H}}^{[K]} = \max_{1 \le i \le n} \widehat{\Delta}_i^{[K]}$, which is essentially the maximum over the residual sums of squares of the various individuals $i$. Its competitors, in contrast, are based on statistics which evaluate averages rather than maxima. Hartigan's rule, for instance, relies on a statistic which is essentially a scaled version of the ratio $\text{RSS}(K)/\text{RSS}(K+1)$, where $\text{RSS}(K)$ is defined as in \eqref{RSS-K} and denotes the average residual sum of squares for a partition with $K$ clusters. Averaging the residual sums of squares reduces the noise in the data much more strongly than taking the maximum. This is the reason why Hartigan's rule tends to perform better than the CluStErr method in a balanced setting with high noise. On the other hand, the average residual sum of squares hardly reacts to changes in the residual sums of squares of a few individuals that form a small cluster. Hence, small clusters are effectively ignored when taking the average of the residual sums of squares. This is the reason why Hartigan's statistic is not able to deal adequately with unbalanced settings. Taking the maximum of the residual sums of squares instead allows us to handle even highly unbalanced cluster sizes.

\section{Extensions}\label{sec-ext}

In this paper, we have developed an approach for estimating the number of clusters with statistical error control. We have derived a rigorous mathematical theory for a model with convex spherical clusters. This model is widely used in practice and is suitable for a large number of applications. Nevertheless, it of course has some limitations. In particular, it is not suitable for applications where the clusters have non-convex shapes. An interesting question is how to extend our ideas to the case of general, potentially non-convex clusters. Developing theory for this general case is a very challenging problem. We have made a first step into this direction by providing theory for the case of spherical clusters.

There are several ways to modify and extend our estimation methods in the model setting \eqref{model-eq1}--\eqref{model-eq2}. So far, we have based our methods on the maximum statistic $\widehat{\mathcal{H}}^{[K]} = \max_{1 \le i \le n} \widehat{\Delta}_i^{[K]}$. However, we are not bound to this choice. Our approach can be based on any test statistic $\widehat{\mathcal{H}}^{[K]}$ that fulfills the higher-order property \eqref{prop-stat}. The maximum statistic serves as a baseline which may be modified and improved in several directions. 
The building blocks of the maximum statistic are the individual statistics $\widehat{\Delta}_i^{[K]}$. Their stochastic behaviour has been analyzed in detail in Section \ref{subsec-est-stat}. Under the null hypothesis $H_0: K = K_0$, the statistics $\widehat{\Delta}_i^{[K]}$ are approximately independent and distributed as $(\chi_p^2 - p)/ \sqrt{2p}$ variables. Under the alternative $H_1: K < K_0$ in contrast, they have an explosive behaviour at least for some $i$. This difference in behaviour suggests to test $H_0$ by checking whether the hypothesis 
\[ H_{0,i}: \widehat{\Delta}_i^{[K]} \text{ has a } (\chi_p^2 - p)/ \sqrt{2p} \text{ distribution} \]
holds for all subjects $i = 1,\ldots,n$. We are thus faced with a multiple testing problem. A maximum statistic is a classical tool to tackle this problem. However, as is well known from the field of multiple testing, maximum statistics tend to be fairly conservative. When the noise level in the data is high, a test based on the maximum statistic $\widehat{\mathcal{H}}^{[K]} = \max_{1 \le i \le n} \widehat{\Delta}_i^{[K]}$ can thus be expected to have low power against the alternative $H_1: K < K_0$. As a consequence, the repeated test procedure on which the estimator $\widehat{K}_0$ is based tends to stop too soon, thus underestimating the true number of clusters. This is exactly what we have seen in the high-noise scenarios of the simulation study from Section \ref{subsec-app-2}. We now present two ways how to construct a statistic $\widehat{\mathcal{H}}^{[K]}$ with better power properties. 

\vspace{10pt}

\noindent \textbf{A blocked maximum statistic.} Let $w_0 = 0$ and define $w_k = \sum\nolimits_{r=1}^k \# \widehat{G}_r^{[K]}$ for $1 \le k \le K$. Moreover, write $\widehat{G}_k^{[K]} = \{ i_{w_{k-1}+1},\ldots,i_{w_k}\}$ with $i_{w_{k-1}+1} < \ldots < i_{w_k}$ for any $k$. To start with, we order the indices $\{ 1, \ldots, n \}$ clusterwise. In particular, we write them as $\{i_1,i_2,\ldots,i_n\}$, which yields the ordering  
\[ \overbrace{i_1 < \ldots < i_{w_1}}^{\widehat{G}_1^{[K]}} \quad \overbrace{i_{w_1+1} < \ldots < i_{w_2}}^{\widehat{G}_2^{[K]}} \quad \ldots\ldots \quad \overbrace{i_{w_{K-1}+1} < \ldots < i_{w_K}}^{\widehat{G}_K^{[K]}}. \]
We next partition the ordered indices into blocks 
\[ B_{\ell}^{[K]} = \big\{ i_{(\ell-1)N + 1},\ldots, i_{\ell N \wedge n} \big\} \quad \text{for } 1 \le \ell \le L, \]
where $N$ is the block length and $L = \lceil n/N \rceil$ is the number of blocks. With this notation at hand, we construct blockwise averages 
\begin{equation*}
\widehat{\Lambda}_{\ell}^{[K]} = \frac{1}{\sqrt{N}} \sum\limits_{i \in B_{\ell}^{[K]}} \widehat{\Delta}_i^{[K]} 
\end{equation*}
of the individual statistics $\widehat{\Delta}_i^{[K]}$ and aggregate them by taking their maximum, thus defining 
\begin{equation*}
\widehat{\mathcal{H}}^{[K]}_B = \max_{1 \le \ell \le L} \widehat{\Lambda}_\ell^{[K]}. 
\end{equation*}
In addition, we let $q_B(\alpha)$ be the $(1-\alpha)$-quantile of 
\begin{equation*}
\mathcal{H}_B = \max_{1 \le \ell \le L} \Lambda_\ell \quad \text{with} \quad \Lambda_\ell = \frac{1}{\sqrt{N}} \sum\limits_{i = (\ell-1)N + 1}^{\ell N} Z_i, 
\end{equation*}
where $Z_i$ are i.i.d.\ variables with the distribution $(\chi_p^2 - p)/\sqrt{2p}$. Note that this definition of $\widehat{\mathcal{H}}^{[K]}_B$ nests the maximum statistic $\widehat{\mathcal{H}}^{[K]} = \max_{1 \le i \le n} \widehat{\Delta}_i^{[K]}$ as a special case with the block length $N = 1$.

Under appropriate restrictions on the block length $N$, the estimators that result from applying the CluStErr method with the blocked statistic $\widehat{\mathcal{H}}^{[K]}_B$ can be shown to have the theoretical properties stated in Theorems \ref{theo1}--\ref{theo3}. More specifically, Theorems \ref{theo1}--\ref{theo3} can be shown to hold true for the blocked statistic $\widehat{\mathcal{H}}^{[K]}_B$ if the following two restrictions are satisfied: (i) $N/p^{1-\eta} = O(1)$ for some small $\eta > 0$, that is, the block length $N$ diverges more slowly than $p$. (ii) $\# G_k / n \rightarrow c_k > 0$ for all $k$, that is, the cluster sizes $\# G_k$ all grow at the same rate. Condition (ii) essentially rules out strongly differing cluster sizes. It is not surprising that we require such a restriction: To construct the blocked statistic $\widehat{\mathcal{H}}^{[K]}_B$, we average over the individual statistics $\widehat{\Delta}_i^{[K]}$. As already discussed in the context of the simulation study of Section \ref{subsec-app-2}, averaging has the effect that small clusters are effectively ignored. Hence, in contrast to the maximum statistic $\widehat{\mathcal{H}}^{[K]} = \max_{1 \le i \le n} \widehat{\Delta}_i^{[K]}$, the blocked statistic $\widehat{\mathcal{H}}^{[K]}_B$ with a large block size $N$ can be expected not to perform adequately when the cluster sizes are highly unbalanced.

In balanced settings, however, the blocked statistic $\widehat{\mathcal{H}}^{[K]}_B$ can be shown to have better power properties than the maximum statistic when the block size $N$ is chosen sufficiently large. To see this, we examine the behaviour of $\widehat{\mathcal{H}}^{[K]}_B$ for different block lengths $N$. Our heuristic discussion of the individual statistics $\widehat{\Delta}_i^{[K]}$ from Section \ref{subsec-est-stat} directly carries over to the blocked versions $\widehat{\Lambda}_{\ell}^{[K]}$: With the help of \eqref{delta-K0}, it is easy to see that
\[ \pr\Big( \widehat{\mathcal{H}}_B^{[K_0]} \le q_B(\alpha) \Big) \approx (1 - \alpha). \]
Moreover, \eqref{delta-K} together with some additional arguments suggests that $\widehat{\mathcal{H}}_B^{[K]}$ has an explosive behaviour for $K < K_0$. Specifically, under our technical conditions from Section \ref{subsec-asym-ass} and the two additional restrictions (i) and (ii) from above, we can show that  
\begin{equation}\label{rate-H-hat-B}
\widehat{\mathcal{H}}_B^{[K]} \ge c \sqrt{Np} \quad \text{ for some } c > 0 \text{ with prob.\ tending to } 1. 
\end{equation} 
As the quantile $q_B(\alpha)$ grows at the slower rate $\sqrt{\log L} \ (\le \sqrt{\log n})$, we can conclude that    
\[ \pr\Big( \widehat{\mathcal{H}}_B^{[K]} \le q_B(\alpha) \Big) = o(1) \]
for $K < K_0$. As a result, $\widehat{\mathcal{H}}_B^{[K]}$ should satisfy the higher-order property \eqref{prop-stat}. Moreover, according to \eqref{rate-H-hat-B}, the statistic $\widehat{\mathcal{H}}_B^{[K]}$ explodes at the rate $\sqrt{Np}$ for $K < K_0$. Hence, the faster the block size $N$ grows, the faster $\widehat{\mathcal{H}}_B^{[K]}$ diverges to infinity. Put differently, the larger $N$, the more power we have to detect that $K < K_0$. This suggests to select $N$ as large as possible. According to restriction (i) from above, we may choose any $N$ with $N/p^{1-\eta} = O(1)$ for some small $\eta > 0$. Ideally, we would thus like to pick $N$ so large that it grows at the same rate as $p^{1 - \eta}$. In practice, we neglect the small constant $\eta > 0$ and set $N = p$ as a simple rule of thumb.

\begin{table}[!t]
\caption{Simulation results for the blocked maximum statistic.}\label{sim-table4} 
\scriptsize
\centering

\begin{tabular}{@{}lcccccccccccccc@{}} 
\toprule
    & \multicolumn{14}{c}{Estimated number of clusters} \\
\cmidrule{2-15}
NSR & 1 & 2 & 3 & 4 & 5 & 6 & 7 & 8 & 9 & 10 & 11 & 12 & 13 & 14 \\
\cmidrule{1-15}
1   & 0 & 0 & 0 & 0 & 0 & 0 & 0 & 0 & 0 & 98 & 1 & 1 & 0 & 0 \\
1.5 & 0 & 0 & 0 & 0 & 0 & 0 & 0 & 0 & 0 & 79 & 12 & 5 & 2 & 2 \\
2   & 0 & 0 & 0 & 0 & 0 & 0 & 0 & 0 & 0 & 73 & 21 & 2 & 3 & 1 \\
\bottomrule
\end{tabular}

\end{table}

According to the heuristic arguments from above, the blocked maximum statistic $\widehat{\mathcal{H}}^{[K]}_B$ with block length $N = p$ should be more powerful than the maximum statistic $\widehat{\mathcal{H}}^{[K]} = \max_{1 \le i \le n} \widehat{\Delta}_i^{[K]}$ in settings with balanced cluster sizes. We examine this claim with the help of some simulations. To do so, we return to the balanced scenario of the comparison study in Section \ref{subsec-app-2}. For each of the data samples that were simulated for this scenario, we compute the CluStErr estimate of $K_0$ based on the blocked statistic $\widehat{\mathcal{H}}^{[K]}_B$ with $N=p$. Table \ref{sim-table4} presents the results. It shows that the blocked CluStErr method yields accurate estimates of $K_0$ for all three noise-to-signal ratios. Comparing the results to those in Table \ref{sim-table3}(a), the blocked method can be seen to perform very well in comparison to the other procedures even in the high-noise setting with $\textnormal{NSR} = 2$. This clearly shows the gain in power induced by the block structure of the statistic. 

\vspace{10pt}

\noindent \textbf{An FDR-based statistic.} There are several approaches in the literature how to construct multiple testing procedures that have better power properties than the classical maximum statistic. Prominent examples are methods that control the false discovery rate (FDR) or the higher criticism procedure by \cite{Donoho2004}. We may try to exploit ideas from these approaches to construct a more powerful statistic $\widehat{\mathcal{H}}^{[K]}$. As an example, we set up a test statistic which uses ideas from the FDR literature: Rather than only taking into account the maximum of the statistics $\widehat{\Delta}_i^{[K]}$, we may try to exploit the information in all of the ordered statistics $\widehat{\Delta}_{(1)}^{[K]} \ge \ldots \ge \widehat{\Delta}_{(n)}^{[K]}$. In particular, following \cite{Simes1986} and \cite{Benjamini1995}, we may set up our test for a given number of clusters $K$ as follows: Reject $H_0: K = K_0$ if 
\[ \widehat{\Delta}_{(i)}^{[K]} > q_{\chi} \Big( \frac{i}{n} \, \alpha \Big) \quad \text{for some } i \in \{1,\ldots,n\}, \]
where $q_{\chi}(\beta)$ is the $(1 - \beta)$-quantile of the distribution $(\chi_p^2 - p)/\sqrt{2p}$. This procedure can be rephrased as follows: Define the statistic
\[ \widehat{\mathcal{H}}^{[K]}_{\text{FDR}} = \max_{1 \le i \le n} \frac{\widehat{\Delta}_{(i)}^{[K]}}{q_{\chi}(i \alpha /n)} \]
and reject $H_0$ if $\widehat{\mathcal{H}}^{[K]}_{\text{FDR}} > 1$. Instead of taking the maximum over the original statistics $\widehat{\Delta}_i^{[K]}$, we thus take the maximum over rescaled versions of the ordered statistics $\widehat{\Delta}_{(i)}^{[K]}$. Developing theory for the FDR-type statistic $\widehat{\mathcal{H}}^{[K]}_{\text{FDR}}$ is a very interesting topic which is however far from trivial.

\bibliographystyle{ims}
{\small
\setlength{\bibsep}{0.45em}
\bibliography{bibliography}}

\newpage
\vspace*{15pt}

\begin{center}
{\LARGE \textbf{Supplement}}
\end{center}
\vspace{25pt}

\noindent In this supplement, we provide the proofs that are omitted in the paper. In particular, we derive Theorems \ref{theo1}--\ref{theo3} from Section \ref{sec-asym}. Throughout the supplement, we use the symbol $C$ to denote a universal real constant which may take a different value on each occurrence.

\renewcommand{\baselinestretch}{1.2}\normalsize
\numberwithin{equation}{section}
\def\theequation{S.\arabic{equation}}
\setcounter{equation}{0}
\def\thelemma{S.\arabic{lemma}}
\setcounter{lemma}{0}
\allowdisplaybreaks[2]

\section*{Auxiliary results}

In the proofs of Theorems \ref{theo1}--\ref{theo3}, we frequently make use of the following uniform convergence result. 
\begin{lemma}\label{lemmaA1}
Let $\mathcal{Z}_s = \{Z_{st}: 1 \le t \le T\}$ be sequences of real-valued random variables for $1 \le s \le S$ with the following properties: (i) for each $s$, the random variables in $\mathcal{Z}_s$ are independent of each other, and (ii) $\ex[Z_{st}] = 0$ and $\ex[|Z_{st}|^\phi] \le C < \infty$ for some $\phi > 2$ and $C > 0$ that depend neither on $s$ nor on $t$. Suppose that $S = T^q$ with $0 \le q < \phi/2 - 1$. Then
\[ \pr \Big( \max_{1 \le s \le S} \Big| \frac{1}{\sqrt{T}} \sum\limits_{t=1}^T Z_{st} \Big| > T^\eta \Big) = o(1), \]
where the constant $\eta > 0$ can be chosen as small as desired. 
\end{lemma}

\noindent \textbf{Proof of Lemma \ref{lemmaA1}.} Define $\tau_{S,T} = (ST)^{1/\{(2 + \delta)(q+1)\}}$ with some sufficiently small $\delta > 0$. In particular, let $\delta > 0$ be so small that $(2+\delta)(q+1) < \phi$. Moreover, set 
\begin{align*}
Z_{st}^{\le} & = Z_{st} \ind (|Z_{st}| \le \tau_{S,T}) - \ex\big[ Z_{st} \ind (|Z_{st}| \le \tau_{S,T}) \big] \\
Z_{st}^{>} & = Z_{st} \ind (|Z_{st}| > \tau_{S,T}) - \ex\big[ Z_{st} \ind (|Z_{st}| > \tau_{S,T}) \big]
\end{align*}
and write
\[ \frac{1}{\sqrt{T}} \sum\limits_{t=1}^T Z_{st} = \frac{1}{\sqrt{T}} \sum\limits_{t=1}^T Z_{st}^{\le} + \frac{1}{\sqrt{T}} \sum\limits_{t=1}^T Z_{st}^>. \]
In what follows, we show that 
\begin{align} 
 & \pr \Big(\max_{1 \le s \le S} \Big| \frac{1}{\sqrt{T}} \sum\limits_{t=1}^T Z_{st}^{>} \Big| > C T^{\eta} \Big) = o(1) \label{lemmaA1-claim1} \\
 & \pr \Big(\max_{1 \le s \le S} \Big| \frac{1}{\sqrt{T}} \sum\limits_{t=1}^T Z_{st}^{\le} \Big| > C T^{\eta} \Big) = o(1) \label{lemmaA1-claim2} 
\end{align}
for any fixed constant $C > 0$. Combining \eqref{lemmaA1-claim1} and \eqref{lemmaA1-claim2} immediately yields the statement of Lemma \ref{lemmaA1}.

We start with the proof of \eqref{lemmaA1-claim1}: It holds that   
\[ \pr \Big(\max_{1 \le s \le S} \Big| \frac{1}{\sqrt{T}} \sum\limits_{t=1}^T Z_{st}^{>} \Big| > C T^{\eta} \Big) \le Q_1^> + Q_2^>, \]
where 
\begin{align*}
Q_1^> 
 := & \, \sum\limits_{s=1}^S \pr \Big( \frac{1}{\sqrt{T}} \sum\limits_{t=1}^T |Z_{st}| \ind (|Z_{st}| > \tau_{S,T}) > \frac{C}{2} T^\eta \Big) \\
 \le & \, \sum\limits_{s=1}^S \pr \Big( |Z_{st}| > \tau_{S,T} \text{ for some } 1 \le t \le T \Big) \\
 \le & \, \sum\limits_{s=1}^S \sum\limits_{t=1}^T \pr \big( |Z_{st}| > \tau_{S,T} \big) \le \sum\limits_{s=1}^S \sum\limits_{t=1}^T \ex \Big[ \frac{|Z_{st}|^\phi}{\tau_{S,T}^\phi} \Big] \\[0.1cm]
 \le & \, \frac{CST}{\tau_{S,T}^\phi} = o(1)
\end{align*}
and
\[ Q_2^> := \sum\limits_{s=1}^S \pr \Big( \frac{1}{\sqrt{T}} \sum\limits_{t=1}^T \ex \big[ |Z_{st}| \ind (|Z_{st}| > \tau_{S,T}) \big] > \frac{C}{2} T^\eta \Big) = 0 \]
for $S$ and $T$ sufficiently large, since 
\begin{align*} 
 & \frac{1}{\sqrt{T}} \sum\limits_{t=1}^T \ex \big[ |Z_{st}| \ind (|Z_{st}| > \tau_{S,T}) \big] \\
 & \le \frac{1}{\sqrt{T}} \sum\limits_{t=1}^T \ex \Big[ \frac{|Z_{st}|^\phi}{\tau_{S,T}^{\phi-1}} \ind (|Z_{st}| > \tau_{S,T}) \Big] \\ 
 & \le \frac{C \sqrt{T}}{\tau_{S,T}^{\phi-1}} = o(T^\eta).
\end{align*}
This yields \eqref{lemmaA1-claim1}.

We next turn to the proof of \eqref{lemmaA1-claim2}: We apply the crude bound 
\[ \pr \Big(\max_{1 \le s \le S} \Big| \frac{1}{\sqrt{T}} \sum\limits_{t=1}^T Z_{st}^{\le} \Big| > C T^{\eta} \Big) \le \sum\limits_{s=1}^S \pr \Big( \Big| \frac{1}{\sqrt{T}} \sum\limits_{t=1}^T Z_{st}^{\le} \Big| > C T^{\eta} \Big) \]
and show that for any $1 \le s \le S$, 
\begin{equation}\label{lemmaA1-claim2-a}
\pr \Big( \Big| \frac{1}{\sqrt{T}} \sum\limits_{t=1}^T Z_{st}^{\le} \Big| > C T^{\eta} \Big) \le C_0 T^{-\rho},
\end{equation}
where $C_0$ is a fixed constant and $\rho > 0$ can be chosen as large as desired by picking $\eta$ slightly larger than $1/2 - 1/(2+\delta)$. Since $S = O(T^q)$, this immediately implies \eqref{lemmaA1-claim2}. To prove \eqref{lemmaA1-claim2-a}, we make use of the following facts:
\begin{enumerate}[label=(\roman*),leftmargin=0.85cm]
\item For a random variable $Z$ and $\lambda > 0$, Markov's inequality says that 
\[ \pr \big( \pm Z > \delta \big) \le \frac{\ex \exp(\pm \lambda Z)}{\exp(\lambda \delta)}. \]
\item Since $|Z_{st}^{\le} / \sqrt{T}| \le 2 \tau_{S,T} / \sqrt{T}$, it holds that $\lambda_{S,T} |Z_{st}^{\le} / \sqrt{T}| \le 1 /2$, where we set $\lambda_{S,T} = \sqrt{T} / (4 \tau_{S,T})$. As $\exp(x) \le 1 + x + x^2$ for $|x| \le 1/2$, this implies that 
\[ 
\ex \Big[ \exp \Big( \pm \lambda_{S,T} \frac{Z_{st}^{\le}}{\sqrt{T}} \Big) \Big] \le 1 + \frac{\lambda_{S,T}^2}{T} \ex \big[ (Z_{st}^{\le})^2 \big] \le \exp \Big( \frac{\lambda_{S,T}^2}{T} \ex \big[ (Z_{st}^{\le})^2 \big] \Big). \]
\item By definition of $\lambda_{S,T}$, it holds that 
\[ \lambda_{S,T} = \frac{\sqrt{T}}{4 (ST)^{\frac{1}{(2 + \delta)(q+1)}}} = \frac{\sqrt{T}}{4 (T^{q+1})^{\frac{1}{(2 + \delta)(q+1)}}} = \frac{T^{\frac{1}{2} - \frac{1}{2+\delta}}}{4}. \] 
\end{enumerate}
Using (i)--(iii) and writing $\ex (Z_{st}^{\le})^2 \le C_Z < \infty$, we obtain that  
\begin{align*}
 & \pr \Big(\Big| \frac{1}{\sqrt{T}} \sum\limits_{t=1}^T Z_{st}^{\le} \Big| > C T^{\eta} \Big) \\
 & \le \pr \Big(\frac{1}{\sqrt{T}} \sum\limits_{t=1}^T Z_{st}^{\le} > C T^{\eta} \Big) + \pr \Big(-\frac{1}{\sqrt{T}} \sum\limits_{t=1}^T Z_{st}^{\le} > C T^{\eta} \Big) \\
 & \le \exp \big( -\lambda_{S,T} CT^\eta \big) \left\{ \ex \Big[ \exp \Big( \lambda_{S,T} \sum\limits_{t=1}^T \frac{Z_{st}^{\le}}{\sqrt{T}} \Big) \Big] + \ex \Big[ \exp \Big( -\lambda_{S,T} \sum\limits_{t=1}^T \frac{Z_{st}^{\le}}{\sqrt{T}} \Big) \Big] \right\} \\
 & = \exp \big( -\lambda_{S,T} CT^\eta \big) \left\{ \prod\limits_{t=1}^T \ex \Big[ \exp \Big( \lambda_{S,T} \frac{Z_{st}^{\le}}{\sqrt{T}} \Big) \Big] + \prod\limits_{t=1}^T \ex \Big[ \exp \Big( -\lambda_{S,T} \frac{Z_{st}^{\le}}{\sqrt{T}} \Big) \Big] \right\} \\
 & \le 2 \exp \big( -\lambda_{S,T} CT^\eta \big) \prod\limits_{t=1}^T \exp \Big( \frac{\lambda_{S,T}^2}{T} \ex \big[ (Z_{st}^{\le})^2 \big] \Big) \\[0.1cm]
 & = 2 \exp \big( C_Z \lambda_{S,T}^2 - C \lambda_{S,T} T^\eta \big) \\[0.25cm]
 & = 2 \exp \Big( \frac{C_Z}{16} \big( T^{\frac{1}{2} - \frac{1}{2+\delta}} \big)^2 - \frac{C}{4} T^{\frac{1}{2} - \frac{1}{2+\delta}} \, T^\eta \Big) \\
 & \le C_0 T^{-\rho}, 
\end{align*}
where $\rho > 0$ can be chosen arbitrarily large if we pick $\eta$ slightly larger than $1/2 - 1/(2+\delta)$. \qed

\section*{Proof of Theorem \ref{theo1}}

We first prove that 
\begin{equation}\label{theo1-res1}
\pr \Big( \widehat{\mathcal{H}}^{[K_0]} \le q(\alpha) \Big) = (1-\alpha) + o(1). 
\end{equation}
To do so, we derive a stochastic expansion of the individual statistics $\widehat{\Delta}_i^{[K_0]}$. 
\begin{lemma}\label{lemmaA2}
It holds that 
\[ \widehat{\Delta}_i^{[K_0]} = \Delta_i^{[K_0]} + R_i^{[K_0]}, \]
where 
\[ \Delta_i^{[K_0]} = \frac{1}{\sqrt{p}} \sum\limits_{j=1}^p \Big\{ \frac{\varepsilon_{ij}^2}{\sigma^2} - 1 \Big\} \big/ \kappa \]
and the remainder $R_i^{[K_0]}$ has the property that 
\begin{equation}\label{theo1-rem-prop1}
\pr \Big( \max_{1 \le i \le n} \big| R_i^{[K_0]} \big| > p ^{-\xi} \Big) = o(1) 
\end{equation}
for some $\xi > 0$. 
\end{lemma}
\noindent The proof of Lemma \ref{lemmaA2} as well as those of the subsequent Lemmas \ref{lemmaA3}--\ref{lemmaA5} are postponed until the proof of Theorem \ref{theo1} is complete. With the help of Lemma \ref{lemmaA2}, we can bound the probability of interest
\[ P_{\alpha} := \pr \Big( \widehat{\mathcal{H}}^{[K_0]} \le q(\alpha) \Big) = \pr \Big( \max_{1 \le i \le n} \widehat{\Delta}_i^{[K_0]} \le q(\alpha) \Big) \] 
as follows: Since 
\[ \max_{1 \le i \le n} \widehat{\Delta}_i^{[K_0]} 
\begin{cases}
\le \max_{1 \le i \le n} \Delta_i^{[K_0]} + \max_{1 \le i \le n}  |R_i^{[K_0]}| \\
\ge \max_{1 \le i \le n} \Delta_i^{[K_0]} - \max_{1 \le i \le n}  |R_i^{[K_0]}|,
\end{cases}
\]
it holds that 
\[ P_{\alpha}^< \le P_{\alpha} \le P_{\alpha}^>, \]
where
\begin{align*}
P_{\alpha}^< & = \pr \Big( \max_{1 \le i \le n} \Delta_i^{[K_0]} \le q(\alpha) - \max_{1 \le i \le n} |R_i^{[K_0]}| \Big) \\
P_{\alpha}^> & = \pr \Big( \max_{1 \le i \le n} \Delta_i^{[K_0]} \le q(\alpha) + \max_{1 \le i \le n} |R_i^{[K_0]}| \Big). 
\end{align*}
As the remainder $R_i^{[K_0]}$ has the property \eqref{theo1-rem-prop1}, we further obtain that 
\begin{equation}\label{theo1-bound1}
P_{\alpha}^{\ll} +o(1) \le P_{\alpha} \le P_{\alpha}^{\gg} + o(1), 
\end{equation}
where
\begin{align*}
P_{\alpha}^{\ll} & = \pr \Big( \max_{1 \le i \le n} \Delta_i^{[K_0]} \le q(\alpha) - p^{-\xi} \Big) \\
P_{\alpha}^{\gg} & = \pr \Big( \max_{1 \le i \le n} \Delta_i^{[K_0]} \le q(\alpha) + p^{-\xi} \Big). 
\end{align*}
With the help of strong approximation theory, we can derive the following result on the asymptotic behaviour of the probabilities $P_{\alpha}^{\ll}$ and $P_{\alpha}^{\gg}$.  
\begin{lemma}\label{lemmaA3}
It holds that 
\begin{align*}
P_{\alpha}^{\ll} & = (1 - \alpha) + o(1) \\
P_{\alpha}^{\gg} & = (1 - \alpha) + o(1). 
\end{align*}
\end{lemma}
\noindent Together with \eqref{theo1-bound1}, this immediately yields that $P_{\alpha} = (1 - \alpha) + o(1)$, thus completing the proof of \eqref{theo1-res1}.

We next show that for any $K < K_0$,
\begin{equation}\label{theo1-res2}
\pr \Big( \widehat{\mathcal{H}}^{[K]} \le q(\alpha) \Big) = o(1). 
\end{equation}
Consider a fixed $K < K_0$ and let $S \in \{ \widehat{G}_k^{[K]}: 1 \le k \le K \}$ be any cluster with the following property: 
\begin{equation}\label{theo1-res2-property}
\begin{minipage}[c][1cm][c]{12.9cm}
$\# S \ge \underline{n} := \min_{1 \le k \le K_0} \# G_k$, and $S$ contains elements from at least two different classes $G_{k_1}$ and $G_{k_2}$. 
\end{minipage}
\end{equation}
It is not difficult to see that a cluster with the property \eqref{theo1-res2-property} must always exist under our conditions. By $\mathscr{C} \subseteq \{ \widehat{G}_k^{[K]}: 1 \le k \le K \}$, we denote the collection of clusters that have the property \eqref{theo1-res2-property}. With this notation at hand, we can derive the following stochastic expansion of the individual statistics $\widehat{\Delta}_i^{[K]}$. 
\begin{lemma}\label{lemmaA4}
For any $i \in S$ and $S \in \mathscr{C}$, it holds that 
\[ \widehat{\Delta}_i^{[K]} = \frac{1}{\kappa \sigma^2 \sqrt{p}} \sum\limits_{j=1}^p d_{ij}^2 + R_i^{[K]}, \]
where $d_{ij} = \mu_{ij} - (\#S)^{-1} \sum\nolimits_{i^\prime \in S} \mu_{i^\prime j}$ and the remainder $R_i^{[K]}$ has the property that 
\begin{equation}\label{theo1-rem-prop2}
\pr \Big( \max_{S \in \mathscr{C}} \max_{i \in S} \big| R_i^{[K]} \big| > p^{\frac{1}{2} - \xi} \Big) = o(1)
\end{equation}
for some small $\xi > 0$. 
\end{lemma}
\noindent Using \eqref{theo1-rem-prop2} and the fact that 
\[ \max_{S \in \mathscr{C}} \max_{i \in S} \widehat{\Delta}_i^{[K]} \ge  \max_{S \in \mathscr{C}} \max_{i \in S} \Big\{ \frac{1}{\kappa \sigma^2 \sqrt{p}} \sum\limits_{j=1}^p d_{ij}^2 \Big\} - \max_{S \in \mathscr{C}} \max_{i \in S} \big| R_i^{[K]} \big|, \] 
we obtain that 
\begin{align}
\pr \Big( \widehat{\mathcal{H}}^{[K]} \le q(\alpha) \Big)
 & = \pr \Big( \max_{1 \le i \le n} \widehat{\Delta}_i^{[K]} \le q(\alpha) \Big) \nonumber \\
 & \le \pr \Big( \max_{S \in \mathscr{C}} \max_{i \in S} \widehat{\Delta}_i^{[K]} \le q(\alpha) \Big) \nonumber \\
 & \le \pr \Big( \max_{S \in \mathscr{C}} \max_{i \in S} \Big\{ \frac{1}{\kappa \sigma^2 \sqrt{p}} \sum\limits_{j=1}^p d_{ij}^2 \Big\} - \max_{S \in \mathscr{C}} \max_{i \in S} \big| R_i^{[K]} \big| \le q(\alpha) \Big) \nonumber \\
 & \le \pr \Big( \max_{S \in \mathscr{C}} \max_{i \in S} \Big\{ \frac{1}{\kappa \sigma^2 \sqrt{p}} \sum\limits_{j=1}^p d_{ij}^2 \Big\} \le q(\alpha) + p^{\frac{1}{2} - \xi} \Big) + o(1). \label{theo1-res2-bounding}
\end{align}
The arguments from the proof of Lemma \ref{lemmaA3}, in particular \eqref{lemmaA3-res2}, imply that $q(\alpha) \le C \sqrt{\log n}$ for some fixed constant $C > 0$ and sufficiently large $n$. Moreover, we can prove the following result.
\begin{lemma}\label{lemmaA5}
It holds that 
\[ \max_{S \in \mathscr{C}} \max_{i \in S} \Big\{ \frac{1}{\sqrt{p}} \sum\limits_{j=1}^p d_{ij}^2 \Big\} \ge c \sqrt{p} \]
for some fixed constant $c > 0$. 
\end{lemma}
\noindent Since $q(\alpha) \le C \sqrt{\log n}$ and $\sqrt{\log n} / \sqrt{p} = o(1)$ by \ref{A3}, Lemma \ref{lemmaA5} allows us to infer that 
\[ \pr \Big( \max_{S \in \mathscr{C}} \max_{i \in S} \Big\{ \frac{1}{\kappa \sigma^2 \sqrt{p}} \sum\limits_{j=1}^p d_{ij}^2 \Big\} \le q(\alpha) + p^{\frac{1}{2} - \xi} \Big) = o(1). \]
Together with \eqref{theo1-res2-bounding}, this yields that $\pr( \widehat{\mathcal{H}}^{[K]} \le q(\alpha) ) = o(1)$. \qed

\vspace{12pt}

\noindent \textbf{Proof of Lemma \ref{lemmaA2}.} Let $n_k = \# G_k$ and write $\overline{\varepsilon}_i = p^{-1} \sum\nolimits_{j=1}^p \varepsilon_{ij}$ along with  $\overline{\mu}_i = p^{-1} \sum\nolimits_{j=1}^p \mu_{ij}$. Since 
\[ \pr \Big( \big\{ \widehat{G}_k^{[K_0]}: 1 \le k \le K_0 \big\} = \big\{ G_k: 1 \le k \le K_0 \big\} \Big) \rightarrow 1 \]
by \eqref{prop-clusters}, we can ignore the estimation error in the clusters $\widehat{G}_k^{[K_0]}$ and replace them by the true classes $G_k$. For $i \in G_k$, we thus get
\[ \widehat{\Delta}_i^{[K_0]} = \Delta_i^{[K_0]} + R_{i,A}^{[K_0]} + R_{i,B}^{[K_0]} - R_{i,C}^{[K_0]} + R_{i,D}^{[K_0]}, \]
where
\begin{align*}
R_{i,A}^{[K_0]} & = \Big( \frac{1}{\widehat{\kappa}} - \frac{1}{\kappa} \Big) \, \frac{1}{\sqrt{p}} \sum\limits_{j=1}^p \Big\{ \frac{\varepsilon_{ij}^2}{\sigma^2} - 1 \Big\} \\
R_{i,B}^{[K_0]} & = \frac{1}{\widehat{\kappa}} \Big( \frac{1}{\widehat{\sigma}^2} - \frac{1}{\sigma^2} \Big) \, \frac{1}{\sqrt{p}} \sum\limits_{j=1}^p \varepsilon_{ij}^2 \\
R_{i,C}^{[K_0]} & = \Big(\frac{2}{\widehat{\kappa} \widehat{\sigma}^2}\Big) \, \frac{1}{\sqrt{p}} \sum\limits_{j=1}^p  \varepsilon_{ij} \Big\{ \overline{\varepsilon}_i + \frac{1}{n_k} \sum\limits_{i^\prime \in G_k} \big( \varepsilon_{i^\prime j} - \overline{\varepsilon}_{i^\prime} \big) \Big\} \\ 
R_{i,D}^{[K_0]} & = \Big(\frac{1}{\widehat{\kappa} \widehat{\sigma}^2}\Big) \, \frac{1}{\sqrt{p}} \sum\limits_{j=1}^p \Big\{ \overline{\varepsilon}_i + \frac{1}{n_k} \sum\limits_{i^\prime \in G_k} \big( \varepsilon_{i^\prime j} - \overline{\varepsilon}_{i^\prime} \big) \Big\}^2.
\end{align*}
We now show that $\max_{i \in G_k} |R_{i,\ell}^{[K_0]}| = o_p(p^{-\xi})$ for any $k$ and $\ell = A,\ldots,D$. This implies that $\max_{1 \le i \le n} |R_{i,\ell}^{[K_0]}| = \max_{1 \le k \le K_0} \max_{i \in G_k} |R_{i,\ell}^{[K_0]}| = o_p(p^{-\xi})$ for $\ell = A,\ldots,D$, which in turn yields the statement of Lemma \ref{lemmaA2}. Throughout the proof, we use the symbol $\eta > 0$ to denote a sufficiently small constant which results from applying Lemma \ref{lemmaA1}.

By assumption, $\widehat{\sigma}^2 = \sigma^2 + O_p(p^{-(1/2 + \delta)})$ and $\widehat{\kappa} = \kappa + O_p(p^{-\delta})$ for some $\delta > 0$. Applying Lemma \ref{lemmaA1} and choosing $\xi > 0$ such that $\xi < \delta - \eta$, we obtain that 
\begin{align*}
\max_{i \in G_k} \big| R_{i,A}^{[K_0]} \big| 
 & \le \Big| \frac{1}{\widehat{\kappa}} - \frac{1}{\kappa} \Big| \, \max_{i \in G_k} \Big| \frac{1}{\sqrt{p}} \sum\limits_{j=1}^p \Big\{ \frac{\varepsilon_{ij}^2}{\sigma^2} - 1 \Big\} \Big| \\
 & = \Big| \frac{1}{\widehat{\kappa}} - \frac{1}{\kappa} \Big| O_p(p^\eta) = O_p(p^{-(\delta-\eta)}) = o_p(p^{-\xi})
\end{align*}
and 
\begin{align*}
\max_{i \in G_k} \big| R_{i,B}^{[K_0]} \big| 
 & \le \Big| \frac{1}{\widehat{\kappa}} \Big( \frac{1}{\widehat{\sigma}^2} - \frac{1}{\sigma^2} \Big) \Big| \, \Big\{ \max_{i \in G_k} \Big| \frac{1}{\sqrt{p}} \sum\limits_{j=1}^p \big(\varepsilon_{ij}^2 - \sigma^2 \big) \Big| + \sigma^2 \sqrt{p} \Big\} \\
 & = \Big| \frac{1}{\widehat{\kappa}} \Big( \frac{1}{\widehat{\sigma}^2} - \frac{1}{\sigma^2} \Big) \Big| \big\{ O_p(p^{\eta}) + \sigma^2 \sqrt{p} \big\} = o_p(p^{-\xi}).
\end{align*}

We next show that 
\begin{equation}\label{R-iC}
\max_{i \in G_k} \big| R_{i,C}^{[K_0]} \big| = o_p \big(p^{-\frac{1}{4}}\big). 
\end{equation}
To do so, we work with the decomposition $R_{i,C}^{[K_0]} = \{ 2 \widehat{\kappa}^{-1} \widehat{\sigma}^{-2} \} \{ R_{i,C,1}^{[K_0]} + R_{i,C,2}^{[K_0]} - R_{i,C,3}^{[K_0]} \}$, where 
\vspace{-1cm}

\begin{align*}
R_{i,C,1}^{[K_0]} & = \frac{1}{\sqrt{p}} \sum\limits_{j=1}^p \varepsilon_{ij} \overline{\varepsilon}_i \\
R_{i,C,2}^{[K_0]} & = \frac{1}{\sqrt{p}} \sum\limits_{j=1}^p \varepsilon_{ij} \Big( \frac{1}{n_k} \sum\limits_{i^\prime \in G_k} \varepsilon_{i^\prime j} \Big) \\
R_{i,C,3}^{[K_0]} & = \Big( \frac{1}{\sqrt{p}} \sum\limits_{j=1}^p \varepsilon_{ij} \Big) \Big( \frac{1}{n_k} \sum\limits_{i^\prime \in G_k} \overline{\varepsilon}_{i^\prime} \Big).
\end{align*}
With the help of Lemma \ref{lemmaA1}, we obtain that 
\begin{equation}\label{R-iC1}
\max_{i \in G_k} \big|R_{i,C,1}^{[K_0]}\big| \le \frac{1}{\sqrt{p}} \Big( \max_{i \in G_k} \Big| \frac{1}{\sqrt{p}} \sum\limits_{j=1}^p \varepsilon_{ij} \Big| \Big)^2 = O_p \Big(\frac{p^{2\eta}}{\sqrt{p}}\Big). 
\end{equation}
Moreover,
\begin{equation}\label{R-iC2}
\max_{i \in G_k} \big|R_{i,C,2}^{[K_0]}\big| = O_p\big(n_k^{-\frac{1}{4}}\big), 
\end{equation}
since
\[ R_{i,C,2}^{[K_0]} = \frac{1}{n_k \sqrt{p}} \sum\limits_{j=1}^p \big\{ \varepsilon_{ij}^2 - \sigma^2 \big\} +  \sigma^2 \frac{\sqrt{p}}{n_k} + \frac{1}{n_k} \sum\limits_{\substack{i^\prime \in G_k \\ i^\prime \ne i}} \frac{1}{\sqrt{p}} \sum\limits_{j=1}^p \varepsilon_{ij} \varepsilon_{i^\prime j}, \]
$p \ll n_k$ and 
\begin{align} 
\max_{i \in G_k} \Big| \frac{1}{n_k \sqrt{p}} \sum\limits_{j=1}^p \big\{ \varepsilon_{ij}^2 - \sigma^2 \big\} \Big| & = O_p\Big( \frac{p^{\eta}}{n_k} \Big) \label{R-iC2-1} \\
\max_{i \in G_k} \Big| \frac{1}{n_k} \sum\limits_{\substack{i^\prime \in G_k \\ i^\prime \ne i}} \frac{1}{\sqrt{p}} \sum\limits_{j=1}^p \varepsilon_{ij} \varepsilon_{i^\prime j} \Big| & = O_p\big( n_k^{-\frac{1}{4}} \big). \label{R-iC2-2}
\end{align}
\eqref{R-iC2-1} is an immediate consequence of Lemma \ref{lemmaA1}. \eqref{R-iC2-2} follows upon observing that for any constant $C_0 > 0$, 
\begin{align*}
 & \pr \Big( \max_{i \in G_k} \Big| \frac{1}{n_k} \sum\limits_{\substack{i^\prime \in G_k \\ i^\prime \ne i}} \frac{1}{\sqrt{p}} \sum\limits_{j=1}^p \varepsilon_{ij} \varepsilon_{i^\prime j} \Big| > \frac{C_0}{n_k^{1/4}} \Big) \\*
 & \le \sum\limits_{i \in G_k} \pr \Big( \Big| \frac{1}{n_k} \sum\limits_{\substack{i^\prime \in G_k \\ i^\prime \ne i}} \frac{1}{\sqrt{p}} \sum\limits_{j=1}^p \varepsilon_{ij} \varepsilon_{i^\prime j} \Big| > \frac{C_0}{n_k^{1/4}} \Big) \\
 & \le \sum\limits_{i \in G_k} \ex \Big\{ \frac{1}{n_k} \sum\limits_{\substack{i^\prime \in G_k \\ i^\prime \ne i}} \frac{1}{\sqrt{p}} \sum\limits_{j=1}^p \varepsilon_{ij} \varepsilon_{i^\prime j} \Big\}^4 \Big/ \Big\{ \frac{C_0}{n_k^{1/4}} \Big\}^4 \\ 
 & \le \sum\limits_{i \in G_k} \Big\{ \frac{1}{n_k^4 p^2} \sum\limits_{\substack{i_1^\prime,\ldots,i_4^\prime \in G_k \\ i_1^\prime,\ldots,i_4^\prime \ne i}} \sum\limits_{j_1,\ldots,j_4 = 1}^p \ex \big[ \varepsilon_{i \, j_1} \ldots \varepsilon_{i \, j_4} \varepsilon_{i_1^\prime j_1} \ldots \varepsilon_{i_4^\prime j_4} \big] \Big\}  \Big/ \Big\{ \frac{C_0}{n_k^{1/4}} \Big\}^4 \\
 & \le \frac{C}{C_0^4},  
\end{align*}
the last inequality resulting from the fact that the mean $\ex [ \varepsilon_{i \, j_1} \ldots \varepsilon_{i \, j_4} \varepsilon_{i_1^\prime j_1} \ldots \varepsilon_{i_4^\prime j_4} ]$ can only be non-zero if some of the index pairs $(i_\ell^\prime,j_\ell)$ for $\ell =1,\ldots,4$ are identical. Finally, with the help of Lemma \ref{lemmaA1}, we get that
\begin{equation}\label{R-iC3}
\max_{i \in G_k} \big|R_{i,C,3}^{[K_0]}\big| \le \Big| \frac{1}{n_k} \sum\limits_{i^\prime \in G_k} \overline{\varepsilon}_{i^\prime} \Big| \, \max_{i \in G_k} \Big| \frac{1}{\sqrt{p}} \sum\limits_{j=1}^p \varepsilon_{ij} \Big| = O_p \Big( \frac{p^\eta}{\sqrt{n_k p}} \Big). 
\end{equation}
Combining \eqref{R-iC1}, \eqref{R-iC2} and \eqref{R-iC3}, we arrive at the statement \eqref{R-iC} on the remainder $R_{i,C}^{[K_0]}$.

We finally show that  
\begin{equation}\label{R-iD}
\max_{i \in G_k} \big|R_{i,D}^{[K_0]}\big| = O_p\Big( \frac{p^{2\eta}}{\sqrt{p}} \Big). 
\end{equation}
For the proof, we write $R_{i,D}^{[K_0]} = \{ \widehat{\kappa}^{-1} \widehat{\sigma}^{-2} \} \{ R_{i,D,1}^{[K_0]} + R_{i,D,2}^{[K_0]} \}$, where
\begin{align*}
R_{i,D,1}^{[K_0]} & = \frac{1}{\sqrt{p}} \Big( \frac{1}{\sqrt{p}} \sum\limits_{j=1}^p \varepsilon_{ij} \Big)^2 \\
R_{i,D,2}^{[K_0]} & = \frac{1}{\sqrt{p}} \sum\limits_{j=1}^p \Big\{ \frac{1}{n_k} \sum\limits_{i^\prime \in G_k} \big( \varepsilon_{i^\prime j} - \overline{\varepsilon}_{i^\prime} \big) \Big\}^2.
\end{align*}
With the help of Lemma \ref{lemmaA1}, we obtain that 
\begin{equation}\label{R-iD1}
\max_{i \in G_k} \big|R_{i,D,1}^{[K_0]}\big| = O_p\Big( \frac{p^{2\eta}}{\sqrt{p}} \Big). 
\end{equation}
Moreover, straightforward calculations yield that 
\begin{equation}\label{R-iD2}
\max_{i \in G_k} \big|R_{i,D,2}^{[K_0]}\big| = O_p\Big( \frac{\sqrt{p}}{n_k} \Big). 
\end{equation}
\eqref{R-iD} now follows upon combining \eqref{R-iD1} and \eqref{R-iD2}. \qed 

\vspace{12pt}
\newpage

\noindent \textbf{Proof of Lemma \ref{lemmaA3}.} We make use of the following three results:
\begin{enumerate}[label=(R\arabic*),leftmargin=1.05cm]
\item \label{R1} Let $\{W_i: 1 \le i \le n\}$ be independent random variables with a standard normal distribution and define $a_n = 1/\sqrt{2 \log n}$ together with 
\[ b_n = \sqrt{2 \log n} - \frac{\log \log n + \log (4\pi)}{2 \sqrt{2 \log n}}. \]
Then for any $w \in \reals$, 
\[ \lim_{n \rightarrow \infty} \pr \Big( \max_{1 \le i \le n} W_i \le a_nw + b_n \Big) = \exp(-\exp(-w)). \]
In particular, for $w(\alpha \pm \varepsilon) = -\log(-\log(1 - \alpha \pm \varepsilon))$, we get 
\[ \lim_{n \rightarrow \infty} \pr \Big( \max_{1 \le i \le n} W_i \le a_nw(\alpha \pm \varepsilon) + b_n \Big) = 1 - \alpha \pm \varepsilon. \]
\end{enumerate}
The next result is known as Khintchine's Theorem. 
\begin{enumerate}[label=(R\arabic*),leftmargin=1.05cm]
\setcounter{enumi}{1}
\item \label{R2} Let $F_n$ be distribution functions and $G$ a non-degenerate distribution function. Moreover, let $\alpha_n > 0$ and $\beta_n \in \reals$ be such that 
\[ F_n(\alpha_n x + \beta_n) \rightarrow G(x) \]
for any continuity point $x$ of $G$. Then there are constants $\alpha_n^\prime > 0$ and $\beta_n^\prime \in \reals$ as well as a non-degenerate distribution function $G_*$ such that
\[ F_n(\alpha_n^\prime x + \beta_n^\prime) \rightarrow G_*(x) \]
at any continuity point $x$ of $G_*$ if and only if 
\[ \alpha_n^{-1} \alpha_n^\prime \rightarrow \alpha_*, \quad \frac{\beta_n^\prime - \beta_n}{\alpha_n} \rightarrow \beta_* \quad \text{and} \quad G_*(x) = G(\alpha_* x + \beta_*). \]
\end{enumerate}
The final result exploits strong approximation theory and is a direct consequence of the so-called KMT Theorems; see \cite{KMT1975, KMT1976}: 
\begin{enumerate}[label=(R\arabic*),leftmargin=1.05cm]
\setcounter{enumi}{2}
\item \label{R3} Write
\[ \Delta_i^{[K_0]} = \frac{1}{\sqrt{p}} \sum\limits_{j=1}^p X_{ij} \quad \text{with } X_{ij} = \Big\{ \frac{\varepsilon_{ij}^2}{\sigma^2} - 1 \Big\} \Big/ \kappa \]
and let $F$ denote the distribution function of $X_{ij}$. It is possible to construct i.i.d.\ random variables $\{ \widetilde{X}_{ij}: 1 \le i \le n, \, 1 \le j \le p \}$ with the distribution function $F$ and independent standard normal random variables $\{ Z_{ij}: 1 \le i \le n, \, 1 \le j \le p \}$ such that
\begin{align*}
\widetilde{\Delta}_i^{[K_0]} = \frac{1}{\sqrt{p}} \sum\limits_{j=1}^p \widetilde{X}_{ij} \quad \text{and} \quad \Delta_i^* = \frac{1}{\sqrt{p}} \sum\limits_{j=1}^p Z_{ij} 
\end{align*}
have the following property:
\[ \pr \Big( \big| \widetilde{\Delta}_i^{[K_0]} - \Delta_i^* \big| > C p^{\frac{1}{2 + \delta} - \frac{1}{2}} \Big) \le p^{1 - \frac{\theta/2}{2 + \delta}} \]
for some arbitrarily small but fixed $\delta > 0$ and some constant $C > 0$ that does not depend on $i$, $p$ and $n$. 
\end{enumerate}

\noindent We now proceed as follows: 
\begin{enumerate}[label=(\roman*),leftmargin=0.85cm]
\item We show that for any $w \in \reals$, 
\begin{equation}\label{lemmaA3-res1a}
\pr \Big( \max_{1 \le i \le n} \Delta_i^{[K_0]} \le a_n w + b_n \Big) \rightarrow \exp(-\exp(-w)).
\end{equation}
This in particular implies that 
\begin{equation}\label{lemmaA3-res1b}
\pr \Big( \max_{1 \le i \le n} \Delta_i^{[K_0]} \le w_n(\alpha \pm \varepsilon) \Big) \rightarrow 1 - \alpha \pm \varepsilon,
\end{equation}
where $w_n(\alpha \pm \varepsilon) = a_n w(\alpha \pm \varepsilon) + b_n$ with $a_n$, $b_n$ and $w(\alpha \pm \varepsilon)$ as defined in \ref{R1}. The proof of \eqref{lemmaA3-res1a} is postponed until the arguments for Lemma \ref{lemmaA3} are complete. 
\item The statement \eqref{lemmaA3-res1b} in particular holds in the special case that $\varepsilon_{ij} \sim \normal(0,\sigma^2)$. In this case, $q(\alpha)$ is the $(1-\alpha)$-quantile of $\max_{1 \le i \le n} \Delta_i^{[K_0]}$. Hence, we have 
\begin{align*}
 & \pr \Big( \max_{1 \le i \le n} \Delta_i^{[K_0]} \le w_n(\alpha - \varepsilon) \Big) \rightarrow 1 - \alpha - \varepsilon \\
 & \pr \Big( \max_{1 \le i \le n} \Delta_i^{[K_0]} \le q(\alpha) \Big) = 1 - \alpha \\
 & \pr \Big( \max_{1 \le i \le n} \Delta_i^{[K_0]} \le w_n(\alpha + \varepsilon) \Big) \rightarrow 1 - \alpha + \varepsilon,
\end{align*}
which implies that
\begin{equation}\label{lemmaA3-res2}
w_n(\alpha - \varepsilon) \le q(\alpha) \le w_n(\alpha + \varepsilon)
\end{equation}
for sufficiently large $n$. 
\item Since $p^{-\xi} / a_n = p^{-\xi} \sqrt{2 \log n} = o(1)$ by \ref{A3}, we can use \eqref{lemmaA3-res1a} together with \ref{R2} to obtain that 
\begin{equation}\label{lemmaA3-res3-equation1}
\pr \Big( \max_{1 \le i \le n} \Delta_i^{[K_0]} \le w_n(\alpha \pm \varepsilon) \pm p^{-\xi} \Big) \rightarrow 1 - \alpha \pm \varepsilon. 
\end{equation}
As $w_n(\alpha - \varepsilon) - p^{-\xi} \le q(\alpha) - p^{-\xi} \le  q(\alpha) + p^{-\xi} \le w_n(\alpha + \varepsilon) + p^{-\xi}$ for sufficiently large $n$, it holds that 
\begin{align*}
 & P_{\alpha,\varepsilon}^{\ll} := \pr \Big( \max_{1 \le i \le n} \Delta_i^{[K_0]} \le w_n(\alpha - \varepsilon) - p^{-\xi} \Big) \\
 & \le P_{\alpha}^{\ll} = \pr \Big( \max_{1 \le i \le n} \Delta_i^{[K_0]} \le q(\alpha) - p^{-\xi} \Big) \\
 & \le P_{\alpha}^{\gg} = \pr \Big( \max_{1 \le i \le n} \Delta_i^{[K_0]} \le q(\alpha) + p^{-\xi} \Big) \\
 & \le P_{\alpha,\varepsilon}^{\gg} := \pr \Big( \max_{1 \le i \le n} \Delta_i^{[K_0]} \le w_n(\alpha + \varepsilon) + p^{-\xi} \Big) 
\end{align*}
for large $n$. Moreover, since $P_{\alpha,\varepsilon}^{\ll} \rightarrow 1 - \alpha - \varepsilon$ and $P_{\alpha,\varepsilon}^{\gg} \rightarrow 1 - \alpha + \varepsilon$ for any fixed $\varepsilon > 0$ by \eqref{lemmaA3-res3-equation1}, we can conclude that $P_{\alpha}^{\ll} = (1 - \alpha) + o(1)$ and $P_{\alpha}^{\gg} = (1 - \alpha) + o(1)$, which is the statement of Lemma \ref{lemmaA3}. 
\end{enumerate}
It remains to prove \eqref{lemmaA3-res1a}: Using the notation from \ref{R3} and the shorthand $w_n = a_n w + b_n$, we can write 
\begin{equation}\label{lemmaA3-res1a-1}
\pr \Big( \max_{1 \le i \le n} \Delta_i^{[K_0]} \le w_n \Big) = \pr \Big( \max_{1 \le i \le n} \widetilde{\Delta}_i^{[K_0]} \le w_n \Big) = \prod\limits_{i=1}^n \pi_i
\end{equation} 
with
\[ \pi_i = \pr \Big( \widetilde{\Delta}_i^{[K_0]} \le w_n \Big). \]
The probabilities $\pi_i$ can be decomposed into two parts as follows: 
\[ \pi_i = \pr \Big( \Delta_i^* \le w_n + \big\{\Delta_i^* - \widetilde{\Delta}_i^{[K_0]}\big\} \Big) = \pi_i^{\le} + \pi_i^{>}, \]
where 
\begin{align*}
\pi_i^{\le} & = \pr \Big( \Delta_i^* \le w_n + \big\{\Delta_i^* - \widetilde{\Delta}_i^{[K_0]}\big\}, \big| \Delta_i^* - \widetilde{\Delta}_i^{[K_0]} \big| \le C p^{\frac{1}{2+\delta} - \frac{1}{2}} \Big) \\
\pi_i^{>} & = \pr \Big( \Delta_i^* \le w_n + \big\{\Delta_i^* - \widetilde{\Delta}_i^{[K_0]}\big\}, \big| \Delta_i^* - \widetilde{\Delta}_i^{[K_0]} \big| > C p^{\frac{1}{2+\delta} - \frac{1}{2}} \Big). 
\end{align*}
With the help of \ref{R3} and the assumption that $n \ll p^{(\theta/4)-1}$, we can show that  
\begin{equation}\label{lemmaA3-res1a-2}
\prod\limits_{i=1}^n \pi_i = \prod\limits_{i=1}^n \pi_i^{\le} + R_n, 
\end{equation}
where $R_n$ is a non-negative remainder term with 
\[ R_n \le \sum\limits_{i=1}^n \binom{n}{i} \Big(p^{1 - \frac{\theta/2}{2+\delta}} \Big)^i = o(1). \]
Moreover, the probabilities $\pi_i^{\le}$ can be bounded by 
\[ \pi_i^{\le} 
\begin{cases}
\le \pr \Big( \Delta_i^* \le w_n + C p^{\frac{1}{2+\delta} - \frac{1}{2}} \Big) \\
\ge \pr \Big( \Delta_i^* \le w_n - C p^{\frac{1}{2+\delta} - \frac{1}{2}} \Big) - p^{1 - \frac{\theta/2}{2+\delta}}, 
\end{cases}
\]
the second line making use of \ref{R3}. From this, we obtain that 
\begin{equation}\label{lemmaA3-res1a-3}
\prod\limits_{i=1}^n \pi_i^{\le} 
\begin{cases}
\le \overline{\Pi}_n \\
\ge \underline{\Pi}_n + o(1),
\end{cases}
\end{equation}
where
\begin{align*}
\overline{\Pi}_n & = \pr \Big( \max_{1 \le i \le n} \Delta_i^* \le w_n + C p^{\frac{1}{2+\delta} - \frac{1}{2}} \Big) \\
\underline{\Pi}_n & = \pr \Big( \max_{1 \le i \le n} \Delta_i^* \le w_n - C p^{\frac{1}{2+\delta} - \frac{1}{2}} \Big).
\end{align*}
By combining \eqref{lemmaA3-res1a-1}--\eqref{lemmaA3-res1a-3}, we arrive at the intermediate result that
\begin{equation}\label{lemmaA3-res1a-ir}
\underline{\Pi}_n  + o(1) \le \pr \Big( \max_{1 \le i \le n} \Delta_i^{[K_0]} \le w_n \Big) \le \overline{\Pi}_n + o(1). 
\end{equation}
Since $p^{\frac{1}{2+\delta} - \frac{1}{2}} / a_n = p^{\frac{1}{2+\delta} - \frac{1}{2}} \sqrt{2 \log n} = o(1)$, we can use \ref{R1} together with \ref{R2} to show that 
\begin{equation}\label{lemmaA3-res1a-Piconv}
\overline{\Pi}_n \rightarrow \exp(-\exp(-w)) \qquad \text{and} \qquad \underline{\Pi}_n \rightarrow \exp(-\exp(-w)).
\end{equation} 
Plugging \eqref{lemmaA3-res1a-Piconv} into \eqref{lemmaA3-res1a-ir} immediately yields that 
\[ \pr \Big( \max_{1 \le i \le n} \Delta_i^{[K_0]} \le w_n \Big) \rightarrow \exp(-\exp(-w)), \]
which completes the proof.
\qed

\vspace{12pt}

\noindent \textbf{Proof of Lemma \ref{lemmaA4}.} We use the notation $n_S = \# S$ along with $\overline{\varepsilon}_i = p^{-1} \sum\nolimits_{j=1}^p \varepsilon_{ij}$, $\overline{\mu}_i = p^{-1} \sum\nolimits_{j=1}^p \mu_{ij}$ and $\overline{d}_i = p^{-1} \sum\nolimits_{j=1}^p d_{ij}$. For any $i \in S$ and $S \in \mathscr{C}$, we can write 
\[ \widehat{\Delta}_i^{[K]} = \frac{1}{\kappa \sigma^2 \sqrt{p}} \sum\limits_{j=1}^p d_{ij}^2 + R_{i,A}^{[K]} + R_{i,B}^{[K]} + R_{i,C}^{[K]} + R_{i,D}^{[K]} - R_{i,E}^{[K]} + R_{i,F}^{[K]} + R_{i,G}^{[K]}, \]
where
\begin{align*}
R_{i,A}^{[K]} & = \Big( \frac{1}{\widehat{\kappa} \widehat{\sigma}^2} - \frac{1}{\kappa \sigma^2} \big) \frac{1}{\sqrt{p}} \sum\limits_{j=1}^p d_{ij}^2 \\
R_{i,B}^{[K]} & = \frac{1}{\sqrt{p}} \sum\limits_{j=1}^p \Big\{ \frac{\varepsilon_{ij}^2}{\sigma^2} - 1 \Big\} \Big/ \kappa \\
R_{i,C}^{[K]} & = \Big( \frac{1}{\widehat{\kappa}} - \frac{1}{\kappa} \Big) \, \frac{1}{\sqrt{p}} \sum\limits_{j=1}^p \Big\{ \frac{\varepsilon_{ij}^2}{\sigma^2} - 1 \Big\} \\
R_{i,D}^{[K]} & = \frac{1}{\widehat{\kappa}} \Big( \frac{1}{\widehat{\sigma}^2} - \frac{1}{\sigma^2} \Big) \, \frac{1}{\sqrt{p}} \sum\limits_{j=1}^p \varepsilon_{ij}^2 \\
R_{i,E}^{[K]} & = \Big(\frac{2}{\widehat{\kappa} \widehat{\sigma}^2}\Big) \, \frac{1}{\sqrt{p}} \sum\limits_{j=1}^p  \varepsilon_{ij} \Big\{ \overline{\varepsilon}_i + \frac{1}{n_S} \sum\limits_{i^\prime \in S} \big( \varepsilon_{i^\prime j} - \overline{\varepsilon}_{i^\prime} \big) \Big\} \\ 
R_{i,F}^{[K]} & = \Big(\frac{1}{\widehat{\kappa} \widehat{\sigma}^2}\Big) \, \frac{1}{\sqrt{p}} \sum\limits_{j=1}^p \Big\{ \overline{\varepsilon}_i + \frac{1}{n_S} \sum\limits_{i^\prime \in S} \big( \varepsilon_{i^\prime j} - \overline{\varepsilon}_{i^\prime} \big) \Big\}^2 \\
R_{i,G}^{[K]} & = \Big(\frac{2}{\widehat{\kappa} \widehat{\sigma}^2}\Big) \, \frac{1}{\sqrt{p}} \sum\limits_{j=1}^p  \Big\{ \varepsilon_{ij} - \overline{\varepsilon}_i - \frac{1}{n_S} \sum\limits_{i^\prime \in S} \big( \varepsilon_{i^\prime j} - \overline{\varepsilon}_{i^\prime} \big) \Big\} d_{ij}.  
\end{align*}
We now show that $\max_{S \in \mathscr{C}} \max_{i \in S} |R_{i,\ell}^{[K]}| = o_p(p^{1/2-\xi})$ for $\ell = A,\ldots,G$. This immediately yields the statement of Lemma \ref{lemmaA4}. Throughout the proof, $\eta > 0$ denotes a sufficiently small constant that results from applying Lemma \ref{lemmaA1}.

With the help of Lemma \ref{lemmaA1} and our assumptions on $\widehat{\sigma}^2$ and $\widehat{\kappa}$, it is straightforward to see that $\max_{S \in \mathscr{C}} \max_{i \in S} |R_{i,\ell}^{[K]}| \le \max_{1 \le i \le n} |R_{i,\ell}^{[K]}| = o_p(p^{1/2 - \xi})$ for $\ell = A,B,C,D$ with some sufficiently small $\xi > 0$. We next show that 
\begin{equation}\label{R-iD-K}
\max_{S \in \mathscr{C}} \max_{i \in S} \big| R_{i,E}^{[K]} \big| = O_p\big(p^{\eta}\big). 
\end{equation}
To do so, we write $R_{i,E}^{[K]} = \{ 2 \widehat{\kappa}^{-1} \widehat{\sigma}^{-2} \} \{ R_{i,E,1}^{[K]} + R_{i,E,2}^{[K]} - R_{i,E,3}^{[K]} \}$, where 
\begin{align*}
R_{i,E,1}^{[K]} & = \frac{1}{\sqrt{p}} \sum\limits_{j=1}^p \varepsilon_{ij} \overline{\varepsilon}_i \\
R_{i,E,2}^{[K]} & = \frac{1}{\sqrt{p}} \sum\limits_{j=1}^p \varepsilon_{ij} \Big( \frac{1}{n_S} \sum\limits_{i^\prime \in S} \varepsilon_{i^\prime j} \Big) \\
R_{i,E,3}^{[K]} & = \Big( \frac{1}{\sqrt{p}} \sum\limits_{j=1}^p \varepsilon_{ij} \Big) \Big( \frac{1}{n_S} \sum\limits_{i^\prime \in S} \overline{\varepsilon}_{i^\prime} \Big).
\end{align*}
Lemma \ref{lemmaA1} yields that $\max_{S \in \mathscr{C}} \max_{i \in S} |R_{i,E,1}^{[K]}| \le \max_{1 \le i \le n} |R_{i,E,1}^{[K]}| = O_p(p^{2\eta}/\sqrt{p})$. Moreover, it holds that
\[ \max_{S \in \mathscr{C}} \max_{i \in S} \big|R_{i,E,2}^{[K]}\big| = O_p\big(p^{\eta}\big), \]
since
\[ R_{i,E,2}^{[K]} = \frac{1}{n_S \sqrt{p}} \sum\limits_{j=1}^p \big\{ \varepsilon_{ij}^2 - \sigma^2 \big\} +  \sigma^2 \frac{\sqrt{p}}{n_S} + \frac{1}{n_S} \sum\limits_{\substack{i^\prime \in S \\ i^\prime \ne i}} \frac{1}{\sqrt{p}} \sum\limits_{j=1}^p \varepsilon_{ij} \varepsilon_{i^\prime j} \]
and
\begin{align*} 
 & \max_{S \in \mathscr{C}} \max_{i \in S} \Big| \frac{1}{n_S \sqrt{p}} \sum\limits_{j=1}^p \big\{ \varepsilon_{ij}^2 - \sigma^2 \big\} \Big| \le \frac{1}{\underline{n}} \max_{1 \le i \le n} \Big| \frac{1}{\sqrt{p}} \sum\limits_{j=1}^p \big\{ \varepsilon_{ij}^2 - \sigma^2 \big\} \Big| = O_p\Big( \frac{p^{\eta}}{\underline{n}} \Big) \\
 & \max_{S \in \mathscr{C}} \max_{i \in S} \Big| \frac{1}{n_S} \sum\limits_{\substack{i^\prime \in S \\ i^\prime \ne i}} \frac{1}{\sqrt{p}} \sum\limits_{j=1}^p \varepsilon_{ij} \varepsilon_{i^\prime j} \Big| \le \max_{1 \le i < i^\prime \le n} \Big| \frac{1}{\sqrt{p}} \sum\limits_{j=1}^p \varepsilon_{ij} \varepsilon_{i^\prime j} \Big| = O_p\big(p^{\eta}\big), 
\end{align*}
which follows upon applying Lemma \ref{lemmaA1}. Finally, 
\[ \max_{S \in \mathscr{C}} \max_{i \in S} \big| R_{i,E,3}^{[K]} \big| \le \frac{1}{\sqrt{p}} \Big\{ \max_{1 \le i \le n} \Big| \frac{1}{\sqrt{p}} \sum\limits_{j=1}^p \varepsilon_{ij} \Big| \Big\}^2 = O_p\Big( \frac{p^{2\eta}}{\sqrt{p}} \Big), \]
which can again be seen by applying Lemma \ref{lemmaA1}. Putting everything together, we arrive at \eqref{R-iD-K}. Similar arguments show that 
\begin{equation}\label{R-iE-K}
\max_{S \in \mathscr{C}} \max_{i \in S} \big| R_{i,F}^{[K]} \big| = O_p\big(p^{\eta}\big) 
\end{equation}
as well.

To analyze the term $R_{i,G}^{[K]}$, we denote the signal vector of the group $G_k$ by $\boldsymbol{m}_k = (m_{1,k},\ldots,m_{p,k})^\top$ and write
\[ \frac{1}{n_S} \sum\limits_{i \in S} \mu_{ij} = \sum\limits_{k=1}^{K_0} \lambda_{S,k} m_{j,k} \]
with $\lambda_{S,k} = \# (S \cap G_k) / n_S$. With this notation, we get
\[ R_{i,G}^{[K]} = \{ 2 \widehat{\kappa}^{-1} \widehat{\sigma}^{-2} \} \{ R_{i,G,1}^{[K]} - R_{i,G,2}^{[K]} - R_{i,G,3}^{[K]} - R_{i,G,4}^{[K]} + R_{i,G,5}^{[K]} \}, \]
where
\begin{align*}
R_{i,G,1}^{[K]} & = \frac{1}{\sqrt{p}} \sum\limits_{j=1}^p \mu_{ij} \varepsilon_{ij} \\
R_{i,G,2}^{[K]} & = \sum\limits_{k=1}^{K_0} \lambda_{S,k} \frac{1}{\sqrt{p}} \sum\limits_{j=1}^p m_{j,k} \varepsilon_{ij} \\
R_{i,G,3}^{[K]} & = \frac{1}{\sqrt{p}} \sum\limits_{j=1}^p \overline{\varepsilon}_i d_{ij} \\
R_{i,G,4}^{[K]} & = \frac{1}{n_S} \sum\limits_{i^\prime \in S} \frac{1}{\sqrt{p}} \sum\limits_{j=1}^p \mu_{ij} \varepsilon_{i^\prime j} \\
R_{i,G,5}^{[K]} & = \frac{1}{n_S} \sum\limits_{i^\prime \in S} \sum\limits_{k=1}^{K_0} \lambda_{S,k} \frac{1}{\sqrt{p}} \sum\limits_{j=1}^p m_{j,k} \varepsilon_{i^\prime j}. 
\end{align*}
With the help of Lemma \ref{lemmaA1}, it can be shown that $\max_{S \in \mathscr{C}} \max_{i \in S} \big| R_{i,G,\ell}^{[K]} \big| = O_p(p^\eta)$ for $\ell = 1,\ldots,5$. For example, it holds that
\[ \max_{S \in \mathscr{C}} \max_{i \in S} \big| R_{i,G,4}^{[K]} \big| \le \max_{1 \le i < i^\prime \le n} \Big| \frac{1}{\sqrt{p}} \sum\limits_{j=1}^p \mu_{ij} \varepsilon_{i^\prime j} \Big| = O_p(p^\eta). \]
As a result, we obtain that  
\begin{equation}\label{R-iF-K}
\max_{S \in \mathscr{C}} \max_{i \in S} \big| R_{i,G}^{[K]} \big| = O_p\big(p^{\eta}\big). 
\end{equation} 
This completes the proof. \qed

\vspace{12pt}

\noindent \textbf{Proof of Lemma \ref{lemmaA5}.} Let $S \in \mathscr{C}$. In particular, suppose that $S \cap G_{k_1} \ne \emptyset$ and $S \cap G_{k_2} \ne \emptyset$ for some $k_1 \ne k_2$. We show the following claim: there exists some $i \in S$ such that 
\begin{equation}\label{lemmaA5-claim}
\frac{1}{\sqrt{p}} \sum\limits_{j=1}^p d_{ij}^2 \ge c \sqrt{p}, 
\end{equation}
where $c = (\sqrt{\delta_0}/2)^2$ with $\delta_0$ defined in assumption \ref{A2}. From this, the statement of Lemma \ref{lemmaA5} immediately follows.

For the proof of \eqref{lemmaA5-claim}, we denote the Euclidean distance between vectors $v = (v_1,\ldots,v_p)^\top$ and $w = (w_1,\ldots,w_p)^\top$ by $d(v,w) = ( \sum\nolimits_{j=1}^p |v_j - w_j|^2 )^{1/2}$. Moreover, as in Lemma \ref{lemmaA4}, we use the notation
\[ \frac{1}{n_S} \sum\limits_{i \in S} \mu_{ij} = \sum\limits_{k=1}^{K_0} \lambda_{S,k} m_{j,k}, \]
where $n_S = \# S$, $\lambda_{S,k} = \# (S \cap G_k) / n_S$ and $\boldsymbol{m}_k = (m_{1,k},\ldots,m_{p,k})^\top$ is the signal vector of the class $G_k$.

Take any $i \in S \cap G_{k_1}$. If 
\[ d\Big( \boldsymbol{\mu}_i, \sum\limits_{k=1}^{K_0} \lambda_{S,k} \boldsymbol{m}_k \Big) = d\Big( \boldsymbol{m}_{k_1}, \sum\limits_{k=1}^{K_0} \lambda_{S,k} \boldsymbol{m}_k \Big) \ge \frac{\sqrt{\delta_0 p}}{2}, \]
the proof is finished, as \eqref{lemmaA5-claim} is satisfied for $i$. Next consider the case that 
\[ d\Big( \boldsymbol{m}_{k_1}, \sum\limits_{k=1}^{K_0} \lambda_{S,k} \boldsymbol{m}_k \Big) < \frac{\sqrt{\delta_0 p}}{2}. \]
By assumption \ref{A2}, it holds that $d(\boldsymbol{m}_k,\boldsymbol{m}_{k^\prime}) \ge \sqrt{\delta_0 p}$ for $k \ne k^\prime$. Hence, by the triangle inequality, 
\begin{align*}
\sqrt{\delta_0 p} 
 & \le d\big(\boldsymbol{m}_{k_1},\boldsymbol{m}_{k_2}\big) \\
 & \le d\Big( \boldsymbol{m}_{k_1}, \sum\limits_{k=1}^{K_0} \lambda_{S,k} \boldsymbol{m}_k \Big) + d\Big( \sum\limits_{k=1}^{K_0} \lambda_{S,k} \boldsymbol{m}_k, \boldsymbol{m}_{k_2} \Big) \\
 & < \frac{\sqrt{\delta_0 p}}{2} + d\Big( \sum\limits_{k=1}^{K_0} \lambda_{S,k} \boldsymbol{m}_k, \boldsymbol{m}_{k_2} \Big), 
\end{align*}
implying that 
\[ d\Big( \sum\limits_{k=1}^{K_0} \lambda_{S,k} \boldsymbol{m}_k, \boldsymbol{m}_{k_2} \Big) > \frac{\sqrt{\delta_0 p}}{2}. \]
This shows that the claim \eqref{lemmaA5-claim} is fulfilled for any $i^\prime \in S \cap G_{k_2}$. \qed

\section*{Proof of Theorem \ref{theo2}}

By Theorem \ref{theo1},
\begin{align*} 
 & \pr \big( \widehat{K}_0 > K_0 \big) \\
 & = \pr \Big( \widehat{\mathcal{H}}^{[K]} > q(\alpha) \text{ for all } K \le K_0 \Big) \\
 & = \pr \Big( \widehat{\mathcal{H}}^{[K_0]} > q(\alpha) \Big) - \pr \Big( \widehat{\mathcal{H}}^{[K_0]} > q(\alpha), \widehat{\mathcal{H}}^{[K]} \le q(\alpha) \text{ for some } K < K_0 \Big) \\
 & = \pr \Big( \widehat{\mathcal{H}}^{[K_0]} > q(\alpha) \Big) + o(1) \\
 & = \alpha + o(1) 
\end{align*}
and
\begin{align*}
\pr \big( \widehat{K}_0 < K_0 \big) 
 & = \pr \Big( \widehat{\mathcal{H}}^{[K]} \le q(\alpha) \text{ for some } K < K_0 \Big) \\
 & \le \sum\limits_{K=1}^{K_0-1} \pr \Big( \widehat{\mathcal{H}}^{[K]} \le q(\alpha) \Big) \\
 & = o(1). 
\end{align*}
Moreover, 
\begin{align*} 
 & \pr \Big( \big\{ \widehat{G}_k: 1 \le k \le \widehat{K}_0 \big\} \ne \big\{ G_k: 1 \le k \le K_0 \big\} \Big) \\
 & = \pr \Big( \big\{ \widehat{G}_k: 1 \le k \le \widehat{K}_0 \big\} \ne \big\{ G_k: 1 \le k \le K_0 \big\}, \widehat{K}_0 = K_0 \Big) \\
 & \quad + \pr \Big( \big\{ \widehat{G}_k: 1 \le k \le \widehat{K}_0 \big\} \ne \big\{ G_k: 1 \le k \le K_0 \big\}, \widehat{K}_0 \ne K_0 \Big) \\
 & = \alpha + o(1),
\end{align*}
since 
\begin{align*}
 & \pr \Big( \big\{ \widehat{G}_k: 1 \le k \le \widehat{K}_0 \big\} \ne \big\{ G_k: 1 \le k \le K_0 \big\}, \widehat{K}_0 = K_0 \Big) \\
 & = \pr \Big( \big\{ \widehat{G}_k^{[K_0]}: 1 \le k \le K_0 \big\} \ne \big\{ G_k: 1 \le k \le K_0 \big\}, \widehat{K}_0 = K_0 \Big) \\
 & \le \pr \Big( \big\{ \widehat{G}_k^{[K_0]}: 1 \le k \le K_0 \big\} \ne \big\{ G_k: 1 \le k \le K_0 \big\} \Big) \\
 & = o(1)
\end{align*}
by the consistency property \eqref{prop-clusters} and
\begin{align*} 
 & \pr \Big( \big\{ \widehat{G}_k: 1 \le k \le \widehat{K}_0 \big\} \ne \big\{ G_k: 1 \le k \le K_0 \big\}, \widehat{K}_0 \ne K_0 \Big) \\
 & = \pr \big( \widehat{K}_0 \ne K_0 \big) = \alpha + o(1).
\end{align*}
\vspace{-1.6cm}

\qed

\section*{Proof of Theorem \ref{theo3}}

With the help of Lemma \ref{lemmaA1}, we can show that 
\begin{align}
\widehat{\rho}(i,i^\prime) 
 & = 2 \sigma^2 + \frac{1}{p}\sum\limits_{j=1}^p \big( \mu_{ij} - \mu_{i^\prime j} \big)^2 + o_p(1) \label{theo3-eq1} 
\end{align}
uniformly over $i$ and $i^\prime$. This together with \ref{A2} allows us to prove the following claim:
\begin{equation}\label{theo3-eq2}
\begin{minipage}[c][1.5cm][c]{12.9cm}
With probability tending to $1$, the indices $i_1,\ldots,i_K$ belong to $K$ different classes in the case that $K \le K_0$ and to $K_0$ different classes in the case that $K > K_0$.
\end{minipage}
\end{equation}
Now let $K = K_0$. With the help of \eqref{theo3-eq1} and \eqref{theo3-eq2}, the starting values $\mathscr{C}_1^{[K_0]},\ldots$ $\ldots,\mathscr{C}_{K_0}^{[K_0]}$ can be shown to have the property that
\begin{equation}\label{theo3-eq3}
\pr \Big( \big\{ \mathscr{C}_k^{[K_0]}: 1 \le k \le K_0 \big\} = \big\{ G_k: 1 \le k \le K_0 \big\} \Big) \rightarrow 1. 
\end{equation}
Together with Lemma \ref{lemmaA1}, \eqref{theo3-eq3} yields that 
\[ \widehat{\rho}_k^{(1)}(i) = \sigma^2 + \frac{1}{p}\sum\limits_{j=1}^p \big( \mu_{ij} - m_{j,k} \big)^2 + o_p(1) \]
uniformly over $i$ and $k$. Combined with \ref{A2}, this in turn implies that the $k$-means algorithm converges already after the first iteration step with probability tending to $1$ and $\widehat{G}_k^{[K_0]}$ are consistent estimators of the classes $G_k$ in the sense of \eqref{prop-clusters}. \qed

\section*{Proof of (\ref{rate-sigma-RSS})}

Suppose that \ref{A1}--\ref{A3} along with \eqref{ass-sigma-RSS} are satisfied. As already noted in Section \ref{subsec-est-sigma}, the $k$-means estimators $\{ \widehat{G}_k^A: 1 \le k \le K_{\max} \}$ can be shown to satisfy \eqref{add-prop-clusters}, that is,
\begin{equation}\label{prop-clusters-A}
\pr \Big( \widehat{G}_k^A \subseteq G_{k^\prime} \text{ for some } 1 \le k^\prime \le K_0 \Big) \rightarrow 1 
\end{equation}
for any $k = 1,\ldots,K_{\max}$. This can be proven by very similar arguments as the consistency property \eqref{prop-clusters}. We thus omit the details. Let $E^A$ be the event that 
\begin{equation*}
\widehat{G}_k^A \subseteq G_{k^\prime} \text{ for some } 1 \le k^\prime \le K_0 
\end{equation*}
holds for all clusters $\widehat{G}_k^A$ with $k = 1,\ldots,K_{\max}$. $E^A$ can be regarded as the event that the partition $\{ \widehat{G}_k^A: 1 \le k \le K_{\max} \}$ is a refinement of the class structure $\{G_k: 1 \le k \le K_0\}$. By \eqref{prop-clusters-A}, the event $E^A$ occurs with probability tending to $1$.

Now consider the estimator  
\[ \widehat{\sigma}^2_{\text{RSS}} = \frac{1}{n \lfloor p/2 \rfloor} \sum\limits_{k=1}^{K_{\max}} \sum\limits_{i \in \widehat{G}_k^A} \Big\| \widehat{\boldsymbol{Y}}_i^B - \frac{1}{\#\widehat{G}_k^A} \sum\limits_{i^\prime \in \widehat{G}_k^A} \widehat{\boldsymbol{Y}}_{i^\prime}^B \Big\|^2. \]
Since the random variables $\widehat{\boldsymbol{Y}}_i^B$ are independent of the estimators $\widehat{G}_k^A$, it is not difficult to verify the following: for any $\delta > 0$, there exists a constant $C_{\delta} > 0$ (that does not depend on $\{ \widehat{G}_k^A: 1 \le k \le K_{\max} \}$) such that on the event $E^A$, 
\[ \pr \Big( \big| \widehat{\sigma}^2_{\text{RSS}} - \sigma^2 \big| \ge \frac{C_{\delta}}{p} \, \Big| \, \big\{ \widehat{G}_k^A: 1 \le k \le K_{\max} \big\} \Big) \le \delta. \] 
From this, the first statement of \eqref{rate-sigma-RSS} easily follows. The second statement can be obtained by similar arguments. \qed

\end{document}